\documentclass[a4paper, 11pt]{article}

\usepackage{latexsym,amsfonts, amsmath, amssymb, amsthm}
\usepackage[utf8]{inputenc}
\usepackage[final=true,protrusion=true,expansion=true]{microtype}
\usepackage{enumitem}
\usepackage{color}
\usepackage{mathtools}
\usepackage{booktabs}
\usepackage{multirow}
\usepackage{subcaption}
\usepackage[font=small,labelfont=bf]{caption}
\usepackage[affil-it]{authblk}
\usepackage{geometry}
\usepackage{hyperref}
\hypersetup{
	final,
    colorlinks=true,
    linkcolor=blue,
    filecolor=magenta,
    urlcolor=cyan,
}
\usepackage{cleveref}

\theoremstyle{plain}
\newtheorem{theorem}{Theorem}[section]

\newtheorem{lemma}[theorem]{Lemma}
\newtheorem{corollary}[theorem]{Corollary}
\newtheorem{assumption}[theorem]{Assumption}

\theoremstyle{plain}

\theoremstyle{definition}
\newtheorem{remark}[theorem]{Remark}

\numberwithin{equation}{section}

\providecommand*{\abs}[1]{\left|{#1}\right|}
\providecommand*{\absbig}[1]{\big|{#1}\big|}% Double bar norm
\providecommand*{\absnormal}[1]{|{#1}|}% Double bar norm
\providecommand*{\N}[1]{\left\|{#1}\right\|} % Double bar norm
\providecommand*{\Nbig}[1]{\big\|{#1}\big\|} % Double bar norm
\providecommand{\CL}{{\cal L}}

\providecommand{\bbR}{{\mathbb{R}}}
\newcommand{\x}{\mathbf{x}}

\newcommand{\dummy}{\mathord{\color{black!33}\bullet}}

\providecommand*{\vol}[1]{\operatorname{vol}\!\left({#1}\right)}
\providecommand{\ker}{\operatorname{ker}}

\providecommand*{\ker}[1]{{\rm ker}({#1})}                  % Kernel
\providecommand{\Id}{\mathrm{Id}}                     % Identity operator

\providecommand{\CC}{{\cal C}}
\providecommand{\CD}{{\cal D}}

\providecommand{\CH}{{\cal H}}

\providecommand{\CJ}{{\cal J}}

\providecommand{\CL}{{\cal L}}
\providecommand{\CM}{{\cal M}}

\providecommand{\CP}{{\cal P}}

\providecommand{\bbE}{\mathbb{E}}

\providecommand{\bbN}{\mathbb{N}}

\providecommand{\bbP}{\mathbb{P}}

\providecommand{\bbR}{\mathbb{R}}

\providecommand*{\abs}[1]{\left|{#1}\right|}
\providecommand*{\absbig}[1]{\big|{#1}\big|}% Double bar norm
\providecommand*{\absnormal}[1]{|{#1}|}% Double bar norm
\providecommand*{\N}[1]{\left\|{#1}\right\|} % Double bar norm
\providecommand*{\Nbig}[1]{\big\|{#1}\big\|} % Double bar norm
\providecommand*{\abs}[1]{\left|{#1}\right|} % absolute value, same as \SN
\usepackage{bbm}
\usepackage{arydshln}
\usepackage{graphicx}
\usepackage{xcolor}
\usepackage{ifdraft}
\usepackage{scalerel}
\usepackage{mathrsfs}  
\usepackage{etoc}
\usepackage{comment}
\usepackage{natbib}
\DeclareFontEncoding{LS1}{}{}
\DeclareFontSubstitution{LS1}{stix2}{m}{n}
\DeclareMathAlphabet{\mathscr}{LS1}{stix2scr}{m}{n}
\providecommand{\CB}{\mathscr{B}}

\newif\ifrevised

\revisedfalse
\makeatother

\title{\usefont{OT1}{bch}{b}{n}
	\huge
    Global Convergence of Second-Order Methods\\for Physics-Informed Neural Networks\\in the Overparameterized Limit%
    \footnote{This article is part of the project ``DMS-EPSRC: Asymptotic Analysis of Online Training Algorithms in Machine Learning: Recurrent, Graphical, and Deep Neural Networks'' (NSF DMS-2311500).}
}
\date{}
\author[1]{Konstantin Riedl\thanks{Email: \texttt{Konstantin.Riedl@maths.ox.ac.uk}}}
\author[1]{Justin Sirignano\thanks{Email: \texttt{Justin.Sirignano@maths.ox.ac.uk}}}
\author[2]{Konstantinos Spiliopoulos\thanks{Email: \texttt{kspiliop@bu.edu}}}

\affil[1]{University of Oxford, Mathematical Institute}
\affil[2]{Boston University, Department of Mathematics \& Statistics}

\begin{document}
\maketitle

\begin{abstract}
\noindent
    Solving high-dimensional partial differential equations (PDEs) constitutes a core task in many engineering and scientific disciplines. Physics-informed neural networks and the deep Galerkin method effectively leverage machine learning models to approximate solutions of known PDEs by minimizing residuals corresponding to the PDE operator alongside external constraints. Although second-order optimizers
    are becoming increasingly popular for training neural networks (NNs) to solve PDEs and empirically outperform their first-order counterparts in both speed and accuracy, the convergence theory for training such models is essentially confined to first-order methods.
    In this paper, we close this gap and develop a convergence analysis for the regularized Newton and Gauss--Newton methods. As the number of hidden units tends to infinity, the training dynamics converge in probability, with explicit rates, to the solution of a deterministic limit equation involving a Newton neural tangent kernel operator associated with the PDE. Although the limiting kernel operator admits no spectral gap, we prove that, in the overparameterized limit, the PDE residual vanishes as training proceeds, and, under mild assumptions, the NN converges to the PDE solution. We further show that, provided a sufficiently small regularization parameter, the projection of the residual onto any fixed number of leading eigenfunctions decays exponentially fast, demonstrating that this hyperparameter acts as a spectral filter: it lifts the convergence rates of high-frequency components toward those of low-frequency ones, thereby mitigating the spectral bias inherent in gradient descent.
    Mathematical challenges arise from the implicit parameter update involving a linear system whose dimension increases with the NN width and, for Newton's method, a potentially indefinite Hessian, as well as from the PDE operator entering the kernel.
\end{abstract}



%65K10 Numerical optimization and variational techniques
%68T07 Artificial neural networks and deep learning
%65N12 Stability and convergence for numerical PDE boundary-value problems
%49M15 Newton-type methods
%90C26 Nonconvex programming, global optimization
%65J15 Numerical methods for nonlinear operator equations
%65N15 Error estimates for numerical PDE boundary-value problems
%68T05 Learning and adaptive systems in artificial intelligence
%90C30 Nonlinear programming
%35J25 Boundary value problems for second-order elliptic equations
%%%%%%%%%%%%%%%%%%%%%%%%%%%%%%%%%%%%%%%%%%%%%%%%%%%%%
%%%%%%%%%%%%%%%%%%%%%%%%%%%%%%%%%%%%%%%%%%%%%%%%%%%%%

%%%%%%%%%%%%%%%%%%%%%%%%%%%%%%%%%%%%%%%%%%%%%%%%%%%%%
%%%%%%%%%%%%%%%%%%%%%%%%%%%%%%%%%%%%%%%%%%%%%%%%%%%%%

\section{Introduction}
\label{sec:introduction}

Leveraging the property of neural networks (NNs) as universal approximators~\cite{hornik1989multilayer,hornik1991approximation}, methods such as physics-informed neural networks (PINNs)~\citep{raissi2019physics,karniadakis2021physics,lu2021deepxde,cuomo2022scientific,wang2023expert,jiang2023global,zhao2025convergence}, the deep Galerkin method \citep{sirignano2018dgm}, the deep Ritz method \citep{yu2018deep}, and neural Q-learning~\citep{cohen2023neural} employ NNs to approximate solutions of known partial differential equations (PDEs).
To achieve this, the NN parameters are optimized to minimize residuals corresponding to the differential operator, alongside initial and boundary conditions.
Typically, these models are trained using first-order optimization methods, such as gradient descent, RMSProp, or Adam, which rely on automatic differentiation to compute the gradient of the loss function.
Several recent practical works \cite{muller2023achieving,dangel2024kronecker,jnini2025gauss,wang2025discovery,kiyani2025optimizing,daryakenari2026representation,jnini2026curvature,webb2026optimisation}, however, have demonstrated the superiority of second-order optimization methods, such as the Gauss--Newton method, Newton's method, quasi-Newton methods (like BFGS), or natural gradient descent,
over conventional first-order optimizers 
in the realm of scientific machine learning
as they accelerate convergence and better navigate complex and severely ill-conditioned loss landscapes.

In this work,
we develop convergence analyses for the classical regularized Gauss--Newton and Newton method for training PINNs in the overparameterized regime.
Our results extend the analysis of first-order PINN training \cite{jiang2023global,gao2026global} to second-order optimization methods.
To establish these results, we leverage techniques developed for second-order training in regression problems \cite{riedl2026convergence} and elevate them to a function-space setting with physics-informed objectives, where the presence of differential operators acting on the NN introduces substantial analytical challenges.
For our mathematical analysis,
we consider the problem of solving a second-order linear PDE with Dirichlet boundary conditions of the form
\begin{alignat}{2}
\label{eq:PDE}
\begin{aligned}
    \CL u 
    &= h,
    \qquad&&\text{in }
    U\subset\bbR^d, \\
    u
    &= 0,
    \qquad&&\text{on }
    \partial U,
\end{aligned}
\end{alignat}
where $\CL$ denotes a second-order linear PDE operator,
by employing a fully-connected single-hidden-layer NN~$g_\theta^N:\bbR^d\rightarrow\bbR$ with $N$ neurons of the form
\begin{equation}
    \label{eq:NN}
    g_\theta^N(x)
    = \chi(x)f_\theta^N(x)
    \quad\text{ with }\quad
    f_\theta^N(x)
    = \frac{1}{N^\beta} \sum_{i=1}^N c^i \sigma(w^i \cdot x + \eta^i).
\end{equation}
The vector $\theta=(\theta^i)_{i=1}^N$ with $\theta^i=(c^i,w^i,\eta^i)$ collects all $N(d+2)$ trainable NN parameters
and $\chi\in\CC_b^3(\bbR^d)$ is a fixed function that vanishes on the boundary $\partial U$ ensuring that $g_\theta^N$ automatically satisfies the boundary conditions in \eqref{eq:PDE}.
The factor~$1/N^\beta$ in \eqref{eq:NN} with $\beta\in(1/2,1)$ is a normalization/scaling~\cite{spiliopoulos2025mathematical,sirignano2019scaling,sirignano2023pde,riedl2025global,riedl2026convergence}.
Throughout this paper, we consider strong Sobolev solutions to the PDE~\eqref{eq:PDE}, meaning functions
$u\in\CH^2_{(0)}(U):= \CH^2(U) \cap \CH^1_{0}(U)$ such that $\CL u=h$ almost everywhere in $U$,
with derivatives understood in the weak sense, the Dirichlet boundary condition satisfied in the trace sense.

\begin{remark}[Extensions of the PDE~\eqref{eq:PDE}]
    The restriction to homogeneous Dirichlet boundary conditions is made for notational convenience.
    For prescribed boundary data $u=b$ on $\partial U$,
    let $\bar b\in\CH^2(U)$ denote a lifting with trace $b$.
    Then $v\coloneqq u-\bar b$ satisfies
    \begin{alignat*}{2}
    \begin{aligned}
        \CL v
        &= h-\CL\bar b,
        \qquad&&\text{in }
        U\subset\bbR^d, \\
        v
        &= 0,
        \qquad&&\text{on }
        \partial U.
    \end{aligned}
    \end{alignat*}
    Since $\CL:\CH^2(U)\to L^2(U)$ is bounded, the modified source $h-\CL\bar b$ belongs to $L^2(U)$.
    Thus our analysis applies to this homogeneous problem, and the original solution is given by $\bar b+g_\theta^N$.
    In particular, linearity gives $\CL(\bar b+g_\theta^N)-h=\CL g_\theta^N-(h-\CL\bar b)$,
    so the loss and its parameter derivatives coincide with those of the transformed problem.
\end{remark}

\paragraph{Training with Newton's and the Gauss--Newton method.}
We train the NN parameters~$\theta$ of the NN \eqref{eq:NN} such that they minimize the loss (squared PDE residual)
\begin{equation}
    \label{eq:loss}
    L^N(\theta)
    = \frac{1}{2} \N{\CL g^N_{\theta}-h}_{L^2(U)}^2
    = \frac{1}{2} \int_U \left(\CL g_\theta^N(x) -h(x)\right)^2 dx
\end{equation}
by employing the regularized (Gauss--)Newton method
\begin{equation}
    \label{eq:regularizedNewton}
    \theta_{k+1}
    = \theta_k + \alpha z^N_k
    \quad\text{ with }\quad
    \theta_0^i \stackrel{\text{i.i.d.\@}}{\sim} \mu_0
    \text{ for all }i=1,\dots,N,
\end{equation}
where $\alpha>0$ is a fixed step size satisfying $\alpha\widetilde\lambda_{\max}<2$, with $\widetilde\lambda_{\max} \coloneqq\|\CL T_{\CB_0^*}\|_{\mathrm{op}}<1$.
Here, $T_{\CB_0^*}$ is the PDE-associated Newton neural tangent kernel (PDE--NNTK) integral operator defined in \eqref{eq:T_NNTK*}, see Remark~\ref{rem:spectralrelationNNTK}.
The admissible step-size interval contains $(0,2]$, and in particular the full step $\alpha=1$.
Furthermore, $z^N_k$ denotes the \mbox{(Gauss--)Newton} step, given implicitly by the  $N(d+2)\times N(d+2)$-dimensional linear system
\begin{equation}
    \label{eq:regularizedNewtonStep}
    \left(\gamma^N \Id_{N(d+2)} + P^N_{\theta_k}\right) z^N_k
    = -\nabla_\theta L^N(\theta_k),
\end{equation}
whenever a unique solution exists, and set to zero otherwise.
The preconditioner $P^N_\theta$ in \eqref{eq:regularizedNewtonStep} is given by
\begin{equation}
    P^N_\theta
    = \begin{cases}
        G^N_\theta := 
        \int_U \CL\nabla_\theta g_\theta^N(x) \left(\CL\nabla_\theta g_\theta^N(x)\right)^\top dx,
        & \quad\text{ for the Gauss--Newton method},\\
        \nabla^2_\theta L^N(\theta),
        & \quad\text{ for the Newton method},
    \end{cases}
\end{equation}
and the regularizer $\gamma^N$ scales as $\gamma^N = \gamma/N^{2\beta-1}$, where the regularizer parameter~$\gamma$ is assumed to be strictly positive ($\gamma>0$) throughout the manuscript.
For notational convenience,
let us note that the Hessian decomposes as
\begin{equation*}
    \nabla^2_\theta L^N(\theta)
    = G^N_\theta + S^N_\theta
    = G^N_\theta - \int_U \left(h(x)-\CL g_\theta^N(x)\right) \CL\nabla_\theta^2 g_\theta^N(x) \, dx,
\end{equation*}
where we denote by
$G^N_\theta$ the Gram--Gauss--Newton matrix
and by $S^N_\theta$ the non-convex contribution of the Hessian,
see Section~\ref{sec:evolution} and in particular \eqref{eq:G} and \eqref{eq:S} for more details.

\subsection{Contributions}
\label{sec:contributions}

To provide a theoretical foundation for the recent empirical successes of training PINNs with second-order optimization methods~\cite{muller2023achieving,jnini2026curvature,webb2026optimisation},
we develop a rigorous convergence analysis of the regularized Newton and Gauss--Newton methods~\eqref{eq:regularizedNewton}--\eqref{eq:regularizedNewtonStep} for solving second-order linear PDEs of the form~\eqref{eq:PDE}.
For an NN as in \eqref{eq:NN}, we study the asymptotic convergence behavior of the discrete-time algorithm as the number of hidden units~$N$ becomes large and subsequently as the number of training steps~$k$ tends to infinity.
The results of the present paper build upon our earlier convergence analysis of Newton’s method for NN regression~\cite{riedl2026convergence} to PINNs, establishing convergence of regularized Newton and Gauss--Newton training in the overparameterized regime.
This provides a second-order counterpart to the first-order PINN analysis of~\cite{jiang2023global}.
The presence of the PDE operator in the loss~\eqref{eq:loss} and the infinite-dimensional nature of the PDE solution space introduce substantial additional mathematical challenges, requiring in particular control of spatial derivatives of the NN and its associated kernel operators, as well as an analysis of the limiting training dynamics without assuming a spectral gap.
Our main contributions are:
\begin{itemize}[label=\large \textbullet,labelsep=8pt,leftmargin=25pt]%,topsep=-4pt,itemsep=-1pt
    \item \emph{PDE--NNTK.}
    We derive the PDE--NNTK integral operator $T_{\CB_{\theta}^N}$ as defined in \eqref{eq:NNTK}, which governs the discrete-time training evolution~\eqref{eq:NNfunction_evolution_prelimit} of the NN function~\eqref{eq:NN} under the regularized Newton and Gauss--Newton methods.
    This operator incorporates both the differential operator of the PDE~\eqref{eq:PDE} and the solution of the regularized linear system~\eqref{eq:regularizedNewtonStep} defining the parameter update.
    At finite width, it is random at initialization and varies during training.
    We identify its deterministic infinite-width limit~\eqref{eq:T_NNTK*}, which remains constant during training and is shared by both optimization methods.
    
    Let us remark that we distinguish the PDE--NNTK integral operator $T_{\CB_{\theta}^N}$ from its associated residual evolution operator $\CL T_{\CB_{\theta}^N}$.
    In the infinite-width limit, the spectral properties of $\CL T_{\CB_0^*}$ govern the convergence of the PDE residual.
    \item \emph{Convergence to the infinite-width limit.}
    Our first main result, Theorem~\ref{thm:NTK_limit}, establishes that, as the number of hidden units $N\rightarrow\infty$, the finite-width NN training dynamics~\eqref{eq:NNfunction_evolution_prelimit} converge in probability to the deterministic limit dynamics~\eqref{eq:NNfunction_evolution_limit}.
    This convergence holds in $\CH^2(U)$, uniformly over any fixed number $K$ of training steps, and is quantified by explicit rates in $N$.
    The result applies to both the regularized Newton and Gauss--Newton methods, showing that they share the same infinite-width training dynamics. For Newton's method, the additional non-convex Hessian contribution becomes negligible relative to the regularization as $N\to\infty$, uniformly over any fixed number of training steps, explaining the equivalence of the two limiting dynamics.
    \item \emph{Global convergence and convergence rates.}
    Our second main result concerns the long-time behavior of the limiting training dynamics~\eqref{eq:NNfunction_evolution_limit}.
    Theorem~\ref{thm:convergence} establishes that the PDE residual converges to zero in $L^2(U)$ as the number of training steps $k\rightarrow\infty$, and hence that the limiting loss vanishes, without assuming a spectral gap for the residual evolution operator.
    This result holds for a range of constant step sizes that includes the full step $\alpha=1$.
    Under an additional stability estimate for the PDE operator, Corollary~\ref{cor:solution_convergence} shows that residual convergence also implies convergence of the limiting NN to the unique PDE solution.
    Theorem~\ref{lem:finite_mode_rate} further provides explicit geometric convergence rates for the residual projected onto any fixed finite number of leading eigenfunctions, quantifying how the regularizer parameter $\gamma$ controls the convergence rates across these modes.
    \item \emph{Mitigation of spectral bias.}
    These convergence rates show, in particular, that choosing the regularizer parameter $\gamma$ sufficiently small reduces the disparity between the decay rates of any fixed finite collection of leading eigenmodes, thereby mitigating the spectral bias inherent in gradient descent.
    This is attributed to the spectral transformation $\CL T_{\CB_0^*} = \CL T_{B_0^*}\left(\gamma\Id_{L^2(U)}+\CL T_{B_0^*}\right)^{-1}$, which improves conditioning on these finitely many modes compared to gradient descent, see Remark~\ref{rem:spectralrelationNNTK}.
    Here, $T_{B_0^*}:L^2(U)\to\CH^2_{(0)}(U)$ denotes the limiting PDE--NTK integral operator associated with gradient descent, defined in \eqref{eq:T_NTK*}.
    The regularizer parameter $\gamma$ thus acts as a spectral filter: as $\gamma$ decreases, the eigenvalues of $\CL T_{\CB_0^*}$ associated with these leading modes approach one.
    This observation applies only to finitely many modes and therefore does not imply a uniform geometric convergence rate across the entire infinite-dimensional spectrum.
\end{itemize}

\paragraph{Mathematical challenges.}
Before outlining the structure of the paper in Section~\ref{sec:organization}, let us first highlight the core mathematical challenges of our analysis.

Firstly, the implicit parameter updates~\eqref{eq:regularizedNewtonStep} of both methods pose a central challenge for our analysis, as each training step requires solving a linear system whose dimension grows with the NN width.
For Newton's method, the potentially indefinite Hessian introduces the additional difficulty of ensuring invertibility of the regularized system throughout training.
To analyze these updates, we exploit the block-diagonal plus Gram structure of both systems and prove an operator push-through identity to obtain neuron-wise representations.
For both methods, our analysis identifies the scaling $\gamma^N=\gamma/N^{2\beta-1}$, which balances regularization with the Gram--Gauss--Newton term.
While the Gauss--Newton system is positive definite for every $\gamma>0$, Newton's method additionally requires control of the non-convex Hessian contribution.
By establishing operator-norm estimates for this contribution, we ensure positive definiteness with high probability over any fixed training horizon for sufficiently large NN widths.

Secondly, the PDE operator introduces further difficulties by requiring control of spatial derivatives of the NN and its associated kernels.
We establish bounds on spatial derivatives up to second order, which are needed to control the effect of the PDE operator throughout training.
These estimates allow us to prove convergence of the finite-width training trajectories to their infinite-width limit in $\CH^2(U)$, uniformly over any fixed number of training steps.

Lastly, the infinite-dimensional nature of the PDE setting poses a central challenge for establishing global convergence, as the residual evolution operator $\CL T_{\CB_0^*}$ lacks a spectral gap and may have a nontrivial kernel.
To overcome these difficulties, we exploit the PDE structure and establish orthogonality of the limiting residuals to the kernel before using the dominated convergence theorem to pass from decay of individual spectral components to decay of the full residual.
    
\subsection{Organization}
\label{sec:organization}

After discussing related works in what follows,
Section~\ref{sec:main} presents the main theoretical results
for training PINNs with the regularized Newton and
Gauss--Newton methods.
We first state the mathematical assumptions in
Section~\ref{sec:assumptions} and then derive the discrete-time
evolution of the NN function in Section~\ref{sec:evolution},
identifying the associated PDE--NNTK integral operator.
Section~\ref{sec:convergenceNinfty} establishes convergence
of the finite-width dynamics to their deterministic
infinite-width limit in Theorem~\ref{thm:NTK_limit},
whose proof is given in
Section~\ref{sec:proof:convergenceNinfty}.
Section~\ref{sec:globalconvergence} studies the long-time
behavior of the limiting dynamics, establishing convergence
to zero PDE residual and loss in
Theorem~\ref{thm:convergence} and analyzing the associated
convergence rates.
The proofs of these results are provided in
Section~\ref{sec:proof:globalconvergence}.

\subsection{Related works}
\label{sec:relatedworks}

The numerical analysis of PINNs has attracted considerable attention in recent years, see the review paper~\cite{de2024numerical} and the references therein.
Our literature review below focuses on works closely related to the objective of this paper, which is to establish global convergence of the regularized Newton and Gauss--Newton methods for training PINNs in the overparameterized regime.

\emph{Approximation error}:
For second-order linear elliptic PDEs, \cite{jiao2021rate} establishes approximation and statistical error estimates with explicit dependence on the NN architecture and the number of training samples.
For broader classes of linear PDEs, the authors of \cite{shin2023error,zeinhofer2025unified} derive error estimates using stability properties of the underlying PDE formulations to relate the terminal loss to the approximation error of the PDE solution.
These analyses assume sufficiently accurate loss minimization but do not establish convergence of the training algorithms used to achieve it.

\emph{Convergence analyses of first-order optimization methods}:
For gradient descent, the NTK analysis of~\cite{wang2022when} characterizes the infinite-width dynamics of PINNs and identifies discrepancies between the convergence rates of different loss components.
For overparameterized two-layer NNs, \cite{gao2023gradient,luo2024two} establish global convergence of gradient flow for the empirical PDE residual losses, with~\cite{gao2023gradient} additionally analyzing discrete-time gradient descent under suitable step-size conditions.
These convergence results rely on positive definiteness of the limiting NTK matrices.
Moving beyond convergence of a fixed empirical loss and closely related to our setting, \cite{jiang2023global,sirignano2026global} analyze continuous-time gradient descent with clipping for continuous PDE residual losses associated with linear and semilinear elliptic PDEs, respectively.
Both works identify deterministic infinite-width dynamics and establish convergence of the limiting PDE residual to zero as training time tends to infinity, without assuming a spectral gap.
Under their respective assumptions on the PDE, they also establish convergence of the limiting NN to the PDE solution.

For stochastic gradient descent, the work \cite{jin2026convergence} studies overparameterized two-layer PINNs for the Poisson equation.
Under suitable width and step-size conditions, the authors control the training trajectories with high probability and establish geometric decay of the expected empirical loss restricted to this event.
The analysis relies on maintaining positive definiteness of the relevant NTK matrices throughout training.

For Adam, the author of \cite{gao2026global} establishes convergence guarantees for overparameterized two-layer PINNs for the Poisson equation, accounting for momentum and adaptive preconditioning.
Under explicit conditions on the network width and optimizer parameters, the empirical loss decays geometrically for full-batch Adam, while the mini-batch result establishes geometric decay of its expectation restricted to a high-probability event.

\emph{Convergence analysis of second-order optimization methods}:
The analysis of~\cite{wang2026gradient} shows how second-order optimization can mitigate conflicts between gradients of different PINN loss components and establishes a connection between SOAP and Newton's method.

For Gauss--Newton methods, the work \cite{hao2024gauss} establishes local convergence for NN discretizations of variational PDE problems, assuming initialization sufficiently close to a stationary point.
For empirical PINN losses, the authors of \cite{bonfanti2024challenges} characterize the residual dynamics under the Gauss--Newton flow and show that, when the residual Jacobian has full row rank throughout training, all residual components decay exponentially at the same rate.
This provides a theoretical explanation for the mitigation of spectral bias.

For natural gradient descent, the paper \cite{xu2026fast} establishes quadratic convergence with full steps for the empirical loss of overparameterized two-layer PINNs with smooth activation functions, under conditions ensuring positive definiteness of the limiting NTK matrix.

For Newton's method, the authors of \cite{rathore2024challenges} analyze a hybrid algorithm combining gradient descent with damped Newton steps.
Under suitable regularity and Polyak--{\L}ojasiewicz assumptions, the gradient descent phase reaches a neighborhood of a minimizer in which the Newton phase achieves geometric convergence of the empirical PINN loss at a rate independent of the condition number.

\subsection{Notation}
\label{sec:notation}

For a multi-index $\nu=(\nu_1,\dots,\nu_d)\in\bbN_0^d$,
we write $\abs{\nu}\coloneqq\sum_{j=1}^d\nu_j$ and $D_\nu\coloneqq\partial_{x_1}^{\nu_1}\cdots\partial_{x_d}^{\nu_d}$,
with derivatives understood in the weak sense.
The spaces $L^p(U)$, $1\leq p\leq\infty$, denote the usual Lebesgue spaces.
We denote by $\CH^m(U)=W^{m,2}(U)$ the Sobolev space of functions whose weak derivatives through order $m$ belong to $L^2(U)$.
Throughout the manuscript, we use the equivalent Sobolev norm $\N{v}_{\CH^2(U)} \coloneqq \sum_{0\leq\abs{\nu}\leq2} \N{D_\nu v}_{L^2(U)}$.
The space $\CH^1_0(U)$ is the closure of $\CC^\infty_c(U)$ in $\CH^1(U)$ consisting of the functions with zero trace.
We write $\CH^2_{(0)}(U) \coloneqq \CH^2(U)\cap\CH^1_0(U)$.

For a vector $v$, $\N{v}_2$ denotes its Euclidean norm, and for a matrix $A\in\bbR^{a\times b}$, $\N{A}_{\mathrm{op}} \coloneqq \sup_{\N{v}_2\leq1}\N{Av}_2$.
For normed spaces $X$ and $Y$ and a bounded linear operator $T:X\to Y$, we write $\N{T}_{X\rightarrow Y} \coloneqq \sup_{\N{v}_X\leq1}\N{Tv}_Y$.
When the underlying spaces are clear, this norm is denoted by $\N{T}_{\mathrm{op}}$.
We write $\Id_n$ for the identity matrix on $\bbR^n$ and $\Id_X$ for the identity operator on a function space $X$.

By $C$ we denote generic positive constants whose values may change from one occurrence to the next.
Unless stated otherwise, these constants are independent of the network width $N$ and the iteration index $k$.
Their dependence on fixed quantities is indicated by notation such as $C(\alpha,\gamma,K,\delta,\sigma,\CD,\mu_0)$.
Here, $\CD$ summarizes the deterministic problem data,
including the dimension $d$, the domain $U$, the PDE operator $\CL$, the source term $h$, and the fixed boundary function $\chi$.
Dependence on $\sigma$ refers to the properties of the activation function, while dependence on $\mu_0$ refers to the initialization.

For a generic single-neuron parameter $\vartheta=(c,w,\eta)\in\bbR^{d+2}$,
define the scalar single-neuron function
\begin{equation}
    \phi_\vartheta(x)\coloneqq\chi(x)c\sigma(w\cdot x+\eta)
\end{equation}
as well as
\begin{align}
    \label{eq:s}
    s_{\vartheta}(x)
    &\coloneqq \left(\sigma(w \cdot x+\eta), c\sigma'(w \cdot x+\eta)x^\top,c\sigma'(w \cdot x+\eta)\right)
\end{align}
and
\begin{align}
    \label{eq:h}
    h_{\vartheta}(x)
    &\coloneqq\left(\begin{matrix}0 & \sigma'(w \cdot x + \eta)x^\top & \sigma'(w \cdot x + \eta) \\
    \sigma'(w \cdot x + \eta)x & c\sigma''(w \cdot x + \eta)xx^\top & c\sigma''(w \cdot x + \eta)x \\
    \sigma'(w \cdot x + \eta) & c\sigma''(w \cdot x + \eta)x^\top & c\sigma''(w \cdot x + \eta) \end{matrix}\right).
\end{align}

\section{Discussion of the Main Results}
\label{sec:main}

This section discusses the main theoretical contributions of the paper.
After stating the assumptions on the PDE, the NN and its initialization in Section~\ref{sec:assumptions},
Section~\ref{sec:evolution} derives the discrete-time training evolution of the network function under the regularized Newton and Gauss--Newton methods.
Section~\ref{sec:convergenceNinfty} then establishes convergence of the finite-width dynamics to their infinite-width limit and quantifies the approximation error,
before Section~\ref{sec:globalconvergence} proves convergence of the PDE residual to zero for the limiting dynamics and analyzes the associated convergence rates.

\subsection{Assumptions}
\label{sec:assumptions}

Let us start by stating all assumptions which hold throughout this manuscript.
We cluster them into assumptions related to the PDE~\eqref{eq:PDE}, which we summarize in Assumption~\ref{def:PDE_assumptions}, and assumptions on the NN and its initialization listed thereafter in Assumption~\ref{def:NN_assumptions}.

\begin{assumption}[PDE~\eqref{eq:PDE}]
    \label{def:PDE_assumptions}
	Throughout we assume that
	\begin{enumerate}[label=A\arabic*,labelsep=10pt,leftmargin=35pt]
        \item\label{asm:solution} the PDE~\eqref{eq:PDE} admits a strong Sobolev solution $u\in\CH^2_{(0)}$,
        \item\label{asm:D} the domain~$U\subset\bbR^d$ of the PDE~\eqref{eq:PDE} is an open connected set with boundary $\partial U$ of class $\CC^{3,\varsigma}$ for some $\varsigma\in(0,1)$, has finite volume $\vol{U}$, and is bounded by $C_U$, i.e., $\sup_{x\in U}\N{x}_2\leq C_U$,
        \item\label{asm:h} the source term $h\in L^2(U)$, and we denote $C_h\coloneqq\N{h}_{L^2(U)}$,
        \item\label{asm:CL} the PDE operator~$\CL$ satisfies for some constant $k_\CL>0$ that
        \begin{equation}
            \label{eq:asm:Lipschitz}
        \begin{split}
            \abs{\CL v_1(x)-\CL v_2(x)}
            \leq k_\CL\!\! \sum_{0\leq\abs{\nu}\leq2}\!\!\abs{D_\nu v_1(x)-D_\nu v_2(x)}
        \end{split}
        \end{equation}
        for all $v_1,v_2\in\CH^2(U)$ and a.e.\@ $x\in U$.
    \end{enumerate}
\end{assumption}

\begin{remark}
    \label{rem:CL_linearity}
    Since $\CL$ is linear, $\CL0=0$, so that \eqref{eq:asm:Lipschitz} yields both the bound $\abs{\CL v(x)}\leq k_\CL\sum_{0\leq\abs{\nu}\leq2}\abs{D_\nu v(x)}$ for a.e.\@ $x\in U$ and, with the convention $\N{v}_{\CH^2(U)}=\sum_{0\leq\abs{\nu}\leq2}\N{D_\nu v}_{L^2(U)}$, the operator bound
    \begin{equation}
        \label{eq:L_H2_L2}
        \N{\CL v}_{L^2(U)}\leq k_\CL\N{v}_{\CH^2(U)}
        \quad\text{ for all }v\in\CH^2(U).
    \end{equation}
    Throughout, pointwise identities and inequalities involving $\CL$ are understood almost everywhere in the spatial variable on which $\CL$ acts.
\end{remark}

\begin{assumption}[Neural network \eqref{eq:NN}]\label{def:NN_assumptions}
	Throughout we assume that the NN is such that
	\begin{enumerate}[label=B\arabic*,labelsep=10pt,leftmargin=35pt]
        \item\label{asm:NN_sigma} the nonlinearity~$\sigma$ of the NN is non-constant, bounded (i.e., $\abs{\sigma}\leq C_\sigma$), and $L_\sigma$-Lipschitz continuous,
        \item\label{asm:NN_sigma'} the derivatives of the nonlinearity $\sigma\in\CC^4_b(\bbR)$ are continuous and bounded up to and including order four, (i.e., $\abs{\sigma^{(\ell)}} \leq C_{\sigma^{(\ell)}}$ for $\ell=1,\dots,4$),
        \item\label{asm:NN_mu0} the randomly initialized NN parameters~$\theta_0^i = (c^i_0,w^i_0,\eta^i_0)$ are i.i.d.\@ and drawn from a distribution $\mu_0\in\CP(\bbR\times\bbR^{d}\times\bbR)$
        which is such that 
        \begin{enumerate}[label=(\roman*),labelsep=10pt,leftmargin=35pt]
            \item $c^i_0$ is independent from $(w^i_0,\eta^i_0)$, \label{asm:NN_mu0i}
            \item the marginal distribution $\mu_{0,c}$ of $c^i_0$ is mean-zero and compactly supported, i.e., $\absnormal{c^i_0}\leq C_{c,0}(\mu_0)$, \label{asm:NN_mu0ii}
            \item the marginal distribution $\mu_{0,(w,\eta)}$ of $(w^i_0,\eta^i_0)$ assigns positive probability to every set with positive Lebesgue measure, \label{asm:NN_mu0iii}
            \item the $q$-th moment $\CM_q\coloneqq\bbE\N{w^1_0}_2^{q}$ is finite for some $q$ such that $q \geq 8$ and $q > \frac{2}{1-\beta}$,
            \label{asm:NN_mu0iv}
        \end{enumerate}
        \item\label{asm:NN_chi} the boundary function $\chi\in\CC^3_b(\bbR^d)$ satisfies $\chi>0$ in $U$, $\chi=0$ on $\partial U$, and $\nabla\chi(x)\cdot\mathbf{n}_x\neq0$ for every $x\in\partial U$, where $\mathbf{n}_x$ denotes an outward unit normal of $\partial U$ at $x$.
    \end{enumerate}
\end{assumption}

\noindent
For notational convenience, we denote $\omega_i\coloneqq 1 + \N{w^i_0}_2^2$ for the NN parameter $\theta^i_0$.

\subsection{Evolution of the NN function during training}
\label{sec:evolution}

To start with, we derive the PDE--NNTK integral operator representation \eqref{eq:NNfunction_evolution_prelimit} of the discrete-time training evolution of the NN function~$g_\theta^N$ when trained with the (Gauss--)Newton method \eqref{eq:regularizedNewton}--\eqref{eq:regularizedNewtonStep}.

To do so, we first derive the NN parameter derivatives of the loss and write the regularized Newton and Gauss--Newton updates~\eqref{eq:regularizedNewtonStep} in operator form as \eqref{eq:zN_reformulation}--\eqref{eq:zetaN_reformulation}.

\paragraph{NN parameter derivatives.}
We define $J_{\theta}^N(x)\in\bbR^{1 \times N(d+2)}$ as the Jacobian of the NN function~$f_{\theta}^N$ for a data point $x$,
and the submatrix $J_{\theta^i}^N(x)\in\bbR^{1 \times (d+2)}$ as the Jacobian w.r.t.\@ the NN parameters $\theta^i$ of the $i$-th neuron.
We further define
$H^N_{\theta}(x)\in\bbR^{N(d+2)\times N(d+2)}$ as the Hessian of the NN function~$f_{\theta}^N$ for a data sample~$x$,
which is block-diagonal with $N$ blocks of size $(d+2)\times (d+2)$,
where individual blocks correspond to different neurons.

\textit{Gradient of the loss.}
For the gradient~$\nabla_\theta L^N(\theta)\in\bbR^{N(d+2)}$ of the loss~$L^N$,
we compute with chain rule that
\begin{equation}
\begin{split}
    \nabla_\theta L^N(\theta)
    &= \int_U \left(\CL g_\theta^N(x) -h(x)\right) \nabla_\theta \left(\CL g_\theta^N(x)\right) dx 
    = \int_U \left(\CL g_\theta^N(x) -h(x)\right) \CL \nabla_\theta g_\theta^N(x)\, dx \\
    &= \int_U \left(\CL g_\theta^N(x) -h(x)\right) \left(\CL \left(\chi(x) J_{\theta}^N(x)\right)\right)^\top dx.
\end{split}
\end{equation}

\textit{Hessian of the loss.}
For the Hessian~$\nabla^2_\theta L^N(\theta)\in\bbR^{N(d+2)\times N(d+2)}$ of the loss~$L^N$,
we have that
\begin{equation}
\begin{split}
    \nabla^2_\theta L^N(\theta)
    &= \int_U \CL \nabla_\theta g_\theta^N(x) \left(\CL \nabla_\theta g_\theta^N(x)\right)^\top + \left(\CL g_\theta^N(x) -h(x)\right) \CL \nabla^2_\theta g_\theta^N(x)\, dx \\
    &= \int_U \left(\CL \left(\chi(x) J_{\theta}^N(x)\right)\right)^\top \CL\left(\chi(x) J_{\theta}^N(x)\right) + \left(\CL g_\theta^N(x)-h(x)\right) \CL  \left(\chi(x) H_\theta^N(x)\right) dx \\
    &= G_\theta^N + S_\theta^N,
\end{split}
\end{equation}
where we denoted for notational simplicity
the Gram--Gauss--Newton matrix by
\begin{equation}
    \label{eq:G}
    G_\theta^N
    = \int_U \CL \nabla_\theta g_\theta^N(x) \left(\CL \nabla_\theta g_\theta^N(x)\right)^\top dx
    = \int_U \left(\CL \left(\chi(x)J_{\theta}^N(x)\right)\right)^\top \CL \left(\chi(x)J_{\theta}^N(x)\right) dx
\end{equation}
and the non-convex contribution by
\begin{equation}
    \label{eq:S}
    S_\theta^N
    = - \int_U \left(h(x)-\CL g_\theta^N(x)\right) \CL \nabla^2_\theta g_\theta^N(x)\, dx
    = - \int_U \left(h(x)-\CL g_\theta^N(x)\right) \CL  \left(\chi(x) H_\theta^N(x)\right) dx.
\end{equation}
That is, the Hessian~$\nabla^2_\theta L^N(\theta) = G_\theta^N + S_\theta^N$ 
of the loss~$L^N(\theta)$ decomposes
into the sum of the Gram--Gauss--Newton matrix $G_\theta^N\in\bbR^{N(d+2)\times N(d+2)}$ and the non-convex contribution $S_\theta^N\in\bbR^{N(d+2)\times N(d+2)}$,
where the latter is block-diagonal with its individual blocks corresponding to different neurons.

\paragraph{Regularized Newton and Gauss--Newton updates.}
For both the regularized Newton and Gauss--Newton methods,
the update~$z^N_k$ in \eqref{eq:regularizedNewtonStep} is therefore the solution of
\begin{equation}
    \label{eq:regularizedNewton_DG}
    \left(D_{\theta_k}^N + G_{\theta_k}^N\right) z^N_k = -\nabla_\theta L^N(\theta_k),
\end{equation}
where the block-diagonal matrix $D_\theta^N$ is defined by
\begin{equation}
    \label{eq:D}
    D_\theta^N
    \coloneqq
    \begin{cases}
        \gamma^N\Id_{N(d+2)},
        & \quad\text{ for the Gauss--Newton method},\\[1mm]
        \gamma^N\Id_{N(d+2)}+S_\theta^N,
        & \quad\text{ for the Newton method}.
    \end{cases}
\end{equation}
Its $N$ blocks of size $(d+2)\times(d+2)$ on the diagonal correspond to individual neurons.

\paragraph{Kernel representation of the NN function training evolution.}
Using a Taylor expansion around $\theta_k$ and denoting by $(z^N_k)^i\in\bbR^{d+2}$ the part of the \mbox{(Gauss--)}Newton update~$z^N_k$ in \eqref{eq:regularizedNewtonStep} that updates the NN parameters $\theta^i_k$ of the $i$-th neuron in \eqref{eq:regularizedNewton},
we obtain for the training evolution of the finite-width NN function~\eqref{eq:NN} evaluated at a single data point $x\in\bbR^d$ that
\begin{equation}
    \label{eq:NTK_limit:evolution_fk+1_1}
\begin{split}
    g_{\theta_{k+1}}^N(x)
    &= g_{\theta_{k}}^N(x) + \nabla_\theta g^N_{\theta_k}(x) \cdot (\theta_{k+1}-\theta_k) + R_k^N(x)
    = g_{\theta_{k}}^N(x) + \alpha\chi(x) J_{\theta_k}^N(x) z^N_k + R^N_k(x) \\
    &= g_{\theta_{k}}^N(x) + \alpha \sum_{i=1}^N \chi(x) J_{\theta^i_k}^N(x) \left(z^N_k\right)^i + R^N_k(x),
\end{split}
\end{equation}
where the Taylor remainder is given by $R_k^N(x) \coloneqq \alpha^2\chi(x)\int_0^1(1-t)\, (z_k^N)^\top H^N_{\theta_k+t\alpha z_k^N}(x) z_k^N\,dt$.
We establish its regularity and quantitative bounds for it in Lemma~\ref{lem:taylor_remainder}.
Following the arguments of \cite[Section~2.2]{riedl2026convergence},
we now leverage the decomposition $\gamma^N \Id_{N(d+2)} + P^N_{\theta_k} = D_{\theta_k}^N+G_{\theta_k}^N$ of the regularized (Gauss--)Newton system~\eqref{eq:regularizedNewtonStep}.
Provided that $D_{\theta_k}^N$ in \eqref{eq:D} is positive definite
(which is automatic for the Gauss--Newton method and, under \eqref{eq:thm:asm:Nlargeenough}, holds for Newton's method for a fixed number $K$ of training steps with probability at least $1-\delta$, see Lemmas~\ref{lem:estimate_Dinv_operatornorm} and~\ref{prop:NTK_limit_bounds_z}), 
the continuous operator analogue of the Woodbury push-through identity, established in Lemma~\ref{lem:pushthrough}, gives
\begin{equation}
    \label{eq:zN_reformulation}
\begin{split}
    z_k^N
    &= \left(D_{\theta_k}^N\right)^{-1} \int_U \left(\CL \left(\chi(x)J_{\theta_k}^N(x)\right)\right)^\top \zeta^N_k(x) \, dx
\end{split}
\end{equation}
with
\begin{equation}
    \label{eq:zetaN_reformulation}
\begin{split}
    \zeta_k^N
    \coloneqq
    \left(\Id_{L^2(U)} + \frac{1}{\gamma}\CL T_{B_{\theta_k}^N} \right)^{-1}
    \left(h - \CL g_{\theta_k}^N\right),
\end{split}
\end{equation}
where the PDE--NTK integral operator $T_{B_\theta^N}: L^2(U) \to \CH^2_{(0)} \subset L^2(U)$ is defined by
\begin{equation}
    \label{eq:T_NTK}
    \left[T_{B_\theta^N} v\right]\!(x)
    \coloneqq \int_U B_\theta^N(x,x') v(x') \, dx'
\end{equation}
with the kernel $B_\theta^N: U \times U \to \bbR$ defined
according to
\begin{equation}
    \label{eq:NTK}
    B_\theta^N(x,x')
    \coloneqq \,\chi(x) J_{\theta}^N(x) \left(\frac1\gamma D_{\theta}^N\right)^{-1} \left(\CL \left(\chi(x') J_{\theta}^N(x')\right)\right)^\top.
\end{equation}
The invertibility of the operator appearing in \eqref{eq:zetaN_reformulation} is established in Lemma~\ref{lem:estimate_zeta}.

Since $D_{\theta_k}^N$ is block-diagonal with blocks $D_{\theta_k^i}^N$,
also its inverse~$\big(D_{\theta_k}^N\big)^{-1}$ is block-diagonal with blocks $\big(D_{\theta_k^i}^N\big)^{-1}$.
Exploiting this, \eqref{eq:zN_reformulation} reads neuron-wise
\begin{equation}
    \label{eq:zN_reformulation_neuronwise}
    \big(z_k^N\big)^i= \big(D_{\theta_k^i}^N\big)^{-1} \int_U \big(\CL \big(\chi(x)J_{\theta_k^i}^N(x)\big)\big)^\top \zeta^N_k(x) \, dx
\end{equation}
for all $i=1\dots,N$.

Recalling the choice of the regularizer $\gamma^N=\frac{\gamma}{N^{2\beta-1}}$ together with the formulas $J^N_{\theta^i}(x) = \frac{1}{N^\beta}s_{\theta^i}(x)$ and $H^N_{\theta^i}(x) = \frac{1}{N^\beta}h_{\theta^i}(x)$,
and denoting by $d_{\theta^i}^N\coloneqq N^{2\beta-1}D^N_{\theta^i}$ the rescaled $i$-th block of $D_\theta^N$, see \eqref{eq:def_d} for a formula,
we can rewrite \eqref{eq:NTK_limit:evolution_fk+1_1} as
\begin{equation}
    \label{eq:NTK_limit:evolution_fk+1_3}
    \begin{split}
        g_{\theta_{k+1}}^N(x)
        &= g_{\theta_{k}}^N(x) + \alpha \int_U \sum_{i=1}^N \chi(x) J_{\theta^i_k}^N(x) \big(D_{\theta_k^i}^N\big)^{-1} \big(\CL \big(\chi(x')J_{\theta_k^i}^N(x')\big)\big)^\top \zeta^N_k(x') \, dx' + R^N_k(x) \\
        &= g_{\theta_{k}}^N(x) + \alpha \int_U \frac{1}{N}\sum_{i=1}^N \chi(x)s_{\theta_k^i}(x) \big(d_{\theta^i_k}^N\big)^{-1}  \big(\CL \big(\chi(x')s_{\theta_k^i}(x')\big)\big)^\top \zeta^N_k(x') \, dx'  + R^N_k(x)\\
        &= g_{\theta_{k}}^N(x) + \frac\alpha\gamma \int_U B^N_{\theta_k}(x,x') \zeta^N_k(x') \, dx'  + R^N_k(x)\\
        &= g_{\theta_{k}}^N(x) + \frac\alpha\gamma \left[T_{B_{\theta_k}^N} \zeta^N_k\right]\!(x) + R^N_k(x),
    \end{split}
\end{equation}
where we used in the next-to-last step that the kernel $B^N_\theta$ from \eqref{eq:NTK} admits the representation
\begin{equation}
    \label{eq:kernel_finite}
    B^N_\theta(x,x')
    =\frac1N\sum_{i=1}^N\chi(x)s_{\theta^i}(x)
    \left(\frac1\gamma d^N_{\theta^i}\right)^{-1}
    \left(\CL\left(\chi(x') s_{\theta^i}(x')\right)\right)^\top,
\end{equation}
whereas the last step uses the definition of the PDE--NTK integral operator $T_{B_\theta^N}$ from \eqref{eq:T_NTK}.
Inserting the expression for $\zeta_k^N$ from \eqref{eq:zetaN_reformulation} into \eqref{eq:NTK_limit:evolution_fk+1_3}
yields
\begin{equation}
    \label{eq:NTK_limit:evolution_fk+1_4}
    \begin{split}
        g_{\theta_{k+1}}^N(x)
        &= g_{\theta_{k}}^N(x) + \frac{\alpha}{\gamma} \left[ T_{B_{\theta_k}^N} \left(\left(\Id_{L^2(U)} + \frac{1}{\gamma}\CL T_{B_{\theta_k}^N} \right)^{-1}\left(h - \CL g_{\theta_k}^N\right)\right)\right](x) + R^N_k(x),
    \end{split}
\end{equation}
or, written in functional form,
\begin{equation}
    \label{eq:NTK_limit:evolution_fk+1_4_functional}
    \begin{split}
        g_{\theta_{k+1}}^N
        &= g_{\theta_{k}}^N + \frac{\alpha}{\gamma} T_{B_{\theta_k}^N} \left(\left(\Id_{L^2(U)} + \frac{1}{\gamma}\CL T_{B_{\theta_k}^N} \right)^{-1}\left(h - \CL g_{\theta_k}^N\right)\right) + R^N_k,
    \end{split}
\end{equation}

Defining the \emph{PDE-associated Newton neural tangent kernel (PDE--NNTK) integral operator} $ T_{\CB_{\theta}^N} : L^2(U) \to \CH^2_{(0)} \subset L^2(U)$ as
\begin{equation}
    \label{eq:NNTK}
    \left[ T_{\CB_{\theta}^N}v\right]\!(x) 
    \coloneqq \left[ \frac{1}{\gamma} T_{B_{\theta}^N} \left(\left(\Id_{L^2(U)} + \frac{1}{\gamma}\CL T_{B_{\theta}^N} \right)^{-1}v\right)\right](x),
\end{equation}
we can rewrite equation \eqref{eq:NTK_limit:evolution_fk+1_4_functional}
and obtain for the training dynamics of the finite-width NN function~\eqref{eq:NN}
\begin{equation}
        \label{eq:NNfunction_evolution_prelimit}
    \begin{split}
        g_{\theta_{k+1}}^N
        &= g_{\theta_{k}}^N + \alpha  T_{\CB_{\theta_k}^N} \left(h - \CL g_{\theta_{k}}^N\right) + R^N_k.
    \end{split}
\end{equation}
The PDE--NNTK integral operator $T_{\CB_\theta^N}$ maps the PDE residual to the NN update direction.
Its terminology emphasizes that the PDE operator enters $T_{\CB_\theta^N}$, both through the kernel $B_\theta^N$ from \eqref{eq:kernel_finite} and through \eqref{eq:NNTK}.

We call the composition $\CL T_{\CB_\theta^N}:L^2(U)\to L^2(U)$ the \emph{associated residual evolution operator}.
Writing $r_k^N\coloneqq h-\CL g_{\theta_k}^N$ for the PDE residual at training step $k$ and applying $\CL$ to \eqref{eq:NNfunction_evolution_prelimit} gives
\begin{equation}
    r_{k+1}^N
    =\left(\Id_{L^2(U)}-\alpha\CL T_{\CB_{\theta_k}^N}\right)r_k^N
    -\CL R_k^N.
\end{equation}

\subsection{Convergence to the overparameterized regime}
\label{sec:convergenceNinfty}

In this section, we investigate the discrete-time training dynamics \eqref{eq:NNfunction_evolution_prelimit} as the NN width $N$ tends to infinity.
For any fixed number $K$ of training steps, 
Theorem~\ref{thm:NTK_limit} below shows that the NN function trajectories~$\big(g_{\theta_{k}}^N\big)_{k=0,\dots,K}$, generated when training the NN~\eqref{eq:NN} with the regularized Newton and Gauss--Newton methods, converge in probability in $\CH^2(U)$ to the same deterministic limit, uniformly over these steps.

These limiting dynamics are governed by the deterministic \emph{limit PDE--NNTK integral operator} $T_{\CB^*_0}:L^2(U)\to\CH^2_{(0)}(U)\subset L^2(U)$,
defined by
\begin{equation}
    \label{eq:T_NNTK*}
    \left[T_{\CB^*_0}v\right]\!(x)
    \coloneqq
    \left[
        \frac{1}{\gamma}T_{B^*_0}
        \left(
            \left(
                \Id_{L^2(U)}
                +\frac{1}{\gamma}\CL T_{B^*_0}
            \right)^{-1}v
        \right)
    \right](x).
\end{equation}
Here,
$T_{B^*_0}:L^2(U)\to\CH^2_{(0)}(U)\subset L^2(U)$
is the \emph{limit PDE--NTK integral operator}
associated with gradient descent, given by
\begin{equation}
    \label{eq:T_NTK*}
    \left[T_{B^*_0}v\right]\!(x)
    \coloneqq
    \int_U B_0^*(x,x')v(x')\,dx',
\end{equation}
whose kernel $B_0^*:U\times U\to\bbR$ is determined by the initialization distribution $\mu_0$ through
\begin{equation}
    \label{eq:NTK*}
    B_0^*(x,x')
    \coloneqq
    \bbE_{\vartheta\sim\mu_0}\Big[
        \chi(x)s_\vartheta(x)
        \big(\CL\big(\chi(x') s_\vartheta(x')\big)\big)^\top
    \Big],
\end{equation}
where $s_\vartheta$ denotes the single-neuron Jacobian from \eqref{eq:s}.

The training dynamics of the corresponding deterministic sequence of limit NN functions $\left(g_k^*\right)_{k\geq0}$ is defined by
\begin{equation}
    \label{eq:NNfunction_evolution_limit}
    g_{k+1}^*
    =
    g_k^*
    +\alpha T_{\CB_0^*}\big(h-\CL g_k^*\big)
    \quad\text{ with }\quad
    g_0^*\equiv0.
\end{equation}
The associated \emph{limit residual evolution operator} is $\CL T_{\CB_0^*}:L^2(U)\to L^2(U)$.
Its spectral properties determine the convergence of the limiting PDE residual, as studied in Section~\ref{sec:globalconvergence}.
The following theorem quantifies the approximation of this sequence by the finite-width NN trajectories over any fixed number $K$ of training steps.
Its proof, including the convergence of the finite-width PDE--NNTK integral operators \eqref{eq:NNTK} to $T_{\CB_0^*}$, is given in Section~\ref{sec:proof:convergenceNinfty}.
The constants $C_\gamma$, $C_z$, and $C_w$ appearing in the statement are defined in the proof of Lemma~\ref{prop:NTK_limit_bounds_z}.

\begin{theorem}[Convergence of the finite-width training dynamics to their overparameterized limit over a fixed training horizon]
    \label{thm:NTK_limit}
    Set $\rho\coloneqq1-\beta-2/q>0$.
    Let $\delta\in(0,1/2]$ and $K\in\bbN$ be a given number of training steps.
    Fix $C_{R,0}>C_h$
    and let $C_\gamma$, $C_z$, and $C_w$ be the constants defined in \eqref{eq:def_constants}.
    Assume that the number of neurons $N$ is large enough such that
    \begin{equation}
        \label{eq:thm:asm:Nlargeenough}
        N^{2\beta-1}
        \geq
        \frac{3C(\sigma,\CD,\mu_0)}
             {\left(C_{R,0}-C_h\right)^2\delta}
        \quad\text{ and }\quad
        \gamma
        \geq
        \frac{
            \max\left\{C_\gamma,\alpha KC_zC_w\right\}
        }{N^\rho}.
    \end{equation}
    Then, we have that
    \begin{equation}
        \label{eq:thm:NTK_limit}
        \bbP\left(
            \max_{k=0,\dots,K}
            \Nbig{g_{\theta_k}^N-g_k^*}_{\CH^2(U)}
            \leq
            \frac C\delta\left(
                \frac{1}{N^{\beta-1/2}}
                +\frac{1}{N^{1-\beta}}
            \right)
        \right)
        \geq 1-2\delta
    \end{equation}
    for a constant
    $C=C(\alpha,\gamma,\beta,q, \sigma,\CD,\mu_0,C_c,C_R,C_z,K)>0$ independent of $N$.
    
    In particular, for fixed $\alpha$, $\gamma$, and $K$, the finite-width trajectory converges in probability to the deterministic limit trajectory in $\CH^2(U)$, uniformly over $k=0,\dots,K$. 
\end{theorem}

Let us note that the choice $\beta=3/4$ balances the two terms in \eqref{eq:thm:NTK_limit},
yielding a convergence rate of $N^{-1/4}$,
provided that Assumption~\ref{asm:NN_mu0}\ref{asm:NN_mu0iv} holds for some $q>8$.
A similar balance occurs in the asymptotic analysis of shallow NNs trained by stochastic gradient descent~\cite[Theorem~2.3]{yu2021normalization}.

\subsection{Global convergence in the overparameterized NN limit}
\label{sec:globalconvergence}

We study in what follows the long-time behavior of the deterministic limit dynamics $\left(g_k^*\right)_{k\geq0}$ given in \eqref{eq:NNfunction_evolution_limit}.
For a range of constant step sizes that includes the full step $\alpha=1$, we show that the corresponding PDE residual converges to zero in $L^2(U)$ as the number~$k$ of training iterations tends to infinity.
Consequently, the limiting loss
\begin{equation}
    \label{eq:loss_limit}
    L_k^*
    \coloneqq
    \frac{1}{2}\N{h-\CL g_k^*}_{L^2(U)}^2
\end{equation}
vanishes.

A central difficulty is that this setting does not guarantee a spectral gap for the limit residual evolution operator $\CL T_{\CB_0^*}$ governing the discrete-time training evolution of the PDE residual. 
To establish convergence, we first collect in the subsequent lemma,
whose proof we provide in Section~\ref{subsec:proof:spectralpropertiesTB},
the spectral properties of $\CL T_{B_0^*}$, the residual evolution operator associated with gradient descent, and discuss their implications for the residual evolution operator $\CL T_{\CB_0^*}$ of the regularized Newton and Gauss--Newton methods.
In particular, we show that every function in $\ker(\CL T_{B_0^*})$ is orthogonal to $\CL\big(\CH^2_{(0)}(U)\big)$.
This orthogonality property is essential for proving convergence of the PDE residual to zero in the overparameterized NN limit.

\begin{lemma}[Spectral properties of $\CL T_{B_0^*}$]
    \label{lem:spectralpropertiesTB}
    The integral operator $\CL T_{B_0^*}:L^2(U)\rightarrow L^2(U)$ is a self-adjoint and positive semi-definite Hilbert--Schmidt (hence, also compact) operator.
    In particular, there exists an orthonormal basis $\{e_\ell\}_{\ell\in I} \cup \{e^0_\ell\}_{\ell\in I_0}$ of $L^2(U)$
    such that
    \begin{equation}
        \CL T_{B_{0}^*} e_\ell = \lambda_\ell e_\ell
        \quad \text{ for all } \ell\in I,
        \quad \text{ and }\quad
        \CL T_{B_{0}^*} e^0_\ell = 0
        \quad \text{ for all } \ell\in I_0,
    \end{equation}
    where $I\subset\bbN$ and $I_0\subset\bbN$ are (at most) countable index sets,
    and where $\{\lambda_\ell\}_{\ell\in I}$ denotes the family of positive eigenvalues, counted with multiplicity and ordered in non-increasing order, i.e.,
    \begin{equation}
        \lambda_\ell\geq\lambda_k>0
        \quad\text{ for all }\ell,k\in I\text{ with }\ell<k.
    \end{equation}
    
    Moreover, if $e\in L^2(U)$ satisfies $\CL T_{B_0^*}e=0$, then $\langle e,\CL f\rangle_{L^2(U)}=0$ for every
    $f\in\CH_{(0)}^2$.
\end{lemma}

\begin{remark}[Spectral relation between $\CL T_{B_0^*}$ and $\CL T_{\CB_0^*}$]
    \label{rem:spectralrelationNNTK}
    The spectral properties of $\CL T_{B_0^*}$ from
    Lemma~\ref{lem:spectralpropertiesTB} directly determine those of
    the residual evolution operator $\CL T_{\CB_0^*}$. Indeed, recalling the definition of the limit PDE--NNTK
    integral operator from \eqref{eq:T_NNTK*}, we have
    \begin{equation}
        \CL T_{\CB_0^*} = \frac{1}{\gamma}\CL T_{B_0^*} \left( \Id_{L^2(U)} + \frac{1}{\gamma}\CL T_{B_0^*} \right)^{-1}.
    \end{equation}
    Consequently, $\CL T_{\CB_0^*}$ has the same eigenfunctions as
    $\CL T_{B_0^*}$. More precisely, for the orthonormal basis
    $\{e_\ell\}_{\ell\in I}\cup\{e_\ell^0\}_{\ell\in I_0}$ from
    Lemma~\ref{lem:spectralpropertiesTB}, it holds
    \begin{equation}
        \CL T_{\CB_0^*}e_\ell
        = \frac{\lambda_\ell}{\gamma+\lambda_\ell} e_\ell
        \quad \text{ for all } \ell\in I,
        \quad \text{ and }\quad
        \CL T_{\CB_0^*}e_\ell^0 = 0
        \quad \text{ for all } \ell\in I_0,
    \end{equation}
    i.e., $\ker\left(\CL T_{\CB_0^*}\right) = \ker\left(\CL T_{B_0^*}\right)$.
    
    For notational simplicity, we write $\widetilde\lambda_\ell\coloneqq\lambda_\ell/(\gamma+\lambda_\ell)$.
    Since the map $\lambda\mapsto\widetilde\lambda=\lambda/(\gamma+\lambda)$ is strictly increasing on $[0,\infty)$, the ordering of the positive eigenvalues is preserved. 
    Defining $\lambda_{\max}\coloneqq\|\CL T_{B_0^*}\|_{\mathrm{op}}$ and $\widetilde\lambda_{\max} \coloneqq\|\CL T_{\CB_0^*}\|_{\mathrm{op}}$,
    we obtain $\widetilde\lambda_{\max} =\lambda_{\max}/(\gamma+\lambda_{\max})<1$.
\end{remark}

The following statements present the main convergence results of this paper.
We provide their proofs in Sections \ref{subsec:proof:globalconvergence} and \ref{subsec:proof:finite_mode_rate}, respectively.
The first, result shows that the PDE residual associated with the deterministic infinite-width dynamics \eqref{eq:NNfunction_evolution_limit} of both the regularized Newton and Gauss--Newton methods converges to zero in $L^2(U)$ as the number of training iterations~$k$ tends to infinity.
Consequently, the limiting loss $L_k^*$ vanishes as $k\rightarrow\infty$, establishing global convergence of both the regularized Newton and Gauss--Newton methods in the overparameterized NN limit.
In the second result, we quantify this convergence.

\begin{theorem}[Global convergence to zero loss
    in the infinite-width limit]
    \label{thm:convergence}
    Let the learning rate satisfy $\alpha>0$ and $\alpha\widetilde\lambda_{\max}<2$.
    Then the PDE residual $r_k^*\coloneqq h-\CL g_k^*$ of the limiting NN function $g_k^*$ governed by \eqref{eq:NNfunction_evolution_limit} converges to zero in $L^2(U)$
    as $k\to\infty$, i.e.,
    \begin{equation}
        \N{r_k^*}_{L^2(U)}
        =
        \N{h-\CL g_k^*}_{L^2(U)}
        \rightarrow0
        \quad\text{ as }k\to\infty.
    \end{equation}
    Consequently, the loss $L_k^*$ defined in
    \eqref{eq:loss_limit} vanishes as $k\to\infty$.
\end{theorem}

Theorem~\ref{thm:convergence} establishes convergence of the full PDE residual.
Under an additional stability assumption on the PDE operator, this also implies convergence to the unique PDE solution.

\begin{corollary}[Convergence to the PDE solution in the infinite-width limit]
    \label{cor:solution_convergence}
    Under the assumptions of Theorem~\ref{thm:convergence},
    assume additionally that there exists
    $C>0$ such that
    \begin{equation}
        \label{eq:PDE_stability}
        \N{v}_{\CH^1(U)}
        \leq C\N{\CL v}_{L^2(U)}
        \quad\text{ for every }v\in\CH^2_{(0)}(U).
    \end{equation}
    Then the limiting NN function $g_k^*$ governed by \eqref{eq:NNfunction_evolution_limit} satisfies
    \begin{equation}
        \label{eq:solution_convergence}
        \N{g_k^*-u}_{\CH^1(U)}
        \leq C\N{r_k^*}_{L^2(U)}
        \rightarrow0
        \quad\text{ as }k\to\infty.
    \end{equation}

    Moreover, if \eqref{eq:PDE_stability} holds with the $\CH^2(U)$ norm in place of the $\CH^1(U)$ norm, then the convergence in \eqref{eq:solution_convergence} also holds in $\CH^2(U)$.
\end{corollary}

\begin{remark}
    \label{rem:PDE_stability}
    Estimate \eqref{eq:PDE_stability} in Corollary~\ref{cor:solution_convergence} follows, for instance,
    if the bilinear form $a$ associated with $\CL$ is bounded and coercive on $\CH^1_0(U)$.
    Indeed, for some $\kappa>0$ and every $v\in\CH^2_{(0)}(U)$,
    \begin{equation*}
        \kappa\N{v}_{\CH^1(U)}^2 \leq a(v,v) = \langle\CL v,v\rangle_{L^2(U)} \leq \N{\CL v}_{L^2(U)}\N{v}_{\CH^1(U)}
    \end{equation*}
    which yields \eqref{eq:PDE_stability}.
    More general conditions can be found in \cite[Chapter~8]{gilbarg1977elliptic}.
    Coercivity is sufficient but not necessary.
    A standard setting in which \eqref{eq:PDE_stability}
    holds is the uniformly elliptic Dirichlet problem for
    \[
        \CL v=-\nabla\cdot(A\nabla v)+b\cdot\nabla v+cv,
    \]
    where $A\in W^{1,\infty}(U;\bbR^{d\times d})$
    is symmetric and uniformly elliptic,
    $b\in L^\infty(U;\bbR^d)$,
    and $c\in L^\infty(U)$ is nonnegative.
    On a bounded Lipschitz domain, weak Dirichlet solvability
    and boundedness of the solution operator give
    \[
        \N{v}_{\CH^1(U)}
        \leq C\N{\CL v}_{L^2(U)}
        \quad\text{ for }v\in \CH^1_0(U)
        \text{ with }\CL v\in L^2(U),
    \]
    see \cite[Chapter~8]{gilbarg1977elliptic}.

    For the corresponding $\CH^2(U)$ stability estimate,
    additional regularity is needed.
    For example, if, in addition, $A\in C^1(\overline U;\bbR^{d\times d})$
    and $\partial U$ is of class $C^2$, global elliptic
    regularity gives
    \[
        \N{v}_{\CH^2(U)}
        \leq C
        \left(
            \N{\CL v}_{L^2(U)}
            +\N{v}_{L^2(U)}
        \right)
        \quad\text{ for every }v\in\CH^2_{(0)}(U),
    \]
    see \cite[Section~6.3.2, Theorem~4]{evans2010partial}.
    Combining this estimate with \eqref{eq:PDE_stability}
    yields $\N{v}_{\CH^2(U)}\leq C\N{\CL v}_{L^2(U)}$.
\end{remark}

For a finite spectral projection, an explicit geometric rate can be obtained from the largest and smallest eigenvalues included in that projection.
The following result quantifies this convergence.

\begin{theorem}[Geometric convergence on finitely many modes]
    \label{lem:finite_mode_rate}
   Let the learning rate satisfy $\alpha>0$ and $\alpha\widetilde\lambda_{\max}<2$.
    Let $L\in\bbN$ be such that $\CL T_{\CB_0^*}$
    has at least $L$ positive eigenvalues, and let
    $e_1,\dots,e_L$ be orthonormal eigenfunctions
    associated with its $L$ largest positive eigenvalues,
    ordered as
    \begin{equation}
        \label{eq:eigenvalue_ordering}
        \widetilde\lambda_1
        \geq\widetilde\lambda_2
        \geq\dots
        \geq\widetilde\lambda_L>0.
    \end{equation}
    Denote by
    $r^*_{k,L}\coloneqq
    \sum_{\ell=1}^L
    \langle r_k^*,e_\ell\rangle_{L^2(U)}e_\ell$
    the projection of the PDE residual onto
    $\mathrm{span}\{e_1,\dots,e_L\}$.
    Then, for every $k\in\bbN_0$, it holds
    \begin{equation}
        \label{eq:finite_mode_rate}
        \Nbig{r^*_{k,L}}_{L^2(U)}
        \leq
        q_L^k\Nbig{r^*_{0,L}}_{L^2(U)},
    \quad\text{ with }\quad
        q_L
        \coloneqq
        \max\left\{
            \abs{1-\alpha\widetilde\lambda_1},
            \abs{1-\alpha\widetilde\lambda_L}
        \right\}
        <1.
    \end{equation}
    In particular, for the full step $\alpha=1$,
    \begin{equation}
        \label{eq:finite_mode_rate_fullstep}
        \Nbig{r^*_{k,L}}_{L^2(U)}
        \leq
        \left(\frac{\gamma}{\gamma+\lambda_L}\right)^k
        \Nbig{r^*_{0,L}}_{L^2(U)}.
    \end{equation}
\end{theorem}

The spectral transformation $\lambda_\ell\mapsto\lambda_\ell/(\gamma+\lambda_\ell)$ explains the acceleration relative to gradient descent.
For gradient descent with step size $\alpha$, each spectral component is multiplied by $1-\alpha\lambda_\ell$,
so components associated with small eigenvalues can decay much more slowly than those associated with large eigenvalues.
For the regularized Newton and Gauss--Newton methods with $\alpha=1$, the corresponding factor is $\gamma/(\gamma+\lambda_\ell)$.
Thus, when $\gamma$ is small compared with $\lambda_L$, all components in the finite spectral projection decay rapidly, with a common contraction bound $\gamma/(\gamma+\lambda_L)$.
This reduces the disparity between convergence rates across these modes, while modes with eigenvalues below the regularization scale may still converge slowly. 

\begin{remark}
    \label{rem:finite_mode_scope}
    Theorem~\ref{lem:finite_mode_rate} gives a geometric rate for the projection $r^*_{k,L}$ onto the first $L$ modes.
    Its contraction factor $q_L$ approaches $1$ as $L\to\infty$ whenever $\widetilde\lambda_L\to0$, so the estimate does not provide a uniform geometric rate across all positive-eigenvalue modes.
    If the initial residual is supported on finitely many such modes, the theorem also gives a geometric rate for the full residual.
\end{remark}

\section{Proof Details for Theorem~\ref{thm:NTK_limit}}
\label{sec:proof:convergenceNinfty}

This section comprises the proof details for Theorem~\ref{thm:NTK_limit}, proving that the training trajectory $(g_{\theta_k}^N(\x))_{k=0,\dots,K}$ of the finite-width NN function $g_\theta^N(\x)$, when trained using the (Gauss--)Newton method \eqref{eq:regularizedNewton}--\eqref{eq:regularizedNewtonStep}, as captured by \eqref{eq:NNfunction_evolution_prelimit}, converges in probability as the number of neurons $N$ tends to infinity to the infinite-width training trajectory $(g_k^\ast(\x))_{k=0,\dots,K}$ given by \eqref{eq:NNfunction_evolution_limit}.

\subsection{Auxiliary estimates}
\label{subsec:proof:auxtechnicalest}

We first collect several auxiliary estimates.
Throughout this section, $\theta=(\theta^i)_{i=1}^N\in\bbR^{N(d+2)}$ denotes a generic vector of NN parameters, while $\vartheta=(c,w,\eta)\in\bbR^{d+2}$ denotes a generic single-neuron parameter.

\begin{lemma}[Regularity of activation function~$\sigma$]
    \label{lem:sigma_consequences}
    We write $\sigma^{(0)}\coloneqq\sigma$.
    Then the following statements hold:
    \begin{enumerate}[label=(\roman*),labelsep=10pt,leftmargin=35pt]
        \item\label{lem:sigma:bounds} The quantities $C_{\sigma^{(j)}}\coloneqq\Nbig{\sigma^{(j)}}_{L^\infty(\bbR)}$, $j=0,1,2,3,4$, are finite.
        For notational convenience, we write $C_\sigma\coloneqq C_{\sigma^{(0)}}$, $C_{\sigma'}\coloneqq C_{\sigma^{(1)}}$, $C_{\sigma''}\coloneqq C_{\sigma^{(2)}}$, $C_{\sigma'''}\coloneqq C_{\sigma^{(3)}}$, $C_{\sigma''''}\coloneqq C_{\sigma^{(4)}}$, and define $C^{\max}_\sigma\coloneqq\max_{0\leq j\leq4}C_{\sigma^{(j)}}<\infty$.
        \item\label{lem:sigma:lipschitz} For $j=0,1,2,3$, the function $\sigma^{(j)}$ is globally Lipschitz continuous with Lipschitz constant $C_{\sigma^{(j+1)}}$.
        In particular, the Lipschitz constant in Assumption~\ref{asm:NN_sigma} may be chosen as $L_\sigma=C_{\sigma'}$. We write $L_{\sigma'}\coloneqq C_{\sigma''}$, $L_{\sigma''}\coloneqq C_{\sigma'''}$, and $L_{\sigma'''}\coloneqq C_{\sigma''''}$.
        \item\label{lem:sigma:continuity} For $j=0,1,2,3,4$, the function $\sigma^{(j)}$ is continuous, and consequently the map $(\vartheta,x) =  \big((c,w,\eta),x\big)\mapsto\sigma^{(j)}(w\cdot x+\eta)$ is continuous on $\bbR^{d+2}\times\bbR^d$.
        \item\label{lem:sigma:C4} For every $x\in\bbR^d$, the map $\vartheta=(c,w,\eta)\mapsto c\,\sigma(w\cdot x+\eta)$ belongs to $\CC^4(\bbR^{d+2})$. Every parameter derivative of order $m\leq4$ is either zero or of the form $r(c)p(x)\sigma^{(j)}(w\cdot x+\eta)$, where $\deg r\leq1$, $j\leq m$, and $p$ is a product of at most $m$ coordinate functions. These derivatives are jointly continuous in $(\vartheta,x)$ and, for $\abs{c}\leq\widetilde C$ and $x\in U$, are bounded in absolute value by $\max\{1,\widetilde C\}\max\{1,C_U^4\}C^{\max}_\sigma$.
    \end{enumerate}
\end{lemma}

\begin{proof}
    Parts~\ref{lem:sigma:bounds} and~\ref{lem:sigma:continuity} follow directly from Assumptions~\ref{asm:NN_sigma}--\ref{asm:NN_sigma'} and continuity of $(\vartheta,x)\mapsto w\cdot x+\eta$.
    For $j\leq3$, the mean value theorem gives
    \begin{equation}
        \absbig{\sigma^{(j)}(a)-\sigma^{(j)}(b)}
        \leq C_{\sigma^{(j+1)}}\abs{a-b},
        \qquad a,b\in\bbR,
    \end{equation}
    which proves Part~\ref{lem:sigma:lipschitz}.
    Finally, Part~\ref{lem:sigma:C4} follows by repeated application of the chain and product rules: since the dependence on $c$ is linear, differentiating once with respect to $c$ removes the factor $c$, while any derivative involving two or more differentiations with respect to $c$ is zero. Each differentiation with respect to a component of $w$ produces the corresponding component of $x$ and raises the order of the derivative of $\sigma$ by one, whereas differentiation with respect to $\eta$ only raises that order. Joint continuity follows from Part~\ref{lem:sigma:continuity}, and the stated bound follows from Part~\ref{lem:sigma:bounds} and $\sup_{x\in U}\N{x}_2\leq C_U$ in Assumption~\ref{asm:D}.
\end{proof}

Recalling the notation of Section~\ref{sec:evolution},
$J^N_{\theta^i}=\frac{1}{N^{\beta}}s_{\theta^i}\in\bbR^{1\times(d+2)}$ and
$H^N_{\theta^i}=\frac{1}{N^{\beta}}h_{\theta^i}\in\bbR^{(d+2)\times(d+2)}$.

Since parameter derivatives act on $\vartheta$, whereas $D_\nu$ and $\CL$ act on the spatial variable, and since the coefficients of $\CL$ are independent of $\vartheta$ and the required mixed derivatives exist by Lemma~\ref{lem:sigma_consequences},
these operations commute.
In particular,
\begin{equation}
    \label{eq:commutation_CL_parameter_gradient}
    \Big(\nabla_\vartheta\left(\CL\phi_\vartheta\right)\!\Big)^\top
    =\Big(\nabla_\vartheta\left(\CL\left(\chi\,c\,\sigma(w\cdot\dummy + \eta)\right)\right)\!\Big)^\top
    =\CL\left(\chi\,s_{\vartheta}\right)
\end{equation}
and
\begin{equation}
    \nabla_\vartheta^2\left(\CL\phi_\vartheta\right)
    =\nabla_\vartheta^2\left(\CL\left(\chi\,c\,\sigma(w\cdot\dummy + \eta)\right)\right)
    =\CL\left(\chi\,h_{\vartheta}\right).
\end{equation}
Here and in what follows, we understand $D_\nu$ and $\CL$ to act component-wise on vector- and matrix-valued functions.

\begin{lemma}[Bounds on the single-neuron function~$\phi_\vartheta$ and its parameter derivatives]
    \label{lem:bound_s}
    \label{lem:bound_h}
    Let $\vartheta=(c,w,\eta)\in\bbR^{d+2}$ and assume that $\abs{c}\leq\widetilde{C}$.
    Then, for every multi-index $\abs{\nu}\leq2$, the following statements hold for a constant $C=C(d,C_U,\N{\chi}_{\CC^2_b},k_\CL,C_\sigma^{\max},\widetilde C)>0$:
    \begin{enumerate}[label=(\roman*),labelsep=10pt,leftmargin=35pt]
        \item\label{lem:bound_s:phi} The single-neuron function satisfies
        \begin{equation}
        \begin{aligned}
            \sup_{x\in U}\absbig{\left[D_\nu\phi_\vartheta\right]\!(x)}
            &\leq C\left(1\!+\!\N{w}_2^2\right)
            \quad\text{and}\quad
            \operatorname*{ess\,sup}_{x\in U}\absbig{\left[\CL\phi_\vartheta\right]\!(x)}
            \leq C\left(1\!+\!\N{w}_2^2\right).
        \end{aligned}
        \end{equation}
        \item\label{lem:bound_s:jacobian} The parameter Jacobian $\chi s_{\vartheta}$ satisfies
        \begin{equation}
            \label{eq:lem:bound_s}
        \begin{aligned}
            \sup_{x\in U}\Nbig{\left[D_\nu\!\left(\chi\,s_{\vartheta}\right)\right]\!(x)}_2&\leq C\left(1\!+\!\N{w}_2^2\right)
            \quad\text{and}\quad
            \operatorname*{ess\,sup}_{x\in U}\Nbig{\left[\CL\!\left(\chi\,s_{\vartheta}\right)\right]\!(x)}_2\leq C\left(1\!+\!\N{w}_2^2\right).
        \end{aligned}
        \end{equation}
        Moreover, $(\vartheta,x)\mapsto \left[D_\nu(\chi s_\vartheta)\right]\!(x)$ is continuous on $\bbR^{d+2}\times\overline U$.
        \item\label{lem:bound_s:hessian} The parameter Hessian $\chi h_{\vartheta}$ satisfies
        \begin{equation}
        \begin{aligned}
            \sup_{x\in U}\Nbig{\left[D_\nu\!\left(\chi\,h_\vartheta\right)\right]\!(x)}_{\mathrm{op}}&\!\leq C\left(1\!+\!\N{w}_2^2\right)
            \quad\text{and}\quad
            \operatorname*{ess\,sup}_{x\in U}\Nbig{\left[\CL\!\left(\chi\,h_\vartheta\right)\right]\!(x)}_{\mathrm{op}}\!\leq C\left(1\!+\!\N{w}_2^2\right).
        \end{aligned}
        \end{equation}
        Moreover, $(\vartheta,x)\mapsto \left[D_\nu(\chi h_\vartheta)\right]\!(x)$ is continuous on $\bbR^{d+2}\times\overline U$.
    \end{enumerate}
    Here the essential suprema are taken with respect to the Lebesgue measure.
    The ordinary suprema in the estimates involving $D_\nu$ refer to the continuous representatives of these derivatives.
    
    If $\chi$ is replaced by the constant function~$1$, so that
    $\widetilde\phi_\vartheta(x)=c\sigma(w\cdot x+\eta)$, Parts~\ref{lem:bound_s:phi}--\ref{lem:bound_s:hessian} remain valid.
\end{lemma}

\begin{proof}
    The scalar function $\phi_\vartheta$ and every component of $\chi s_\vartheta$ as well as every nonzero component of $\chi h_\vartheta$ are of the form $x\mapsto\chi(x)r(c)p(x)\sigma^{(j)}(w\cdot x+\eta)$,
    where $\deg r\leq1$, $j\leq2$, and $p$ is a product of at most two coordinate functions.
    Indeed, for $\phi_\vartheta$ we have $r(c)=c$, $p=1$, and $j=0$, while the possible choices for the components of $\chi s_\vartheta$ and $\chi h_\vartheta$ follow directly from \eqref{eq:s} and \eqref{eq:h}.

    Let $\nu$ be a multi-index with $\abs{\nu}\leq2$.
    By the Leibniz and chain rules, $D_\nu\phi_\vartheta$ and every component of $D_\nu(\chi s_\vartheta)$ and $D_\nu(\chi h_\vartheta)$ are finite sums of terms of the form
    \begin{equation}
        \label{eq:proof:lem:bound_s:generic}
        \big(D_{\nu_1}\chi\big)(x)\,r(c)\,\big(D_{\nu_2}p\big)(x)\,q(w)\,\sigma^{(j')}(w\cdot x+\eta),
    \end{equation}
    where $q$ is a monomial in the components of $w$,
    $\abs{\nu_1}+\abs{\nu_2}+\deg q\leq2$, and $j'=j+\deg q\leq4$.
    Since $\abs{c}\leq\widetilde C$, $\sup_{x\in U}\N{x}_2\leq C_U$, and $\sigma^{(j')}$ and the derivatives of $\chi$ up to order two are bounded, every such term is bounded uniformly in $x\in U$ by a constant times $1+\N{w}_2^2$.
    Consequently, Lemma~\ref{lem:sigma_consequences} and Assumptions~\ref{asm:D} and~\ref{asm:NN_chi} yield
    \begin{equation}
        \sup_{x\in U}\absbig{\left[D_\nu\phi_\vartheta\right]\!(x)}
        \leq C\left(1+\N{w}_2^2\right)
    \end{equation}
    and, after absorbing the finite number $d+2$ of components into $C$,
    \begin{equation}
        \sup_{x\in U}\Nbig{\left[D_\nu(\chi s_\vartheta)\right]\!(x)}_2
        \leq C\left(1+\N{w}_2^2\right)
    \end{equation}
    and, using that the operator norm of a $(d+2)\times(d+2)$ matrix is bounded by $(d+2)$ times the largest absolute value of its entries,
    \begin{equation}
        \sup_{x\in U}\Nbig{\left[D_\nu(\chi h_\vartheta)\right]\!(x)}_{\mathrm{op}}
        \leq C\left(1+\N{w}_2^2\right)
    \end{equation}
    for a constant $C=C(d,C_U,\N{\chi}_{\CC^2_b},C_\sigma^{\max},\widetilde C)>0$.

    For Part~\ref{lem:bound_s:phi}, applying Assumption~\ref{asm:CL} to $\phi_\vartheta$ and taking the essential supremum gives
    \begin{equation}
        \operatorname*{ess\,sup}_{x\in U}\absbig{\left[\CL\phi_\vartheta\right]\!(x)}
        \leq k_\CL\sum_{0\leq\abs{\nu}\leq2}
        \sup_{x\in U}\absbig{\left[D_\nu\phi_\vartheta\right]\!(x)}
        \leq C\left(1+\N{w}_2^2\right),
    \end{equation}
    where we made use of Remark~\ref{rem:CL_linearity} in the first step.

    For Part~\ref{lem:bound_s:jacobian}, applying the same assumption component-wise to $\chi s_\vartheta$ and using the triangle inequality shows
    \begin{equation}
        \begin{split}
        \operatorname*{ess\,sup}_{x\in U}\Nbig{\left[\CL(\chi s_\vartheta)\right]\!(x)}_2
        &\leq k_\CL\sum_{0\leq\abs{\nu}\leq2}
        \sup_{x\in U}\Nbig{\left[D_\nu(\chi s_\vartheta)\right]\!(x)}_2
        \leq C\left(1+\N{w}_2^2\right).
        \end{split}
    \end{equation}
    Every term in \eqref{eq:proof:lem:bound_s:generic} is jointly continuous in $(\vartheta,x)$ by Lemma~\ref{lem:sigma_consequences}\ref{lem:sigma:continuity} and the regularity of $\chi$.
    This proves the continuity statement in Part~\ref{lem:bound_s:jacobian}.

    For Part~\ref{lem:bound_s:hessian}, the same argument applied component-wise to $\chi h_\vartheta$, together with the preceding matrix-norm bound, yields
    \begin{equation}
        \operatorname*{ess\,sup}_{x\in U}\Nbig{\left[\CL(\chi h_\vartheta)\right]\!(x)}_{\mathrm{op}}
        \leq C\left(1+\N{w}_2^2\right).
    \end{equation}
    The joint continuity of the terms in \eqref{eq:proof:lem:bound_s:generic} also proves the continuity statement in Part~\ref{lem:bound_s:hessian}.
    After potentially increasing the constant $C$ to absorb $k_\CL$ and the sum over a finite number of multi-indices, all three parts follow.
    
    Finally, replacing $\chi$ by the constant function~$1$ removes all terms containing a nonzero derivative of $\chi$.
    Hence, the same estimates and continuity argument remain valid, potentially after increasing $C$.
\end{proof}

\begin{lemma}[Lipschitz continuity of the derivatives of the Jacobian~$s_\vartheta$ in the parameters]
    \label{lem:bound_s_Lipschitz}
    Let $\vartheta=(c,w,\eta), \vartheta'=(c',w',\eta')\in\bbR^{d+2}$ and assume that $\abs{c},\abs{c'}\leq\widetilde C$.
    Then, for every multi-index $\abs{\nu}\leq2$, it holds
    \begin{equation}
        \sup_{x\in U}\Nbig{\left[D_\nu\!\left(\chi\,s_\vartheta\right)\right]\!(x)-\left[D_\nu\!\left(\chi\,s_{\vartheta'}\right)\right]\!(x)}_2
        \leq C\left(1+\N{w}_2^2+\N{w'}_2^2\right)\Nbig{\vartheta-\vartheta'}_2
    \end{equation}
    and
    \begin{equation}
        \operatorname*{ess\,sup}_{x\in U}\Nbig{\left[\CL\!\left(\chi\,s_\vartheta\right)\right]\!(x)-\left[\CL\!\left(\chi\,s_{\vartheta'}\right)\right]\!(x)}_2
        \leq C\left(1+\N{w}_2^2+\N{w'}_2^2\right)\Nbig{\vartheta-\vartheta'}_2,
    \end{equation}
    for a constant $C=C(d,C_U,\N{\chi}_{\CC^2_b},k_\CL,C_\sigma^{\max},\widetilde C)>0$.
    
    The first estimate also remains valid with $\chi$ replaced by the constant function~$1$.
\end{lemma}

\begin{proof}
    Differentiating a term of the form \eqref{eq:proof:lem:bound_s:generic} arising from $D_\nu(\chi s_\vartheta)$ once more with respect to $\vartheta=(c,w,\eta)$ produces terms of the same type:
    $\partial_c$ differentiates the factor $r(c)$ and leaves $j'$ unchanged;
    $\partial_\eta$ raises $j'$ by one;
    and $\partial_{w_r}$ either differentiates the monomial $q(w)$ or raises $j'$ by one while contributing an additional factor $x_r$.
    For these terms, $j\leq1$ and hence $j'\leq3$ before the additional parameter derivative is taken.
    Thus the resulting terms contain a monomial in $w$ of degree at most two and a derivative of $\sigma$ of order at most four, which is covered by Lemma~\ref{lem:sigma_consequences}\ref{lem:sigma:bounds}.
    For fixed $x\in U$, set $F_{\nu,x}(\xi)\coloneqq D_\nu(\chi s_\xi)(x)$ as well as $\vartheta_t\coloneqq\vartheta'+t(\vartheta-\vartheta')$ and analogously $c_t\coloneqq c'+t(c-c')$ and $w_t\coloneqq w'+t(w-w')$, $t\in[0,1]$.
    The preceding bounds give
    \begin{equation}
        \sup_{x\in U}\Nbig{D F_{\nu,x}(\vartheta_t)}_{\mathrm{op}}
        \leq C\left(1+\N{w}_2^2+\N{w'}_2^2\right).
    \end{equation}
    Indeed, $\abs{c_t}\leq\widetilde C$ and $\N{w_t}_2^2\leq\N{w}_2^2+\N{w'}_2^2$ by convexity.
    The fundamental theorem of calculus along this segment therefore yields the first estimate.

    Applying Assumption~\ref{asm:CL} component-wise and using the first estimate gives
    \begin{equation}
    \begin{split}
        &\operatorname*{ess\,sup}_{x\in U}\Nbig{\left[\CL(\chi s_\vartheta)\right]\!(x)-\left[\CL(\chi s_{\vartheta'})\right]\!(x)}_2\\
        &\qquad\leq k_\CL\sum_{0\leq\abs{\nu}\leq2}
        \sup_{x\in U}\Nbig{\left[D_\nu(\chi s_\vartheta)\right]\!(x)-\left[D_\nu(\chi s_{\vartheta'})\right]\!(x)}_2,
    \end{split}
    \end{equation}
    and hence the second estimate follows after potentially increasing $C$.
    
    Repeating the argument with $\chi$ replaced by the constant function~$1$ proves the final assertion.
\end{proof}

\begin{lemma}
    \label{lem:estimate_S_operatornorm_i}
    Let $\theta^i=(c^i,w^i,\eta^i)$ and assume that $\absnormal{c^i}\leq\widetilde{C}$.
    Then it holds
    \begin{equation}
        \label{eq:lem:estimate_S_operatornorm_i}
        \N{\int_U\left(h(x)-\CL g^N_\theta(x)\right)\left[\CL\!\left(\chi\,h_{\theta^i}\right)\right]\!(x)\,dx}_{\mathrm{op}}
        \leq C_S\left(1+\N{w^i}_2^2\right)\N{h-\CL g^N_\theta}_{L^2(U)}
    \end{equation}
    for a constant $C_S\coloneqq C\sqrt{\vol{U}}$, where $C=C(d,C_U,\N{\chi}_{\CC^2_b},k_\CL,C_\sigma^{\max},\widetilde C)>0$ is the constant from Lemma~\ref{lem:bound_h}.
\end{lemma}

\begin{proof}
    We compute with Jensen's inequality in the first and Cauchy--Schwarz inequality in the second step that
    \begin{equation}
    \begin{split}
        &\N{\int_U\left(h(x)-\CL g^N_\theta(x)\right)\left[\CL\!\left(\chi\,h_{\theta^i}\right)\right]\!(x)\,dx}_{\mathrm{op}}
        \leq\int_U\abs{h(x)-\CL g^N_\theta(x)}\Nbig{\left[\CL\!\left(\chi\,h_{\theta^i}\right)\right]\!(x)}_{\mathrm{op}}dx \\
        &\qquad\quad\,\leq\N{h-\CL g^N_\theta}_{L^2(U)}\left(\int_U\Nbig{\left[\CL\!\left(\chi\,h_{\theta^i}\right)\right]\!(x)}^2_{\mathrm{op}}dx\right)^{1/2}\\
        &\qquad\quad\,\leq\sqrt{\vol{U}}\,C\left(1+\N{w^i}_2^2\right)\N{h-\CL g^N_\theta}_{L^2(U)},
    \end{split}
    \end{equation}
    where the last step follows from Lemma~\ref{lem:bound_h} and $\vol{U}<\infty$ in Assumption~\ref{asm:D}.
    Choosing $C_S\coloneqq\sqrt{\vol{U}}\,C$ proves \eqref{eq:lem:estimate_S_operatornorm_i}.
\end{proof}

\subsection{Invertibility of the regularized matrix~\texorpdfstring{$D^N_\theta$}{D}}
\label{subsec:invertibility_D}

The invertibility of $D^N_\theta$ is essential for the well-posedness of the regularized (Gauss--)Newton step and for the operator formulation developed in the next subsection.
For the Gauss--Newton method, it is immediate from $D^N_\theta=\gamma^N\Id_{N(d+2)}$, whereas for the Newton method the additional, potentially indefinite nonconvex contribution $S^N_\theta$ makes invertibility nontrivial.
The following lemma establishes the required uniform invertibility result for both methods.

For $\theta=(\theta^i)_{i=1}^N$, let $D^N_{\theta^i}$ denote the $i$-th diagonal block of $D^N_\theta$ and define its rescaled version by
\begin{equation}
    \label{eq:def_d}
\begin{split}
    d^N_{\theta^i}
    &\coloneqq N^{2\beta-1}D^N_{\theta^i} \\
    &=\begin{cases}
        \gamma\Id_{d+2},
        & \quad\text{ for the Gauss--Newton method},\\[1mm]
        \displaystyle
        \gamma\Id_{d+2}-\frac{1}{N^{1-\beta}}
        \int_U\left(h(x)-\CL g^N_\theta(x)\right)
        \left[\CL\left(\chi\,h_{\theta^i}\right)\right]\!(x)\,dx,
        & \quad\text{ for the Newton method}.
    \end{cases}
\end{split}
\end{equation}
The matrices $d^N_{\theta^i}$ are symmetric: this is immediate for Gauss--Newton, while for Newton it follows because $h_{\theta^i}(x)$ is a Hessian and $\CL$ acts entry-wise.
Consequently, the block-diagonal matrix $D^N_\theta$ is symmetric as well.

\begin{lemma}[Uniform invertibility of~$D^N_\theta$]
    \label{lem:estimate_Dinv_operatornorm}
    Let $\theta=(\theta^i)_{i=1}^N$ with $\theta^i=(c^i,w^i,\eta^i)$.
    For the \emph{Newton method}, assume that $\abs{c^i}\leq\widetilde C$ for all $i=1,\dots,N$.
    Assume moreover that
    \begin{equation}
        \label{eq:cond_D}
        \gamma\geq\frac{2C_S\max_{i=1,\dots,N}\left(1+\N{w^i}_2^2\right)}{N^{1-\beta}}\Nbig{h-\CL g^N_\theta}_{L^2(U)}
    \end{equation}
    with $C_S>0$ as in Lemma~\ref{lem:estimate_S_operatornorm_i}.
    Then $D^N_\theta$ and all $d^N_{\theta^i}$ are symmetric positive definite and
    \begin{equation}
        \label{eq:lem:Dinv}
        \N{\big(D^N_\theta\big)^{-1}}_{\mathrm{op}}\leq\frac2\gamma N^{2\beta-1}
        \quad\text{ and }\quad
        \gamma\N{\big(d^N_{\theta^i}\big)^{-1}}_{\mathrm{op}}\leq2\quad\text{ for all }i=1,\dots,N.
    \end{equation}
    
    For the \emph{Gauss--Newton method}, the same conclusions hold for every $\theta$ and every $\gamma>0$, and \eqref{eq:lem:Dinv} holds with the constant $2$ replaced by $1$.
\end{lemma}

\begin{proof}
    For the Newton method, $S^N_\theta$ is block-diagonal and symmetric, and its $i$-th block is
    \begin{equation}
        S^N_{\theta^i}=-\frac{1}{N^\beta}\int_U\left(h(x)-\CL g^N_\theta(x)\right)\left[\CL\!\left(\chi\,h_{\theta^i}\right)\right]\!(x)\,dx.
    \end{equation}
    Its operator norm is therefore the maximum of the operator norms of its blocks.
    Lemma~\ref{lem:estimate_S_operatornorm_i} and condition~\eqref{eq:cond_D} therefore give
    \begin{equation}
        \Nbig{S^N_\theta}_{\mathrm{op}}
        \leq\frac{C_S\max_{i=1,\dots,N}\big(1+\N{w^i}_2^2\big)}{N^{\beta}}\Nbig{h-\CL g^N_\theta}_{L^2(U)}
        \leq\frac{\gamma}{2N^{2\beta-1}}=\frac{\gamma^N}{2}.
    \end{equation}
    Hence
    \begin{equation}
        \lambda_{\min}(D^N_\theta)
        =\gamma^N+\lambda_{\min}(S^N_\theta)
        \geq\gamma^N-\N{S^N_\theta}_{\mathrm{op}}
        \geq\frac{\gamma^N}{2}>0.
    \end{equation}
    Thus $D^N_\theta$ is positive definite and
    $\N{(D^N_\theta)^{-1}}_{\mathrm{op}}\leq2/\gamma^N=(2/\gamma)N^{2\beta-1}$.
    Each block $D^N_{\theta^i}$, and hence each $d^N_{\theta^i}=N^{2\beta-1}D^N_{\theta^i}$, is positive definite as well. Moreover,
    \begin{equation}
        \N{(d^N_{\theta^i})^{-1}}_{\mathrm{op}}
        =\frac{1}{N^{2\beta-1}}\N{(D^N_{\theta^i})^{-1}}_{\mathrm{op}}
        \leq \frac{1}{N^{2\beta-1}}\N{(D^N_\theta)^{-1}}_{\mathrm{op}}
        \leq\frac2\gamma,
    \end{equation}
    which proves the second bound in \eqref{eq:lem:Dinv}.

    For the Gauss--Newton method, both bounds follow directly from $D^N_\theta=\gamma^N\Id_{N(d+2)}$ and $d^N_{\theta^i}=\gamma\Id_{d+2}$.
\end{proof}

\subsection{Operator formulation of the (Gauss--)Newton step}
\label{subsubsec:operator_formulation}

This subsection gives a rigorous operator formulation of \eqref{eq:zN_reformulation}--\eqref{eq:zetaN_reformulation}.
We introduce the finite-width PDE-residual Jacobian operator, which yields factorizations of both the Gram--Gauss--Newton matrix $G_\theta^N$ and the operator $\CL T_{B_\theta^N}$.
For $\theta=(\theta^i)_{i=1}^N$,
let us define
\begin{equation}
    \label{eq:def_CJ}
    \CJ^N_\theta:\bbR^{N(d+2)}\rightarrow L^2(U),
    \qquad
    \left[\CJ^N_\theta u\right]\!(x)
    \coloneqq\left[\CL\left(\chi\,J^N_\theta\right)\right]\!(x)u .
\end{equation}

\begin{lemma}[Finite-width PDE-residual Jacobian operator~$\CJ^N_\theta$]
    \label{lem:operator_CJ}
    Let $\theta\in\bbR^{N(d+2)}$ be such that $D^N_\theta$ is symmetric positive definite.
    Then $\CJ^N_\theta$ defined in \eqref{eq:def_CJ} is well-defined and bounded, with
    \begin{equation}
        \label{eq:bound_CJ}
        \Nbig{\CJ^N_\theta}_{\bbR^{N(d+2)}\rightarrow L^2(U)}
        \leq
        \sqrt{\vol{U}}
        \operatorname*{ess\,sup}_{x\in U}
        \Nbig{
            \left[\CL\left(\chi\,J^N_\theta\right)\right]\!(x)
        }_2
        <\infty .
    \end{equation}
    Its adjoint
    $\left(\CJ^N_\theta\right)^*:L^2(U)\rightarrow\bbR^{N(d+2)}$
    is given by
    \begin{equation}
        \left(\CJ^N_\theta\right)^*v
        =
        \int_U
        \Big(\!
            \left[\CL\left(\chi\,J^N_\theta\right)\right]\!(x)
        \Big)^\top
        v(x)\,dx .
    \end{equation}
    Moreover,
    \begin{equation}
        \label{eq:G_and_LT_as_CJ}
        G^N_\theta
        =
        \left(\CJ^N_\theta\right)^*\CJ^N_\theta
        \quad\text{ and }\quad
        \frac{1}{\gamma}\CL T_{B^N_\theta}
        =
        \CJ^N_\theta
        \left(D^N_\theta\right)^{-1}
        \left(\CJ^N_\theta\right)^* .
    \end{equation}
    In particular, $\CL T_{B^N_\theta}$ is a bounded, self-adjoint
    and positive semi-definite operator on $L^2(U)$.
\end{lemma}

\begin{proof}
    Lemma~\ref{lem:bound_s}, applied with
    $\widetilde C=\max_{i=1,\dots,N}\abs{c^i}$, and the identity
    \begin{equation}
        \left[\CL\left(\chi\,J^N_{\theta^i}\right)\right]\!(x)
        =
        \frac{1}{N^\beta}
        \left[\CL\left(\chi\,s_{\theta^i}\right)\right]\!(x)
    \end{equation}
    show that
    $x\mapsto\left[\CL\left(\chi\,J^N_\theta\right)\right]\!(x)$
    is essentially bounded on $U$.
    Hence, for every $u\in\bbR^{N(d+2)}$,
    by Cauchy--Schwarz inequality
    \begin{equation}
    \begin{split}
        \Nbig{\CJ^N_\theta u}_{L^2(U)}
        &=
        \left(
            \int_U
            \abs{
                \left[\CJ^N_\theta u\right]\!(x)
            }^2
            dx
        \right)^{1/2}
        =
        \left(
            \int_U
            \abs{
                \left[\CL\left(\chi\,J^N_\theta\right)\right]\!(x)u
            }^2
            dx
        \right)^{1/2}
        \\
        &\leq
        \sqrt{\vol{U}}
        \operatorname*{ess\,sup}_{x\in U}
        \Nbig{
            \left[\CL\left(\chi\,J^N_\theta\right)\right]\!(x)
        }_2
        \N{u}_2,
    \end{split}
    \end{equation}
    which proves \eqref{eq:bound_CJ}.

    For $u\in\bbR^{N(d+2)}$ and $v\in L^2(U)$, it holds 
    \begin{equation}
    \begin{split}
        \left\langle\CJ^N_\theta u,v\right\rangle_{L^2(U)}
        &=
        \int_U
        v(x)
        \left[\CL\left(\chi\,J^N_\theta\right)\right]\!(x)u
        \,dx
        =
        \left\langle
            u,
            \int_U
            \left(
                \left[\CL\left(\chi\,J^N_\theta\right)\right]\!(x)
            \right)^\top
            v(x)\,dx
        \right\rangle_{\bbR^{N(d+2)}},
    \end{split}
    \end{equation}
    which proves the asserted formula for the adjoint $\left(\CJ^N_\theta\right)^*$ of the finite-width PDE-residual Jacobian operator $\CJ^N_\theta$.
    The definition of $G^N_\theta$ gives $G^N_\theta = \left(\CJ^N_\theta\right)^*\CJ^N_\theta$.
    Furthermore, \eqref{eq:T_NTK}--\eqref{eq:NTK} imply
    \begin{equation}
        \frac{1}{\gamma}\left[T_{B^N_\theta}v\right]\!(x)
        =
        \chi(x)J^N_\theta(x)
        \left(D^N_\theta\right)^{-1}
        \left(\CJ^N_\theta\right)^*v .
    \end{equation}
    Since $\CL$ is linear and acts component-wise, and since
    $\left(D^N_\theta\right)^{-1}\left(\CJ^N_\theta\right)^*v$
    is independent of $x$, applying $\CL$ to the preceding identity yields
    \begin{equation}
        \frac{1}{\gamma}\CL T_{B^N_\theta}
        =
        \CJ^N_\theta
        \left(D^N_\theta\right)^{-1}
        \left(\CJ^N_\theta\right)^*,
    \end{equation}
    which proves the second identity in
    \eqref{eq:G_and_LT_as_CJ}.
    Since $\CJ^N_\theta$ and its adjoint are bounded and
    $\left(D^N_\theta\right)^{-1}$ is a bounded linear map on
    $\bbR^{N(d+2)}$, this representation shows that
    $\CL T_{B^N_\theta}$ is bounded, with
    \begin{equation}
        \N{\frac{1}{\gamma}\CL T_{B^N_\theta}}_{L^2(U)\rightarrow L^2(U)}
        \leq
        \N{\CJ^N_\theta}_{\bbR^{N(d+2)}\rightarrow L^2(U)}^2
        \N{\left(D^N_\theta\right)^{-1}}_{\mathrm{op}} .
    \end{equation}
    Since $\left(D^N_\theta\right)^{-1}$ is symmetric, the same
    representation gives
    \begin{equation}
    \begin{split}
        \left(\frac{1}{\gamma}\CL T_{B^N_\theta}\right)^*
        &=
        \left(
            \CJ^N_\theta
            \left(D^N_\theta\right)^{-1}
            \left(\CJ^N_\theta\right)^*
        \right)^*
        =
        \CJ^N_\theta
        \left(D^N_\theta\right)^{-1}
        \left(\CJ^N_\theta\right)^*
        =
        \frac{1}{\gamma}\CL T_{B^N_\theta},
    \end{split}
    \end{equation}
    which proves self-adjointness. Moreover, since
    $\left(D^N_\theta\right)^{-1}$ is positive definite,
    \begin{equation}
    \begin{split}
        \left\langle
            \frac{1}{\gamma}\CL T_{B^N_\theta}v,v
        \right\rangle_{L^2(U)}
        &=
        \left\langle
            \CJ^N_\theta
            \left(D^N_\theta\right)^{-1}
            \left(\CJ^N_\theta\right)^*v,
            v
        \right\rangle_{L^2(U)} 
        =
        \left\langle
            \left(D^N_\theta\right)^{-1}
            \left(\CJ^N_\theta\right)^*v,
            \left(\CJ^N_\theta\right)^*v
        \right\rangle_{\bbR^{N(d+2)}} \\
        &\geq0,
    \end{split}
    \end{equation}
    which proves positive semi-definiteness and concludes the proof.
\end{proof}

\begin{remark}[Self-adjointness of the residual operators]
    \label{rem:selfadjointness}
    When viewed as an operator on $L^2(U)$, the PDE--NTK operator $T_{B^N_\theta}$ defined in \eqref{eq:T_NTK} need not be self-adjoint, whereas
    the operator $\CL T_{B^N_\theta}$ is self-adjoint by Lemma~\ref{lem:operator_CJ}.
    Its kernel is given by
    \begin{equation}
        \gamma\,\left[\CL\big(\chi J^N_\theta\big)\right]\!(x)\big(D^N_\theta\big)^{-1}
        \Big(\!\left[\CL\big(\chi J^N_\theta\big)\right]\!(x')\Big)^\top .
    \end{equation}
    Similarly, the PDE--NNTK integral operator $T_{\CB^N_\theta}$ need not be self-adjoint on $L^2(U)$, whereas its associated residual evolution operator $\CL T_{\CB^N_\theta}$ is self-adjoint.
    The same distinction applies to the limit PDE--NTK and PDE--NNTK integral operators $T_{B^*_0}$ and $T_{\CB^*_0}$, respectively: neither is required to be self-adjoint on $L^2(U)$, whereas the associated residual evolution operators $\CL T_{B^*_0}$ and $\CL T_{\CB^*_0}$ are self-adjoint.
\end{remark}

\begin{lemma}[Invertibility of the operators $\Id_{L^2(U)}+\frac{1}{\gamma}\CL T_{B^N_\theta}$ and $\Id_{L^2(U)} + \frac{1}{\gamma}\CL T_{B^*_0}$]
    \label{lem:estimate_zeta}
    Let $\theta\in\bbR^{N(d+2)}$ be such that $D^N_\theta$ is
    symmetric positive definite.
    Then $\Id_{L^2(U)}+\frac{1}{\gamma}\CL T_{B^N_\theta}:L^2(U)\rightarrow L^2(U)$
    is a bounded bijection, and its inverse satisfies
    \begin{equation}
        \label{eq:lem:estimate_zeta}
        \N{
            \left(
                \Id_{L^2(U)}+\frac{1}{\gamma}\CL T_{B^N_\theta}
            \right)^{-1}
        }_{L^2(U)\rightarrow L^2(U)}
        \leq1 .
    \end{equation}
    Hence, for every $r\in L^2(U)$, the equation
    $\big( \Id_{L^2(U)}+\frac{1}{\gamma}\CL T_{B^N_\theta} \big)\zeta^N = r$
    admits a unique solution $\zeta^N\in L^2(U)$ satisfying $\N{\zeta^N}_{L^2(U)} \leq \N{r}_{L^2(U)}$.
    
    The same statement holds for
    \begin{equation}
        \Id_{L^2(U)} + \frac{1}{\gamma}\CL T_{B^*_0}.
    \end{equation}
\end{lemma}

\begin{proof}
    We first study the operator $\Id_{L^2(U)}+\frac{1}{\gamma}\CL T_{B^N_\theta}$.
    Define the bilinear form
    \begin{equation}
        b(v_1,v_2)
        \coloneqq
        \left\langle
            \left(
                \Id_{L^2(U)}+\frac{1}{\gamma}\CL T_{B^N_\theta}
            \right)v_1,
            v_2
        \right\rangle_{L^2(U)},
        \qquad
        v_1,v_2\in L^2(U).
    \end{equation}
    Since $\CL T_{B^N_\theta}$ is bounded, self-adjoint and positive
    semi-definite by Lemma~\ref{lem:operator_CJ},
    $b$ is bounded on $L^2(U)$
    since
    \begin{equation}
    \begin{split}
        \absbig{b(v_1,v_2)}
        &\leq
        \left(
            1+
            \Nbig{
                \frac{1}{\gamma}\CL T_{B^N_\theta}
            }_{L^2(U)\rightarrow L^2(U)}
        \right)
        \N{v_1}_{L^2(U)}
        \N{v_2}_{L^2(U)},
    \end{split}
    \end{equation}
    for $v_1,v_2\in L^2(U)$.
    Moreover, $b$ is coercive on $L^2(U)$ since
    \begin{equation}
    \begin{split}
        b(v,v)
        &=
        \N{v}_{L^2(U)}^2
        +
        \left\langle
            \frac{1}{\gamma}\CL T_{B^N_\theta}v,v
        \right\rangle_{L^2(U)}
        \geq
        \N{v}_{L^2(U)}^2,
    \end{split}
    \end{equation}
    for every $v\in L^2(U)$.
    Hence, by the Lax--Milgram theorem, for every $r\in L^2(U)$ there exists
    a unique $\zeta^N\in L^2(U)$ such that
    \begin{equation}
        \label{eq:proof:lem:estimate_zeta:20}
        b(\zeta^N,v)
        =
        \left\langle r,v\right\rangle_{L^2(U)}
        \qquad
        \text{ for every }v\in L^2(U),
    \end{equation}
    or equivalently, $\big( \Id_{L^2(U)}+\frac{1}{\gamma}\CL T_{B^N_\theta} \big)\zeta^N = r$.
    Thus,
    $\Id_{L^2(U)}+\frac{1}{\gamma}\CL T_{B^N_\theta}$ is a bijection on
    $L^2(U)$.

    Since $b$ is coercive on $L^2(U)$, it holds by using \eqref{eq:proof:lem:estimate_zeta:20} and Cauchy--Schwarz inequality that 
    \begin{equation}
    \begin{split}
        \N{\zeta^N}_{L^2(U)}^2
        &\leq
        b(\zeta^N,\zeta^N)
        =
        \left\langle r,\zeta^N\right\rangle_{L^2(U)}
        \leq
        \N{r}_{L^2(U)}
        \N{\zeta^N}_{L^2(U)}
    \end{split}
    \end{equation}
    and consequently, $\N{\zeta^N}_{L^2(U)} \leq \N{r}_{L^2(U)}$.
    Since $\zeta^N = \big( \Id_{L^2(U)}+\frac{1}{\gamma}\CL T_{B^N_\theta} \big)^{-1}r$
    and the estimate holds for every $r\in L^2(U)$, taking the
    supremum over all nonzero $r\in L^2(U)$ proves
    \eqref{eq:lem:estimate_zeta}.

    Let us now investigate the operator $\Id_{L^2(U)} + \frac{1}{\gamma}\CL T_{B^*_0}$.
    Since Lemma~\ref{lem:spectralpropertiesTB} shows that
    $\CL T_{B^*_0}$ is bounded, self-adjoint and positive
    semi-definite and since $\gamma>0$,
    the same properties hold for
    $\frac{1}{\gamma}\CL T_{B^*_0}$.
    Repeating the preceding argument proves the corresponding
    conclusions for
    $\Id_{L^2(U)}+\frac{1}{\gamma}\CL T_{B^*_0}$.
\end{proof}

\begin{lemma}[Push-through identity]
    \label{lem:pushthrough}
    Let $\theta\in\bbR^{N(d+2)}$ be such that $D^N_\theta$ is symmetric positive definite.
    Then $D^N_\theta+G^N_\theta$ is symmetric positive definite,
    and, for every $r\in L^2(U)$, it holds
    \begin{equation}
        \label{eq:pushthrough}
        \left(
            D^N_\theta+G^N_\theta
        \right)^{-1}
        \left(\CJ^N_\theta\right)^*r
        =
        \left(D^N_\theta\right)^{-1}
        \left(\CJ^N_\theta\right)^*
        \left(
            \Id_{L^2(U)}+\frac{1}{\gamma}\CL T_{B^N_\theta}
        \right)^{-1}r .
    \end{equation}
    The (Gauss--)Newton step in
    \eqref{eq:regularizedNewtonStep} therefore satisfies
    \begin{equation}
        \label{eq:zN_operator}
        \begin{aligned}
            z^N
            &=
            \left(D^N_\theta\right)^{-1}
            \left(\CJ^N_\theta\right)^*\zeta^N
            \quad\text{ with }\quad
            \zeta^N
            \coloneqq
            \left(
                \Id_{L^2(U)}+\frac{1}{\gamma}\CL T_{B^N_\theta}
            \right)^{-1}
            \left(h-\CL g^N_\theta\right).
        \end{aligned}
    \end{equation}
    In particular, for every $i=1,\dots,N$,
    \begin{equation}
        \label{eq:zN_neuronwise}
        \left(z^N\right)^i
        =
        \left(D^N_{\theta^i}\right)^{-1}
        \int_U
        \Big(\!
            \left[
                \CL\left(\chi\,J^N_{\theta^i}\right)
            \right]\!(x)
        \Big)^\top
        \zeta^N(x)\,dx .
    \end{equation}
\end{lemma}

\begin{proof}
    By Lemma~\ref{lem:operator_CJ}, $G^N_\theta = \left(\CJ^N_\theta\right)^*\CJ^N_\theta$.
    Hence, since $D^N_\theta$ is symmetric,
    $D^N_\theta+G^N_\theta$ is symmetric
    Moreover, since $D^N_\theta$ is positive definite, it holds
    \begin{equation}
    \begin{split}
        \left\langle
            \left(D^N_\theta+G^N_\theta\right)u,u
        \right\rangle_{\bbR^{N(d+2)}}
        &=
        \left\langle
            D^N_\theta u,u
        \right\rangle_{\bbR^{N(d+2)}}
        +
        \Nbig{\CJ^N_\theta u}_{L^2(U)}^2
        >0
    \end{split}
    \end{equation}
    for every nonzero
    $u\in\bbR^{N(d+2)}$.
    Therefore,
    $D^N_\theta+G^N_\theta$ is symmetric positive definite and therefore invertible.
    Again by Lemma~\ref{lem:operator_CJ}, $\Id_{L^2(U)} + \CJ^N_\theta\left(D^N_\theta\right)^{-1}\left(\CJ^N_\theta\right)^* = \Id_{L^2(U)} + \frac{1}{\gamma}\CL T_{B^N_\theta}$,
    which is invertible on $L^2(U)$ as of
    Lemma~\ref{lem:estimate_zeta}.
    Expanding both sides gives the operator identity
    \begin{equation}
    \begin{split}
        &
        \left(
            D^N_\theta
            +
            \left(\CJ^N_\theta\right)^*\CJ^N_\theta
        \right)
        \left(D^N_\theta\right)^{-1}
        \left(\CJ^N_\theta\right)^* =
        \left(\CJ^N_\theta\right)^*
        \left(
            \Id_{L^2(U)}
            +
            \CJ^N_\theta
            \left(D^N_\theta\right)^{-1}
            \left(\CJ^N_\theta\right)^*
        \right).
    \end{split}
    \end{equation}
    Multiplying this identity from the left by
    $\big( D^N_\theta + \left(\CJ^N_\theta\right)^*\CJ^N_\theta\big)^{-1}$
    and from the right by $\big(\Id_{L^2(U)} + \CJ^N_\theta \left(D^N_\theta\right)^{-1} \left(\CJ^N_\theta\right)^* \big)^{-1}$
    yields
    \begin{equation}
    \begin{split}
        \left(
            D^N_\theta
            +
            \left(\CJ^N_\theta\right)^*\CJ^N_\theta
        \right)^{-1}
        \left(\CJ^N_\theta\right)^*
        =
        \left(D^N_\theta\right)^{-1}
        \left(\CJ^N_\theta\right)^*
        \left(
            \Id_{L^2(U)}
            +
            \CJ^N_\theta
            \left(D^N_\theta\right)^{-1}
            \left(\CJ^N_\theta\right)^*
        \right)^{-1}.
    \end{split}
    \end{equation}
    Substituting the notations $G^N_\theta = \left(\CJ^N_\theta\right)^*\CJ^N_\theta$ and $\frac{1}{\gamma}\CL T_{B^N_\theta} = \CJ^N_\theta \left(D^N_\theta\right)^{-1}\left(\CJ^N_\theta\right)^*$
    proves \eqref{eq:pushthrough}.
    
    By \eqref{eq:regularizedNewton_DG} and the formula for
    $\left(\CJ^N_\theta\right)^*$ in
    Lemma~\ref{lem:operator_CJ}, the regularized (Gauss--)Newton
    system can be written as
    \begin{equation}
        \label{eq:regularizedNewton_CJ}
        \left(
            D^N_\theta+G^N_\theta
        \right)z^N
        =
        \left(\CJ^N_\theta\right)^*
        \left(h-\CL g^N_\theta\right).
    \end{equation}
    Applying \eqref{eq:pushthrough} with
    $r=h-\CL g^N_\theta$ gives
    \begin{equation}
        z^N
        =
        \left(D^N_\theta\right)^{-1}
        \left(\CJ^N_\theta\right)^*
        \left(
            \Id_{L^2(U)}+\frac{1}{\gamma}\CL T_{B^N_\theta}
        \right)^{-1}
        \left(h-\CL g^N_\theta\right),
    \end{equation}
    which proves \eqref{eq:zN_operator}.
    By the formula for the adjoint in Lemma~\ref{lem:operator_CJ},
    it holds
    \begin{equation}
        \left(
            \left(\CJ^N_\theta\right)^*\zeta^N
        \right)^i
        =
        \int_U
        \Big(\!
            \left[
                \CL\left(\chi\,J^N_{\theta^i}\right)
            \right]\!(x)
        \Big)^\top
        \zeta^N(x)\,dx
    \end{equation}
    for every $i=1,\dots,N$.
    Since $\left(D^N_\theta\right)^{-1}$ is block-diagonal with blocks $\left(D^N_{\theta^i}\right)^{-1}$,
    \eqref{eq:zN_operator} therefore gives
    \begin{equation}
        \left(z^N\right)^i
        =
        \left(D^N_{\theta^i}\right)^{-1}
        \int_U
        \left(
            \left[
                \CL\left(\chi\,J^N_{\theta^i}\right)
            \right]\!(x)
        \right)^\top
        \zeta^N(x)\,dx ,
    \end{equation}
    which proves \eqref{eq:zN_neuronwise} and concludes the proof.
\end{proof}

\subsection{Properties of the kernel operators~\texorpdfstring{$T_{B^N_\theta}$, $T_{\CB^N_\theta}$, $T_{B^*_0}$ and $T_{\CB^*_0}$}{}}
\label{subsec:kernel_mapping}

This section establishes the $L^2(U)\rightarrow\CH^2_{(0)}(U)$ mapping properties and operator-norm bounds for the finite-width and limit kernel operators needed when proving convergence of the finite-width dynamics to its overparameterized limit as $N\rightarrow\infty$.

We first establish an estimate that transfers spatial regularity of a kernel to its associated integral operator.

\begin{lemma}[Kernel-to-operator estimate]
    \label{lem:kernel_to_operator}
    Let $B\in L^2(U\times U)$ and assume that
    $D^x_\nu B\in L^2(U\times U)$ for every multi-index
    $\abs{\nu}\leq2$, where the derivatives are understood in the weak sense with respect to the first variable.
    Then $T_B:L^2(U)\rightarrow\CH^2(U)$ is bounded
    and its derivatives can be represented as
    \begin{equation}
        \label{eq:diff_under_integral_kernel}
        D_\nu\big[T_Bv\big]
        =
        T_{D^x_\nu B}v
        =
        \int_U D^x_\nu B(\dummy,x')v(x')\,dx'
        \qquad
        \text{ for every }\abs{\nu}\leq2,
    \end{equation}
    and
    \begin{equation}
        \label{eq:kernel_to_operator}
        \Nbig{T_B}_{L^2(U)\rightarrow\CH^2(U)}
        \leq
        \sum_{0\leq\abs{\nu}\leq2}
        \Nbig{D^x_\nu B}_{L^2(U\times U)}.
    \end{equation}
\end{lemma}

\begin{proof}
Let $v\in L^2(U)$. For every multi-index $\abs{\nu}\leq2$, the integral operator with kernel $D^x_\nu B$ is Hilbert--Schmidt and therefore
\begin{equation}
\Nbig{T_{D^x_\nu B}v}_{L^2(U)}
\leq
\Nbig{D^x_\nu B}_{L^2(U\times U)}
\N{v}_{L^2(U)}.
\end{equation}
Let $(v_n)_{n\in\bbN}\subset\CC^\infty_c(U)$ satisfy
$v_n\rightarrow v$ in $L^2(U)$.
For every $\varphi\in\CC^\infty_c(U)$, the definition of the distributional derivative gives
\begin{equation}
\begin{split}
\int_U
\big[T_Bv_n\big](x)D_\nu\varphi(x)\,dx
&=
\int_{U\times U}
B(x,x')v_n(x')D_\nu\varphi(x)\,dx'\,dx
\\
&=
(-1)^{\abs{\nu}}
\int_{U\times U}
D^x_\nu B(x,x')v_n(x')\varphi(x)\,dx'\,dx
\\
&=
(-1)^{\abs{\nu}}
\int_U
\big[T_{D^x_\nu B}v_n\big](x)\varphi(x)\,dx.
\end{split}
\end{equation}
Since $T_B$ and $T_{D^x_\nu B}$ are Hilbert--Schmidt operators, their continuity from $L^2(U)$ to $L^2(U)$ allows us to pass to the limit $n\rightarrow\infty$.
Consequently,
\begin{equation}
D_\nu\big[T_Bv\big]
=
T_{D^x_\nu B}v
\end{equation}
in the sense of distributions, proving \eqref{eq:diff_under_integral_kernel}.
Thus $T_Bv\in\CH^2(U)$.
Summing the preceding Hilbert--Schmidt estimate over $\abs{\nu}\leq2$ and using the norm convention of Remark~\ref{rem:CL_linearity} gives
\eqref{eq:kernel_to_operator}.
\end{proof}

\begin{lemma}
    \label{lem:rank_one_kernel}
    Let $\vartheta,\vartheta'\in\bbR^{d+2}$ satisfy
    $\absnormal{\vartheta_c},\absnormal{\vartheta'_c}\leq\widetilde C$, let $M\in\bbR^{(d+2)\times(d+2)}$, and define the kernel
    \begin{equation}
        K_{\vartheta,\vartheta',M}(x,x')
        \coloneqq
        \chi(x)s_{\vartheta}(x)M
        \Big(\!
            \left[
                \CL(\chi\,s_{\vartheta'})
            \right]\!(x')
        \Big)^\top.
    \end{equation}
    Then
    \begin{equation}
        \label{eq:proof:rank_one}
        \sum_{\abs{\nu}\leq2}
        \Nbig{
            D^x_\nu
            K_{\vartheta,\vartheta',M}
        }_{L^2(U\times U)}
        \leq
        C\vol{U}\N{M}_{\mathrm{op}}
        \left(
            1+\N{\vartheta_w}_2^2
        \right)
        \left(
            1+\N{\vartheta'_w}_2^2
        \right)
    \end{equation}
    for a constant $C=C(\sigma,\CD,\widetilde C)>0$.
    
    Moreover, if $\vartheta_1,\vartheta'_1,\vartheta_2,\vartheta'_2\in\bbR^{d+2}$
    satisfy
    $\absnormal{(\vartheta_1)_c},
    \absnormal{(\vartheta'_1)_c},
    \absnormal{(\vartheta_2)_c},
    \absnormal{(\vartheta'_2)_c}
    \leq\widetilde C$, then
    \begin{equation}
        \label{eq:proof:rank_one2}
    \begin{split}
        &\sum_{\abs{\nu}\leq2}
        \N{
            D^x_\nu
            \left[
                K_{\vartheta_1,\vartheta'_1,M}
                -
                K_{\vartheta_2,\vartheta'_2,M}
            \right]
        }_{L^2(U\times U)}
        \\
        &\qquad\,\leq
        C\vol{U}\N{M}_{\mathrm{op}}
        \left[
            \left(
                1+\N{(\vartheta_1)_w}_2^2
                +
                \N{(\vartheta_2)_w}_2^2
            \right)
            \left(
                1+\N{(\vartheta'_1)_w}_2^2
            \right)
            \N{\vartheta_1-\vartheta_2}_2\right.
            \\
            &\qquad\,\quad\,
            \left.+
            \left(
                1+\N{(\vartheta_2)_w}_2^2
            \right)
            \left(
                1+\N{(\vartheta'_1)_w}_2^2
                +
                \N{(\vartheta'_2)_w}_2^2
            \right)
            \N{\vartheta'_1-\vartheta'_2}_2
        \right].
    \end{split}
    \end{equation}
\end{lemma}

\begin{proof}
    Since $D^x_\nu$ differentiates with respect to the first spatial variable $x$, while
    $M\left(\left[\CL(\chi\,s_{\vartheta'})\right]\!(x')\right)^\top$
    is independent of $x$,
    it holds
    \begin{equation}
        \left[
            D^x_\nu
            K_{\vartheta,\vartheta',M}
        \right]\!(x,x')
        =
        \left[
            D^x_\nu
            \big(
                \chi s_{\vartheta}
            \big)
        \right]\!(x)
        M
        \Big(\!
            \left[
                \CL(\chi\,s_{\vartheta'})
            \right]\!(x')
        \Big)^\top.
    \end{equation}
    Therefore, Lemma~\ref{lem:bound_s} gives
    \begin{equation}
    \begin{split}
        \Nbig{
            D^x_\nu
            K_{\vartheta,\vartheta',M}
        }_{L^2(U\times U)}
        &\leq
        \vol{U}
        \N{M}_{\mathrm{op}}
        \sup_{x\in U}
        \Nbig{
            \left[
                D^x_\nu
                \big(
                    \chi s_{\vartheta}
                \big)
            \right]\!(x)
        }_2
        \operatorname*{ess\,sup}_{x'\in U}
        \Nbig{
            \left[
                \CL(\chi\,s_{\vartheta'})
            \right]\!(x')
        }_2
        \\
        &\leq
        C\vol{U}\N{M}_{\mathrm{op}}
        \big(
            1+\N{\vartheta_w}_2^2
        \big)
        \big(
            1+\N{\vartheta'_w}_2^2
        \big).
    \end{split}
    \end{equation}
    Summing over $\abs{\nu}\leq2$ and increasing $C$ to absorb the finite number of multi-indices proves \eqref{eq:proof:rank_one}.
    
    For the second part of the statement, insert the mixed term
    $K_{\vartheta_2,\vartheta'_1,M}$ to obtain
    \begin{equation}
    \begin{split}
        &
        K_{\vartheta_1,\vartheta'_1,M}
        -
        K_{\vartheta_2,\vartheta'_2,M}
        =
        K_{\vartheta_1,\vartheta'_1,M}
        -
        K_{\vartheta_2,\vartheta'_1,M}
        +
        K_{\vartheta_2,\vartheta'_1,M}
        -
        K_{\vartheta_2,\vartheta'_2,M}.
    \end{split}
    \end{equation}
    Consequently, by product rule
    \begin{equation}
    \begin{split}
        \left[
            D^x_\nu
            \left(
                K_{\vartheta_1,\vartheta'_1,M}
                -
                K_{\vartheta_2,\vartheta'_2,M}
            \right)
        \right]\!(x,x')
        &=
        \left[
            D^x_\nu
            \left(
                \chi s_{\vartheta_1}
                -
                \chi s_{\vartheta_2}
            \right)
        \right]\!(x)
        M
        \left(
            \left[
                \CL(\chi\,s_{\vartheta'_1})
            \right]\!(x')
        \right)^\top
        \\
        &\quad\,
        +
        \left[
            D^x_\nu
            \big(
                \chi s_{\vartheta_2}
            \big)
        \right]\!(x)
        M
        \left(
            \left[
                \CL
                \left(
                    \chi\,s_{\vartheta'_1}
                    -
                    \chi\,s_{\vartheta'_2}
                \right)
            \right]\!(x')
        \right)^\top.
    \end{split}
    \end{equation}
    Applying Lemma~\ref{lem:bound_s_Lipschitz} to the difference factors and Lemma~\ref{lem:bound_s} to the remaining factors gives
    \begin{equation}
    \begin{split}
        &
        \Nbig{
            D^x_\nu
            \left[
                K_{\vartheta_1,\vartheta'_1,M}
                -
                K_{\vartheta_2,\vartheta'_2,M}
            \right]
        }_{L^2(U\times U)}
        \\
        &\quad\,\leq
        C\vol{U}\N{M}_{\mathrm{op}}
        \left[
            \left(
                1+\N{(\vartheta_1)_w}_2^2
                +
                \N{(\vartheta_2)_w}_2^2
            \right)
            \left(
                1+\N{(\vartheta'_1)_w}_2^2
            \right)
            \N{\vartheta_1-\vartheta_2}_2
            \right.\\
            &\qquad\qquad\qquad\qquad\qquad\quad
            \left.+
            \left(
                1+\N{(\vartheta_2)_w}_2^2
            \right)
            \left(
                1+\N{(\vartheta'_1)_w}_2^2
                +
                \N{(\vartheta'_2)_w}_2^2
            \right)
            \N{\vartheta'_1-\vartheta'_2}_2
        \right].
    \end{split}
    \end{equation}
    Summing over $\abs{\nu}\leq2$ and increasing $C$ once more proves \eqref{eq:proof:rank_one2}.
\end{proof}

\begin{lemma}
    \label{lem:eta_H10}
    Let $\Psi\in\CC^2(\overline U)$.
    Then $\chi\Psi\in\CH^2_{(0)}(U)$.
\end{lemma}

\begin{proof}
    Assumptions~\ref{asm:D} and~\ref{asm:NN_chi} imply that $\chi\Psi\in\CC^2(\overline U)\subset\CH^2(U)$ and $(\chi\Psi)|_{\partial U}=0$.
    The inclusion in $\CH^2(U)$ follows because $\overline U$ is compact, and therefore all derivatives of $\chi\Psi$ up to order two are bounded and hence square-integrable on $U$.
    Since $\partial U$ is Lipschitz, the trace theorem gives
    $\chi\Psi\in\CH^1_0(U)$.
    Together with $\chi\Psi\in\CH^2(U)$, this shows the assertion.
\end{proof}

\begin{lemma}[Limit kernel operators]
    \label{lem:standardNTKB:operatornorm}
    The operators $T_{B^*_0}$ and $T_{\CB^*_0}$ are bounded from
    $L^2(U)$ into $\CH^2_{(0)}(U)$ and
    \begin{equation}
        \label{eq:bound_TBstar}
        \Nbig{
            T_{\CB^*_0}
        }_{L^2(U)\rightarrow\CH^2(U)}
        \leq
        \frac{1}{\gamma}
        \Nbig{
            T_{B^*_0}
        }_{L^2(U)\rightarrow\CH^2(U)}
        \leq
        \frac{C_B}{\gamma}.
    \end{equation}
    for a constant
    $C_B=C_B(\sigma,\CD,\mu_0)>0$.
\end{lemma}

\begin{proof}
    For every multi-index $\abs{\nu}\leq2$, differentiation inside the expectation in \eqref{eq:NTK*} is justified by the dominated convergence theorem such that it holds
    \begin{equation}
    \label{eq:bound_DxB0}
    \begin{split}
        \absbig{
            \left[
                D^x_\nu B^*_0
            \right]\!(x,x')
        }
        &=
        \abs{\bbE
        \left[
            \left[
                D^x_\nu
                \big(
                    \chi s_\vartheta
                \big)
            \right]\!(x)
            \Big(\!
                \left[
                    \CL(\chi\,s_\vartheta)
                \right]\!(x')
            \Big)^\top
        \right]}
        \\
        &\leq
        \bbE
        \left[
            \Nbig{
                \left[
                    D^x_\nu
                    \big(
                        \chi s_\vartheta
                    \big)
                \right]\!(x)
            }_2
            \Nbig{
                \left[
                    \CL(\chi\,s_\vartheta)
                \right]\!(x')
            }_2
        \right]
        \\
        &\leq
        C\bbE[\omega_1^2]
        <\infty
    \end{split}
    \end{equation}
    for almost every $(x,x')\in U\times U$, with a bound independent of $(x,x')$,
    where the first inequality follows from Jensen's inequality for the expectation and the Cauchy--Schwarz inequality, while the second follows from Lemma~\ref{lem:bound_s}, applied with $\widetilde C=C_{c,0}$, which is permitted by Assumption~\ref{asm:NN_mu0}\ref{asm:NN_mu0ii}.
    Moreover, $C\omega_1^2$ is integrable since
    $\omega_1^2\leq2\big(1+\N{w^1_0}_2^4\big)$ and $q\geq4$ by Assumption~\ref{asm:NN_mu0}\ref{asm:NN_mu0iv}, which provides the domination required above.

    Since $U$ is bounded, \eqref{eq:bound_DxB0} implies that
    $D^x_\nu B^*_0\in L^2(U\times U)$ for every multi-index $\abs{\nu}\leq2$.
    Lemma~\ref{lem:kernel_to_operator}, specifically
    \eqref{eq:kernel_to_operator}, therefore gives
    \begin{equation}
        \Nbig{
            T_{B^*_0}
        }_{L^2(U)\rightarrow\CH^2(U)}
        \leq
        C_B
    \end{equation}
    for a constant
    $C_B=C_B(\sigma,\CD,\mu_0)>0$.

    Let $v\in L^2(U)$ and define
    \begin{equation}
        \Psi_v(x)
        \coloneqq
        \bbE_{\vartheta\sim\mu_0}
        \left[
            s_\vartheta(x)
            \int_U
            \Big(\!
                \left[
                    \CL(\chi\,s_\vartheta)
                \right]\!(x')
            \Big)^\top
            v(x')\,dx'
        \right].
    \end{equation}
    By \eqref{eq:NTK*} and Fubini's theorem, $T_{B^*_0}v = \chi\,\Psi_v$.
    Lemma~\ref{lem:bound_s}, with $\chi$ replaced by the constant function $1$, and the Cauchy--Schwarz inequality give, for every multi-index 
    $\abs{\nu}\leq2$,
    \begin{equation}
    \begin{split}
        &
        \Nbig{
            \left[
                D^x_\nu s_\vartheta
            \right]\!(x)
        }_2
        \N{
            \int_U
            \Big(\!
                \left[
                    \CL(\chi\,s_\vartheta)
                \right]\!(x')
            \Big)^\top
            v(x')\,dx'
        }_2
        \leq
        C\,\sqrt{\vol{U}}\,
        \omega_1^2
        \N{v}_{L^2(U)}.
    \end{split}
    \end{equation}
    Since $\bbE[\omega_1^2]<\infty$, the dominated convergence theorem permits differentiation inside the expectation and shows that
    $\Psi_v\in\CC^2(\overline U)$.
    Lemma~\ref{lem:eta_H10}, applied with
    $\Psi=\Psi_v$, therefore gives $T_{B^*_0}v = \chi\,\Psi_v \in \CH^2_{(0)}(U)$.
    Thus $T_{B^*_0}$ is bounded from
    $L^2(U)$ into $\CH^2_{(0)}(U)$.

    Finally, by \eqref{eq:T_NNTK*},
    \begin{equation}
        T_{\CB^*_0}
        =
        \frac{1}{\gamma}
        T_{B^*_0}
        \left(
            \Id_{L^2(U)}
            +
            \frac{1}{\gamma}\CL T_{B^*_0}
        \right)^{-1}.
    \end{equation}
    Lemma~\ref{lem:estimate_zeta} shows that the inverse has
    $L^2(U)\rightarrow L^2(U)$ operator norm at most one.
    Consequently,
    \begin{equation}
        \Nbig{
            T_{\CB^*_0}
        }_{L^2(U)\rightarrow\CH^2(U)}
        \leq
        \frac{1}{\gamma}
        \Nbig{
            T_{B^*_0}
        }_{L^2(U)\rightarrow\CH^2(U)}
        \leq
        \frac{C_B}{\gamma},
    \end{equation}
    which proves \eqref{eq:bound_TBstar}.
    The same representation shows that the range of
    $T_{\CB^*_0}$ is contained in the range of
    $T_{B^*_0}$ and therefore in $\CH^2_{(0)}(U)$.
\end{proof}

\begin{lemma}[Finite-width kernel operators]
\label{lem:finitewidthNTKB:operatornorm}
    Let $\theta=(\theta^i)_{i=1}^N$ satisfy
    $\absnormal{c^i}\leq\widetilde C$ for every $i=1,\dots,N$, and
    suppose that $D^N_\theta$ is symmetric positive definite.
    Assume that there exist weights $\omega_i>0$ and a constant
    $C_V<\infty$ such that
    $1+\N{w^i}_2^2\leq3\omega_i$ and
    $\gamma\Nbig{\big(d^N_{\theta^i}\big)^{-1}}_{\mathrm{op}}\leq2$
    for every $i=1,\dots,N$, as well as
    $\frac{1}{N}\sum_{i=1}^N\omega_i^2\leq C_V$.
    Then $T_{B^N_\theta}$ and $T_{\CB^N_\theta}$ are bounded from
    $L^2(U)$ into $\CH^2_{(0)}(U)$ and
    \begin{equation}
    \label{eq:bound_TBN_uniform}
    \Nbig{
        T_{\CB^N_\theta}
    }_{L^2(U)\rightarrow\CH^2(U)}
    \leq\frac{1}{\gamma}
    \Nbig{
        T_{B^N_\theta}
    }_{L^2(U)\rightarrow\CH^2(U)}
    \leq
    \frac{C}{\gamma},
    \end{equation}
    for a constant
    $C=C(\sigma,\CD,\widetilde C,C_V)>0$ independent of $N$.
\end{lemma}

\begin{proof}
    For every multi-index $\abs{\nu}\leq2$, differentiation of the finite sum in
    \eqref{eq:kernel_finite} gives
    \begin{equation}
    \label{eq:bound_DxBN}
    \begin{split}
    \absbig{
    \left[
    D^x_\nu B^N_\theta
    \right]\!(x,x')
    }
    &=
    \abs{
    \frac{1}{N}
    \sum_{i=1}^N
    \left[
    D^x_\nu
    \big(
    \chi s_{\theta^i}
    \big)
    \right]\!(x)
    \gamma\big(d^N_{\theta^i}\big)^{-1}
    \Big(\!
    \left[
    \CL(\chi\,s_{\theta^i})
    \right]\!(x')
    \Big)^\top
    }
    \\
    &\leq
    \frac{1}{N}
    \sum_{i=1}^N
    \Nbig{
    \left[
    D^x_\nu
    \big(
    \chi s_{\theta^i}
    \big)
    \right]\!(x)
    }_2
    \gamma\Nbig{
    \big(d^N_{\theta^i}\big)^{-1}
    }_{\mathrm{op}}
    \Nbig{
    \left[
    \CL(\chi\,s_{\theta^i})
    \right]\!(x')
    }_2
    \\
    &\leq
    \frac{C}{N}
    \sum_{i=1}^N
    \omega_i^2
    \leq
    CC_V
    \end{split}
    \end{equation}
    for almost every $(x,x')\in U\times U$, with a bound independent of $(x,x')$, where the first inequality follows from
    the triangle inequality and the Cauchy--Schwarz inequality, while
    the second follows from Lemma~\ref{lem:bound_s} and the assumptions
    $1+\N{w^i}_2^2\leq3\omega_i$ and
    $\gamma\Nbig{\big(d^N_{\theta^i}\big)^{-1}}_{\mathrm{op}}\leq2$.
    
    Since $U$ is bounded, \eqref{eq:bound_DxBN} implies that
    $D^x_\nu B^N_\theta\in L^2(U\times U)$ for every
    $\abs{\nu}\leq2$.
    Lemma~\ref{lem:kernel_to_operator}, specifically
    \eqref{eq:kernel_to_operator}, therefore gives
    \begin{equation}
    \Nbig{
    T_{B^N_\theta}
    }_{L^2(U)\rightarrow\CH^2(U)}
    \leq
    C,
    \end{equation}
    after increasing $C$.

   Let $v\in L^2(U)$ and define
    \begin{equation}
        \Psi^N_{\theta,v}(x)
        \coloneqq
        \frac{1}{N}
        \sum_{i=1}^N
        s_{\theta^i}(x)
        \gamma\big(d^N_{\theta^i}\big)^{-1}
        \int_U
        \Big(\!
            \left[
                \CL(\chi\,s_{\theta^i})
            \right]\!(x')
        \Big)^\top
        v(x')\,dx'.
    \end{equation}
    By \eqref{eq:kernel_finite}, $T_{B^N_\theta}v = \chi\,\Psi^N_{\theta,v}$.
    Lemma~\ref{lem:bound_s} and the Cauchy--Schwarz inequality show that the integrals in the definition of $\Psi^N_{\theta,v}$ are finite.
    Moreover, the definition of $s_{\theta^i}$ in \eqref{eq:s} and Assumption~\ref{asm:NN_sigma'} imply that $s_{\theta^i}\in\CC^2(\overline U)$ for every $i=1,\dots,N$.
    Since $\Psi^N_{\theta,v}$ is a finite linear combination of these functions, $\Psi^N_{\theta,v} \in \CC^2(\overline U)$.
    Lemma~\ref{lem:eta_H10}, applied with
    $\Psi=\Psi^N_{\theta,v}$, therefore gives $T_{B^N_\theta}v = \chi\,\Psi^N_{\theta,v} \in \CH^2_{(0)}(U)$.
    Thus $T_{B^N_\theta}$ is bounded from
    $L^2(U)$ into $\CH^2_{(0)}(U)$.
    
    Finally, by \eqref{eq:NNTK},
    \begin{equation}
        T_{\CB^N_\theta}
        =
        \frac{1}{\gamma}T_{B^N_\theta}
        \left(
            \Id_{L^2(U)}
            +
            \frac{1}{\gamma}\CL T_{B^N_\theta}
        \right)^{-1}.
    \end{equation}
    Since $D^N_\theta$ is symmetric positive definite,
    Lemma~\ref{lem:estimate_zeta} shows that the inverse has
    $L^2(U)\rightarrow L^2(U)$ operator norm at most one.
    Consequently,
    \begin{equation}
        \Nbig{
            T_{\CB^N_\theta}
        }_{L^2(U)\rightarrow\CH^2(U)}
        \leq\frac{1}{\gamma}
        \Nbig{
            T_{B^N_\theta}
        }_{L^2(U)\rightarrow\CH^2(U)}
        \leq
        \frac{C}{\gamma},
    \end{equation}
    which proves \eqref{eq:bound_TBN_uniform}.
    The same representation shows that the range of
    $T_{\CB^N_\theta}$ is contained in the range of
    $T_{B^N_\theta}$ and therefore in $\CH^2_{(0)}(U)$.
\end{proof}

\subsection{High-probability bounds of finite-width NN along training trajectory}
\label{subsec:trajectory_control}

This section establishes uniform bounds that hold with high probability on the finite-width training trajectory $(g_{\theta_k}^N(\x))_{k=0,\dots,K}$ for a fixed number of training iterations $k=0,\dots,K$.
The first lemma controls the NN function at initialization and is then used to construct the high-probability event~$\Omega_R$ on which the trajectory bounds and the invertibility of $D^N_{\theta_k}$ hold simultaneously throughout training.

\begin{lemma}[Bounds on NN function~$g^N_{\theta_0}$ at initialization]
    \label{lem:estimate_NNfunction_init}
    At initialization, the NN function $g^N_{\theta_0}$ satisfies
    \begin{equation}
        \bbE\N{g^N_{\theta_0}}^2_{\CH^2(U)}
        \leq \frac{C}{N^{2\beta-1}}
        \quad\text{ and }\quad
        \bbE\N{\CL g^N_{\theta_0}}^2_{L^2(U)}
        \leq\frac{C}{N^{2\beta-1}}
    \end{equation}
    for a constant $C=C(\sigma,\CD,\mu_0)>0$.
\end{lemma}

\begin{proof}
    For every multi-index $\abs{\nu}\leq2$, define the $L^2(U)$-valued random variables
    \begin{equation}
        X_i^\nu
        \coloneqq D_\nu\phi_{\theta^i_0}
        = D_\nu\!\left(\chi\,c^i_0\sigma(w^i_0\cdot\dummy+\eta^i_0)\right),
        \qquad i=1,\dots,N.
    \end{equation}
    Lemma~\ref{lem:bound_s}\ref{lem:bound_s:phi}, applied with $\widetilde C=C_{c,0}$, gives the bound
    \begin{equation}
        \N{X_i^\nu}_{L^2(U)}^2
        \leq C\vol{U}\omega_i^2.
    \end{equation}
    Since $q\geq4$ by Assumption~\ref{asm:NN_mu0}\ref{asm:NN_mu0iv}, it holds $\bbE[\omega_i^2]<\infty$, and hence $X_i^\nu$ is square-integrable as an $L^2(U)$-valued random variable.
    Moreover, the $X_i^\nu$'s are i.i.d.\@ and centered, since each $X_i^\nu$ is linear in $c^i_0$, while $c^i_0$ is independent of $(w^i_0,\eta^i_0)$ and has mean zero by Assumption~\ref{asm:NN_mu0}\ref{asm:NN_mu0i}--\ref{asm:NN_mu0ii}.
    Independence and centering therefore imply that all cross terms vanish, and
    \begin{equation}
    \begin{split}
        \bbE\N{D_\nu g^N_{\theta_0}}_{L^2(U)}^2
        &=\frac{1}{N^{2\beta}}\sum_{i=1}^N\bbE\N{X_i^\nu}_{L^2(U)}^2
        =\frac{1}{N^{2\beta-1}}\bbE\N{X_1^\nu}_{L^2(U)}^2\\
        &\leq\frac{C\vol{U}}{N^{2\beta-1}}\bbE\big[\omega_1^2\big],
    \end{split}
    \end{equation}
    where the last inequality follows from the preceding bound.

    Let $m_d$ denote the number of multi-indices $\abs{\nu}\leq2$.
    By the convention for the $\CH^2(U)$-norm in Remark~\ref{rem:CL_linearity} and the Cauchy--Schwarz inequality, it holds
    \begin{equation}
        \bbE\N{g^N_{\theta_0}}_{\CH^2(U)}^2
        \leq m_d\sum_{0\leq\abs{\nu}\leq2}\bbE\N{D_\nu g^N_{\theta_0}}_{L^2(U)}^2
        \leq\frac{C}{N^{2\beta-1}},
    \end{equation}
    which proves the first part of the statement.

    For the second part, define the $L^2(U)$-valued random variables
    \begin{equation}
        Y_i\coloneqq \CL\phi_{\theta^i_0}
        = \CL\left(\chi\,c^i_0\sigma(w^i_0\cdot\dummy+\eta^i_0)\right),
        \qquad i=1,\dots,N.
    \end{equation}
    The same argument shows that the $Y_i$ are i.i.d., centered, and square-integrable in $L^2(U)$, while Lemma~\ref{lem:bound_s}\ref{lem:bound_s:phi} gives
    $\bbE\N{Y_1}_{L^2(U)}^2\leq C\vol{U}\bbE[\omega_1^2]$.
    Therefore,
    \begin{equation}
        \bbE\N{\CL g^N_{\theta_0}}_{L^2(U)}^2
        =\frac{1}{N^{2\beta-1}}\bbE\N{Y_1}_{L^2(U)}^2
        \leq\frac{C}{N^{2\beta-1}},
    \end{equation}
    concluding the proof.
\end{proof}

The next lemma provides the path-wise high-probability bounds of the finite-width NN along the training trajectory used throughout the finite-width convergence argument.

\begin{lemma}[Uniform bounds along training trajectory]
    \label{prop:NTK_limit_bounds_z}
    Let $\delta\in(0,1/2]$ and $K\in\bbN$ be a given number of training steps. Fix $C_{R,0}>C_h$.
    Denote by $C_c,C_w,C_V,C_R,C_z,C_\gamma>0$ the constants defined in \eqref{eq:def_constants} below, which depend only on $\alpha$, $\gamma$, $K$, $\delta$, $q$, $\sigma$, $\CD$, $\mu_0$, and $C_{R,0}$.
    Assume that the number of neurons~$N$ is large enough such that \eqref{eq:thm:asm:Nlargeenough} is satisfied.
\begin{comment}
    \begin{equation}
        \label{eq:lem:bounds_z}
        \bbP\left(
        \begin{aligned}
            \max_{i=1,\dots,N}\omega_i
            &\leq C_wN^{2/q},
            &\frac1N\sum_{i=1}^N\omega_i^2
            &\leq C_V, \\
            \max_{\substack{i=1,\dots,N\\k=0,\dots,K}}\absbig{c^i_k}
            &\leq C_c,
            &\max_{\substack{i=1,\dots,N\\k=0,\dots,K}}\N{\theta^i_k-\theta^i_0}_2
            &\leq1, \\
            \max_{k=0,\dots,K}\N{h-\CL g^N_{\theta_k}}_{L^2(U)}
            &\leq C_R,
            &\max_{\substack{i=1,\dots,N\\k=0,\dots,K}}
            \frac{\N{\big(z^N_k\big)^i}_2}{\omega_i}
            &\leq\frac{C_z}{\gamma N^{1-\beta}}, \\
            \min_{k=0,\dots,K}\lambda_{\min}\big(D^N_{\theta_k}\big)
            &>0,
            &\max_{\substack{i=1,\dots,N\\k=0,\dots,K}}
            \gamma\N{(d^N_{\theta^i_k})^{-1}}_{\mathrm{op}}
            &\leq2
        \end{aligned}
        \right)\geq1-\delta.
    \end{equation}
\end{comment}
    Then there exists an event $\Omega_R$ with $\bbP(\Omega_R)\geq1-\delta$ on which
    \begin{subequations}
        \label{eq:lem:bounds_z}
        \begin{equation}
            \max_{i=1,\dots,N}\omega_i\leq C_wN^{2/q}
            \quad\text{ and }\quad
            \frac1N\sum_{i=1}^N\omega_i^2\leq C_V.
        \end{equation}
        as well as
        \begin{equation}
            \absbig{c^i_k}\leq C_c
            \quad\text{ and }\quad
            \N{\theta^i_k-\theta^i_0}_2\leq1
        \end{equation}
        for all $i=1,\dots,N$ and all $k=0,\dots,K$.
        On the same event, it holds
        \begin{equation}
            \N{h-\CL g^N_{\theta_k}}_{L^2(U)}\leq C_R
        \end{equation}
        and 
        \begin{equation}
            \N{\big(z^N_k\big)^i}_2\leq\frac{C_z\,\omega_i}{\gamma N^{1-\beta}}
        \end{equation}
        for all $i=1,\dots,N$ and all $k=0,\dots,K$.
        Furthermore, for all $k=0,\dots,K$, the matrix $D^N_{\theta_k}$ is symmetric positive definite with
        \begin{equation}
            \N{\big(D^N_{\theta_k}\big)^{-1}}_{\mathrm{op}}\leq\frac2\gamma N^{2\beta-1}
            \quad\text{ and }\quad
            \gamma\N{(d^N_{\theta^i_k})^{-1}}_{\mathrm{op}}\leq2
        \end{equation}    
        for all $i=1,\dots,N$.
    \end{subequations}
\end{lemma}

\begin{proof}
    Set $C_c\coloneqq C_{c,0}+1$, choose $C=C(\sigma,\CD,C_c)>0$ sufficiently large for Lemmas~\ref{lem:bound_s} and~\ref{lem:bound_s_Lipschitz}, and let $C_S\coloneqq C\sqrt{\vol{U}}$ be the corresponding constant from Lemma~\ref{lem:estimate_S_operatornorm_i}.
    The remaining constants are fixed as
    \begin{equation}
        \label{eq:def_constants}
    \begin{split}
         c_1&\coloneqq6\sqrt{\vol{U}}\,C,
        \quad  c_2\coloneqq3\sqrt{\vol{U}}\,C c_1, \\
        C_w&\coloneqq1+\left(\frac{3\CM_q}{\delta}\right)^{2/q},
        \quad C_V\coloneqq\frac{3\,\bbE\big[\omega_1^2\big]}{\delta},
        \quad
        C_R\coloneqq C_{R,0}\left(1+\frac{\alpha c_2C_V}{\gamma}\right)^{K},
        \quad C_z\coloneqq c_1C_R, \\
        C_\gamma&\coloneqq6\,C_SC_wC_R .
    \end{split}
    \end{equation}
    These constants are finite by Assumption~\ref{asm:NN_mu0}\ref{asm:NN_mu0iv} and independent of $N$ and $k$.
    Although $C_R$, $C_z$, and $C_\gamma$ depend on the fixed regularizer parameter $\gamma$, they are determined before $N$ is chosen.
    Since $\rho>0$ and $2\beta-1>0$, both conditions in \eqref{eq:thm:asm:Nlargeenough} hold for all sufficiently large $N$.

    \emph{Step 1: Construction of event $\Omega_R$.}
    Define the individual events
    \begin{align}
        \Omega_{R,1}
        &\coloneqq\left\{\N{h-\CL g^N_{\theta_0}}_{L^2(U)}\leq C_{R,0}\right\}, \\
        \Omega_{R,2}
        &\coloneqq\left\{\max_{i=1,\dots,N}\N{w^i_0}_2
        \leq\left(\frac{3N\CM_q}{\delta}\right)^{1/q}\right\}, \\
        \Omega_{R,3}
        &\coloneqq\left\{\frac1N\sum_{i=1}^N\omega_i^2
        \leq\frac{3\,\bbE[\omega_1^2]}{\delta}\right\}.
    \end{align}
    First, for the event $\Omega_{R,1}$, since $\N{h-\CL g^N_{\theta_0}}_{L^2(U)}\leq C_h + \N{\CL g^N_{\theta_0}}_{L^2(U)}$ by the triangle inequality and Assumption~\ref{asm:h}, and since $\bbE\N{\CL g^N_{\theta_0}}^2_{L^2(U)} \leq \frac{C}{N^{2\beta-1}}$ by Lemma~\ref{lem:estimate_NNfunction_init},
    we have by Markov's inequality for any $C_{R,0}>C_h$ that
    \begin{equation}
        \label{eq:proof:Omega_residual}
    \begin{split}
        \bbP\big(\Omega_{R,1}\big)
        &\geq\bbP\left(\N{\CL g^N_{\theta_0}}_{L^2(U)}\leq C_{R,0}-C_h\right) \\
        &\geq1-\frac{\bbE\N{\CL g^N_{\theta_0}}^2_{L^2(U)}}{(C_{R,0}-C_h)^2}
        \geq1-\frac{C}{(C_{R,0}-C_h)^2}\frac{1}{N^{2\beta-1}}
        \geq1-\frac\delta3,
    \end{split}
    \end{equation}
    where the last step is by the first condition in \eqref{eq:thm:asm:Nlargeenough}.
    Second, for the event $\Omega_{R,2}$, since by Assumption~\ref{asm:NN_mu0} the $w^i_0$'s are i.i.d., we obtain with a union bound together with Markov's inequality that
    \begin{equation}
        \label{eq:proof:Omega_weights}
    \begin{split}
        \bbP\big(\Omega_{R,2}\big)
        &= 1 - \bbP\left(\max_{i=1,\dots,N}\N{w^i_0}_2 > \left(\frac{3N\CM_q}{\delta}\right)^{1/q} \right)\\
        &\geq 1 - \sum_{i=1}^N \bbP\left(\N{w^i_0}_2 > \left(\frac{3N\CM_q}{\delta}\right)^{1/q} \right) \\
        &\geq 1 - N\frac{\delta}{3N\CM_q}\bbE \N{w^1_0}_2^q
        = 1 - N\frac{\delta}{3N\CM_q}\CM_q
        = 1 - \frac{\delta}{3},
    \end{split}
    \end{equation}
    where Assumption~\ref{asm:NN_mu0}\ref{asm:NN_mu0iv} was used in the penultimate equality.
    On $\Omega_{R,2}$, by the definition of $\omega_i$, it holds
    \begin{equation}
        \label{eq:proof:omega_max}
        \max_{i=1,\dots,N} \omega_i
        = 1 + \max_{i=1,\dots,N} \N{w^i_0}_2^2
        \leq 1 + \left(\frac{3N\CM_q}{\delta}\right)^{2/q}
        \leq C_w N^{2/q},
    \end{equation}
    where the last step follows from the definition of $C_w$ in \eqref{eq:def_constants} and $N\geq1$.
    Third, for the event $\Omega_{R,3}$, first note that by the definition of $\omega_i$, $(a+b)^2\leq2(a^2+b^2)$, and Lyapunov's inequality,
    \begin{equation}
        \label{eq:proof:omega_second_moment}
        \bbE\big[\omega_i^2\big]
        =\bbE\left[\left(1+\N{w^i_0}_2^2\right)^2\right]
        \leq2\left(1+\bbE\N{w^i_0}_2^4\right)
        \leq2\left(1+\CM_q^{4/q}\right)
        <\infty,
    \end{equation}
    the next-to-last step by Jensen's inequality and $q\geq4$, which holds by Assumption~\ref{asm:NN_mu0}\ref{asm:NN_mu0iv}.
    Since the $\omega_i$ are identically distributed by Assumption~\ref{asm:NN_mu0}, the non-negative random variable $\frac1N\sum_{i=1}^N\omega_i^2$ has expectation $\bbE\big[\frac1N\sum_{i=1}^N\omega_i^2\big]=\bbE[\omega_1^2]$, so that Markov's inequality gives
    \begin{equation}
        \label{eq:proof:Omega_average}
    \begin{split}
        \bbP\big(\Omega_{R,3}\big)
        &=1-\bbP\left(\frac1N\sum_{i=1}^N\omega_i^2>\frac{3\,\bbE[\omega_1^2]}{\delta}\right) \\
        &\geq1-\frac{\delta}{3\,\bbE[\omega_1^2]}\,\bbE\left[\frac1N\sum_{i=1}^N\omega_i^2\right]
        =1-\frac{\delta}{3\,\bbE[\omega_1^2]}\,\bbE\big[\omega_1^2\big]
        =1-\frac{\delta}{3}.
    \end{split}
    \end{equation}
    By the definition of $C_V$ in \eqref{eq:def_constants}, it holds on $\Omega_{R,3}$ that
    \begin{equation}
        \label{eq:proof:omega_avg}
        \frac1N\sum_{i=1}^N\omega_i^2
        \leq C_V.
    \end{equation}
    Finally, setting $\Omega_R\coloneqq\Omega_{R,1}\cap\Omega_{R,2}\cap\Omega_{R,3}$, a union bound over the three complements gives
    \begin{equation}
        \label{eq:proof:OmegaR}
        \bbP\big(\Omega_R\big)
        \geq1-\bbP\big(\Omega_{R,1}^c\big)-\bbP\big(\Omega_{R,2}^c\big)-\bbP\big(\Omega_{R,3}^c\big)
        \geq1-3\cdot\frac{\delta}{3}
        =1-\delta,
    \end{equation}
    and both \eqref{eq:proof:omega_max} and \eqref{eq:proof:omega_avg} hold on $\Omega_R$, as required in \eqref{eq:lem:bounds_z}.

    \emph{Step 2: Induction start $k=0$.}
    The proof proceeds using an inductive argument. We start with $k=0$.
    
    \emph{Step 2a: Bound on parameters at $k=0$.}
    By Assumption~\ref{asm:NN_mu0}\ref{asm:NN_mu0ii} it holds $\absnormal{c^i_0}\leq C_{c,0}\leq C_c$ and trivially $\N{\theta^i_0-\theta^i_0}_2=0\leq1$, which proves the corresponding bounds in \eqref{eq:lem:bounds_z} for $k=0$.

    \emph{Step 2b: Bound on residual at $k=0$.}
    On $\Omega_R$, we further have $\Nbig{h-\CL g^N_{\theta_0}}_{L^2(U)}\leq C_{R,0}\leq C_R$.

    \emph{Step 2c: Positive definiteness of $D^N_{\theta_0}$.}
    By Step~2a the hypothesis $\absnormal{c^i_0}\leq C_c$ of Lemma~\ref{lem:estimate_Dinv_operatornorm} is met, so we may use the constant $C_S$ introduced before \eqref{eq:def_constants}.
    Using \eqref{eq:proof:omega_max}, $\Nbig{h-\CL g^N_{\theta_0}}_{L^2(U)}\leq C_R$, and the second condition in \eqref{eq:thm:asm:Nlargeenough} together with $C_\gamma=6C_SC_wC_R$, we obtain
    \begin{equation}
        \label{eq:proof:thm:NTK_limit:condition_gamma0}
        \frac{2C_S\max_{i=1,\dots,N}\left(1+\N{w^i_0}_2^2\right)}{N^{1-\beta}}\N{h-\CL g^N_{\theta_0}}_{L^2(U)}
        \leq\frac{2C_SC_wC_RN^{2/q}}{N^{1-\beta}}
        =\frac{C_\gamma}{3N^{\rho}}
        \leq\frac\gamma3
        \leq\gamma .
    \end{equation}
    Thus condition~\eqref{eq:cond_D} is satisfied in case of the Newton method; no such condition is needed in case of the Gauss--Newton method.
    Hence Lemma~\ref{lem:estimate_Dinv_operatornorm} applies and yields that $D^N_{\theta_0}$ is symmetric positive definite with
    \begin{equation}
        \label{eq:proof:thm:NTK_limit:bound_D0inv}
        \N{\big(D^N_{\theta_0}\big)^{-1}}_{\mathrm{op}}\leq\frac2\gamma N^{2\beta-1}
        \quad\text{ and }\quad
        \gamma\N{\big(d^N_{\theta^i_0}\big)^{-1}}_{\mathrm{op}}\leq2
        \quad\text{ for all }i=1,\dots,N.
    \end{equation}
    By the coercivity established in the proof of
    Lemma~\ref{lem:estimate_zeta}, it holds that
    $\N{\zeta_0^N}_{L^2(U)}^2 \leq \big\langle \big( \Id_{L^2(U)} +\frac{1}{\gamma}\CL T_{B_{\theta_0}^N} \big)\zeta_0^N, \zeta_0^N \big\rangle_{L^2(U)}$.
    Hence, by Cauchy--Schwarz inequality,
    \begin{equation}
        \label{eq:proof:thm:NTK_limit:bound_zeta0}
        \N{\zeta_0^N}_{L^2(U)}
        \leq
        \N{
            \left(
                \Id_{L^2(U)}
                +\frac{1}{\gamma}\CL T_{B_{\theta_0}^N}
            \right)\zeta_0^N
        }_{L^2(U)}
        =
        \Nbig{h-\CL g_{\theta_0}^N}_{L^2(U)}
        \leq C_R,
    \end{equation}
    where the equality follows from the definition
    of $\zeta_0^N$ in \eqref{eq:zetaN_reformulation}.

    \emph{Step 2d: Bound on (Gauss--)Newton step at $k=0$.}
    Using the neuron-wise representation in \eqref{eq:zN_reformulation_neuronwise}, the definition of the operator norm, the Cauchy--Schwarz inequality, Lemma~\ref{lem:bound_s} together with $\CL(\chi J^N_{\theta^i})=\frac{1}{N^\beta}\CL(\chi s_{\theta^i})$, and finally \eqref{eq:proof:thm:NTK_limit:bound_D0inv}--\eqref{eq:proof:thm:NTK_limit:bound_zeta0}, it follows for every $i=1,\dots,N$ that
    \begin{equation}
        \label{eq:proof:thm:NTK_limit:bound_z0}
    \begin{split}
        \N{\big(z^N_0\big)^i}_2
        &=\N{\big(D^N_{\theta^i_0}\big)^{-1}\int_U\Big(\CL\big(\chi\,J^N_{\theta^i_0}\big)(x)\Big)^\top\zeta^N_0(x)\,dx}_2 \\
        &\leq\N{\big(D^N_{\theta^i_0}\big)^{-1}}_{\mathrm{op}}\,\sqrt{\vol{U}}\,\operatorname*{ess\,sup}_{x\in U}\N{\left[\CL\big(\chi\,J^N_{\theta^i_0}\big)\right]\!(x)}_2\,\N{\zeta^N_0}_{L^2(U)} \\
        &\leq\frac2\gamma N^{2\beta-1}\cdot\sqrt{\vol{U}}\cdot\frac{C\omega_i}{N^\beta}\cdot C_R
        \;\leq\;\frac{ c_1C_R\omega_i}{\gamma N^{1-\beta}}
        =\frac{C_z\,\omega_i}{\gamma N^{1-\beta}},
    \end{split}
    \end{equation}
    where $ c_1=6\sqrt{\vol{U}}\,C\geq2\sqrt{\vol{U}}\,C$ and $C_z= c_1C_R$ as in \eqref{eq:def_constants}.
    This proves the update bound in \eqref{eq:lem:bounds_z} for $k=0$.

    \emph{Step 3: Induction step $k\mapsto k+1$.}
    Let $k\in\{0,\dots,K-1\}$ and assume that the trajectory bounds in \eqref{eq:lem:bounds_z} hold for every index $\ell=0,\dots,k$.
    In particular, $D^N_{\theta_\ell}$ is symmetric positive definite and satisfies
    \begin{equation}
        \label{eq:proof:D_inverseboundes_inductionhypothesis}
        \Nbig{\big(D^N_{\theta_\ell}\big)^{-1}}_{\mathrm{op}}
        \leq\frac2\gamma N^{2\beta-1}
        \quad\text{ and }\quad
        \gamma\Nbig{\big(d^N_{\theta^i_\ell}\big)^{-1}}_{\mathrm{op}}\leq2
    \end{equation}
    for every $i=1,\dots,N$ and every $\ell=0,\dots,k$.
    By Step~2 this is the case for $k=0$.
    We further note that $\Nbig{\theta^i_k-\theta^i_0}_2\leq1$ implies
    \begin{equation}
        \label{eq:proof:omega_transfer}
        1+\N{w^i_k}_2^2\leq1+\left(\N{w^i_0}_2+1\right)^2\leq3\left(1+\N{w^i_0}_2^2\right)=3\,\omega_i ,
    \end{equation}
    where the first step is due to $\N{w^i_k}_2\leq\N{w^i_0}_2 + \N{w^i_k-w^i_0}_2 \leq \N{w^i_0}_2 + \Nbig{\theta^i_k-\theta^i_0}_2 \leq \N{w^i_0}_2+1$ while the second uses $(a+1)^2\leq2a^2+2$.
    The same bound holds for any $\vartheta$ on the segment $[\theta^i_0,\theta^i_k]$, by convexity of the ball.

    \emph{Step 3a: Bound on parameters.}
    According to the parameter update \eqref{eq:regularizedNewton}, a telescopic sum, triangle inequality and the induction hypothesis for the update bound in \eqref{eq:lem:bounds_z} for the steps $0,\dots,k$,
    \begin{equation}
        \label{eq:proof:movement_neuronwise}
        \N{\theta^i_{k+1}-\theta^i_0}_2
        \leq\sum_{\ell=0}^{k} \N{\theta^i_{\ell+1}-\theta^i_\ell}_2
        =\alpha\sum_{\ell=0}^{k}\N{\big(z^N_\ell\big)^i}_2
        \leq\frac{\alpha(k+1)C_z\,\omega_i}{\gamma N^{1-\beta}}
        \leq\frac{\alpha KC_z\omega_i}{\gamma N^{1-\beta}},
    \end{equation}
    where $k+1\leq K$ was used in the last step.
    Combining \eqref{eq:proof:movement_neuronwise} with \eqref{eq:proof:omega_max} and $\frac{1}{N^{1-\beta}}N^{2/q}=\frac{1}{N^{\rho}}$ gives
    \begin{equation}
        \label{eq:proof:movement}
        \Nbig{\theta^i_{k+1}-\theta^i_0}_2
        \leq\frac{\alpha KC_zC_w}{\gamma N^{\rho}}
        \leq1,
    \end{equation}
    where the last step holds due to the second condition in \eqref{eq:thm:asm:Nlargeenough}.
    In particular, $\absnormal{c^i_{k+1}}\leq\absnormal{c^i_0}+\absnormal{c^i_{k+1}-c^i_0}\leq\absnormal{c^i_0}+\Nbig{\theta^i_{k+1}-\theta^i_0}_2\leq C_{c,0}+1=C_c$, so that the corresponding bound in \eqref{eq:lem:bounds_z} holds at $k+1$. Furthermore, since $\Nbig{\theta^i_{k+1}-\theta^i_0}_2\leq1$ as of \eqref{eq:proof:movement}, $1+\N{w^i_{k+1}}_2^2\leq3\omega_i$ by \eqref{eq:proof:omega_transfer}.

    \emph{Step 3b: Bound on residual.}
    Using the neuron-wise representation in \eqref{eq:zN_reformulation_neuronwise}
    and proceeding along the lines of \eqref{eq:proof:thm:NTK_limit:bound_z0}, employing in particular Lemma~\ref{lem:bound_s},
    it follows
    for every $i=1,\dots,N$ that
    \begin{equation}
        \label{eq:proof:z_sharp}
    \begin{split}
        \N{\big(z^N_k\big)^i}_2
        &=\N{\big(D^N_{\theta^i_k}\big)^{-1}\int_U\Big(\CL\big(\chi\,J^N_{\theta^i_k}\big)(x)\Big)^\top\zeta^N_k(x)\,dx}_2 \\
        &\leq\N{\big(D^N_{\theta^i_k}\big)^{-1}}_{\mathrm{op}}\,\sqrt{\vol{U}}\,\operatorname*{ess\,sup}_{x\in U}\N{\CL\big(\chi\,J^N_{\theta^i_k}\big)(x)}_2\,\N{\zeta^N_k}_{L^2(U)} \\
        &\leq\frac2\gamma N^{2\beta-1}\cdot\sqrt{\vol{U}}\cdot\frac{3C\omega_i}{N^\beta}\,\N{\zeta^N_k}_{L^2(U)}
        \leq\frac{ c_1\,\omega_i}{\gamma}\Nbig{h-\CL g^N_{\theta_k}}_{L^2(U)}\frac{1}{N^{1-\beta}},
    \end{split}
    \end{equation}
    where we used that $D^N_{\theta_k}$ is symmetric positive definite and satisfies the required inverse bounds~\eqref{eq:proof:D_inverseboundes_inductionhypothesis} by the induction hypothesis, that $\absnormal{c^i_k}\leq C_c$ and $1+\N{w^i_k}_2^2\leq3\omega_i$ by \eqref{eq:lem:bounds_z} and \eqref{eq:proof:omega_transfer} (the factor $3$ of \eqref{eq:proof:omega_transfer} is absorbed into $ c_1=6\sqrt{\vol{U}}\,C$), and that $\Nbig{\zeta^N_{k}}_{L^2(U)}\leq\Nbig{h-\CL g^N_{\theta_k}}_{L^2(U)}$ thanks to an analogous argument as the one that led to \eqref{eq:proof:thm:NTK_limit:bound_zeta0}.
    Note that the update bound in \eqref{eq:lem:bounds_z} at step $k$ is the special case of \eqref{eq:proof:z_sharp} obtained by inserting $\Nbig{h-\CL g^N_{\theta_k}}_{L^2(U)}\leq C_R$ and recalling $C_z= c_1C_R$.

    Next, by \eqref{eq:NN} and the linearity of $\CL$,
    $\CL g^N_\theta=\frac{1}{N^\beta}\sum_{i=1}^N\CL\phi_{\theta^i}$ in $L^\infty(U)$.
    Assumption~\ref{asm:CL} makes $\CL:\CC^2(\overline U)\to L^\infty(U)$ a bounded linear map.
    Consequently, \eqref{eq:commutation_CL_parameter_gradient} and Lemma~\ref{lem:sigma_consequences}\ref{lem:sigma:C4} show that $\vartheta\mapsto\CL\phi_\vartheta$ is continuously differentiable as an $L^\infty(U)$-valued map, with derivative represented by $\CL(\chi s_\vartheta)$.
    Integrating this derivative along $[\theta^i_k,\theta^i_{k+1}]$, which by Step~3a is contained in $\{\vartheta:\absnormal{c}\leq C_c\}\cap\{\vartheta:1+\N{w}_2^2\leq3\omega_i\}$, and using Lemma~\ref{lem:bound_s} and $\theta^i_{k+1}-\theta^i_k=\alpha(z^N_k)^i$, we obtain
    \begin{equation}
        \label{eq:proof:residual_increment}
    \begin{split}
        \Nbig{\CL g^N_{\theta_{k+1}}-\CL g^N_{\theta_k}}_{L^\infty(U)}
        &\leq\frac{1}{N^\beta}\sum_{i=1}^N
        \sup_{\vartheta\in[\theta^i_k,\theta^i_{k+1}]}
        \left(\operatorname*{ess\,sup}_{x\in U}\N{\left[\CL\big(\chi\,s_\vartheta\big)\right]\!(x)}_2\right)\Nbig{\theta^i_{k+1}-\theta^i_k}_2 \\
        &\leq\frac{3\alpha C}{N^\beta}\sum_{i=1}^N\omega_i\Nbig{\big(z^N_k\big)^i}_2
        \\
        &\leq\frac{3\alpha C c_1}{\gamma}
        \Nbig{h-\CL g^N_{\theta_k}}_{L^2(U)}
        \frac{1}{N^{1-\beta}}\frac{1}{N^\beta}\sum_{i=1}^N\omega_i^2 \\
        &=\frac{3\alpha C c_1}{\gamma}
        \Nbig{h-\CL g^N_{\theta_k}}_{L^2(U)}\frac1N\sum_{i=1}^N\omega_i^2 \\
        &\leq\frac{3\alpha C c_1C_V}{\gamma}
        \Nbig{h-\CL g^N_{\theta_k}}_{L^2(U)}.
    \end{split}
    \end{equation}
    To obtain the inequality in the third line in the former display, we use the bound \eqref{eq:proof:z_sharp}, while the last inequality leverages \eqref{eq:proof:omega_avg}.
    Note that retaining the neuron-wise factor $\omega_i$ in \eqref{eq:proof:z_sharp} was essential for the last step.

    By triangle inequality
    $\Nbig{h-\CL g^N_{\theta_{k+1}}}_{L^2(U)}\leq\Nbig{h-\CL g^N_{\theta_{k}}}_{L^2(U)}+\Nbig{\CL g^N_{\theta_{k+1}}-\CL g^N_{\theta_{k}}}_{L^2(U)}$.
    Using $\N{f}_{L^2(U)}\leq\sqrt{\vol{U}}\N{f}_{L^\infty(U)}$ and \eqref{eq:proof:residual_increment}, we can estimate
    \begin{equation}
    \begin{split}
        \N{\CL g^N_{\theta_{k+1}}-\CL g^N_{\theta_{k}}}_{L^2(U)}
        &\leq \sqrt{\vol{U}}\Nbig{\CL g^N_{\theta_{k+1}}-\CL g^N_{\theta_k}}_{L^\infty(U)} \\
        &\leq \frac{3\sqrt{\vol{U}}\alpha C c_1C_V}{\gamma}
        \Nbig{h-\CL g^N_{\theta_k}}_{L^2(U)}.
    \end{split}
    \end{equation}
    Recalling $ c_2=3\sqrt{\vol{U}}\,C c_1$ from \eqref{eq:def_constants},
    we obtain the multiplicative recursion
    \begin{equation}
        \label{eq:proof:residual_recursion}
        \Nbig{h-\CL g^N_{\theta_{k+1}}}_{L^2(U)}
        \leq\left(1+\frac{\alpha c_2C_V}{\gamma}\right)\Nbig{h-\CL g^N_{\theta_{k}}}_{L^2(U)} .
    \end{equation}

    Since the derivation of \eqref{eq:proof:z_sharp} and \eqref{eq:proof:residual_recursion} used only the hypotheses of Step~3 at index $k$, which by the cumulative form of that hypothesis are available at every index $\ell\in\{0,\dots,k\}$, \eqref{eq:proof:residual_recursion} holds with $k$ replaced by any such $\ell$.
    Iterating \eqref{eq:proof:residual_recursion} recursively for $\ell=0,\dots,k$ and using $\Nbig{h-\CL g^N_{\theta_0}}_{L^2(U)}\leq C_{R,0}$ from Step~2b gives
    \begin{equation}
        \label{eq:proof:residual_closed}
        \Nbig{h-\CL g^N_{\theta_{k+1}}}_{L^2(U)}
        \leq C_{R,0}\left(1+\frac{\alpha c_2C_V}{\gamma}\right)^{k+1}
        \leq C_{R,0}\left(1+\frac{\alpha c_2C_V}{\gamma}\right)^{K}
        =C_R ,
    \end{equation}
    where the last step follows from the definition of $C_R$ in \eqref{eq:def_constants}.
    This proves the residual bound in \eqref{eq:lem:bounds_z} at $k+1$.

    \emph{Step 3c: Positive definiteness of $D^N_{\theta_{k+1}}$.}
    By Step~3a the hypothesis $\absnormal{c^i_{k+1}}\leq C_c$ of Lemma~\ref{lem:estimate_Dinv_operatornorm} is met, so we may use the constant $C_S$ introduced before \eqref{eq:def_constants}.
    Furthermore, $1+\N{w^i_{k+1}}_2^2\leq3\omega_i$ thanks to Step~3a.
    Using \eqref{eq:proof:omega_max}, $\Nbig{h-\CL g^N_{\theta_{k+1}}}_{L^2(U)}\leq C_R$, and the second condition in \eqref{eq:thm:asm:Nlargeenough} together with $C_\gamma=6C_SC_wC_R$, we obtain
    \begin{equation}
        \label{eq:Gamma_bound}
        \frac{2C_S\max_{i=1,\dots,N}\left(1+\N{w^i_{k+1}}_2^2\right)}{N^{1-\beta}}\Nbig{h-\CL g^N_{\theta_{k+1}}}_{L^2(U)}
        \leq\frac{6C_SC_wC_RN^{2/q}}{N^{1-\beta}}
        =\frac{C_\gamma}{N^{\rho}}
        \leq\gamma .
    \end{equation}
    Thus condition~\eqref{eq:cond_D} is satisfied in case of the Newton method; no such condition is needed in case of the Gauss--Newton method.
    Hence Lemma~\ref{lem:estimate_Dinv_operatornorm} applies and yields that $D^N_{\theta_{k+1}}$ is symmetric positive definite with
    \begin{equation}
        \label{eq:proof:thm:NTK_limit:bound_Dk+1inv}
        \N{\big(D^N_{\theta_{k+1}}\big)^{-1}}_{\mathrm{op}}\leq\frac2\gamma N^{2\beta-1}
        \quad\text{ and }\quad
        \gamma\N{\big(d^N_{\theta^i_{k+1}}\big)^{-1}}_{\mathrm{op}}\leq2
        \quad\text{ for all }i=1,\dots,N.
    \end{equation}

    \emph{Step 3d: Bound on (Gauss--)Newton step at $k+1$.}
    Using the neuron-wise representation in \eqref{eq:zN_reformulation_neuronwise}
    and proceeding along the lines of \eqref{eq:proof:thm:NTK_limit:bound_z0}, employing in particular Lemma~\ref{lem:bound_s},
    it follows
    for every $i=1,\dots,N$ that
    \begin{equation}
    \begin{split}
        \N{\big(z^N_{k+1}\big)^i}_2
        &=\N{\big(D^N_{\theta^i_{k+1}}\big)^{-1}\int_U\Big(\CL\big(\chi\,J^N_{\theta^i_{k+1}}\big)(x)\Big)^\top\zeta^N_{k+1}(x)\,dx}_2 \\
        &\leq\N{\big(D^N_{\theta^i_{k+1}}\big)^{-1}}_{\mathrm{op}}\,\sqrt{\vol{U}}\,\operatorname*{ess\,sup}_{x\in U}\N{\CL\big(\chi\,J^N_{\theta^i_{k+1}}\big)(x)}_2\,\N{\zeta^N_{k+1}}_{L^2(U)} \\
        &\leq\frac2\gamma N^{2\beta-1}\cdot\sqrt{\vol{U}}\cdot\frac{3C\omega_i}{N^\beta}\,\N{\zeta^N_{k+1}}_{L^2(U)}\\
        &\leq\frac{ c_1\,\omega_i}{\gamma}\Nbig{h-\CL g^N_{\theta_{k+1}}}_{L^2(U)}\frac{1}{N^{1-\beta}}
        \leq \frac{ c_1C_R\,\omega_i}{\gamma}\frac{1}{N^{1-\beta}} = \frac{C_z\,\omega_i}{\gamma N^{1-\beta}},
    \end{split}
    \end{equation}
    where we used that $D^N_{\theta_{k+1}}$ is symmetric positive definite and satisfies the required inverse bounds~\eqref{eq:proof:thm:NTK_limit:bound_Dk+1inv} thanks to Step 3c, that $\absnormal{c^i_{k+1}}\leq C_c$ and $1+\N{w^i_{k+1}}_2^2\leq3\omega_i$ as of Step 3a (the factor $3$ is absorbed into $ c_1=6\sqrt{\vol{U}}\,C$),
    and that $\Nbig{\zeta^N_{k+1}}_{L^2(U)}\leq\Nbig{h-\CL g^N_{\theta_{k+1}}}_{L^2(U)}\leq C_R$ thanks to an analogous argument as the one that led to \eqref{eq:proof:thm:NTK_limit:bound_zeta0} together with the conclusion of Step 3b to obtain the last inequality.
    This closes the induction.
    
    Thus all properties in \eqref{eq:lem:bounds_z} hold on $\Omega_R$, and Step~1 gives $\bbP(\Omega_R)\geq1-\delta$, which concludes the proof.
\end{proof}

\subsection{Convergence of the PDE--NNTK integral operator \texorpdfstring{$T_{\CB^N_{\theta}}$}{} to its infinite-width limit \texorpdfstring{$T_{\CB^*_0}$}{}}
\label{sec:app:proof_NTK_limit}

\begin{lemma}
    \label{lem:differenceNNTK_finite_infinite}
    Let $\delta\in(0,1/2]$ and $K\in\bbN$ be a given number of training steps.
    Denote by $C_c,C_w,C_V,C_R,C_z,C_\gamma>0$ the constants defined in \eqref{eq:def_constants}.
    Assume that the number of neurons~$N$ is large enough such that \eqref{eq:thm:asm:Nlargeenough} is satisfied.
    Then it holds
    \begin{equation}
        \label{eq:lem:differenceNNTK}
        \bbE\left[
            \max_{0\leq k\leq K}
            \N{T_{\CB^N_{\theta_k}}-T_{\CB^*_0}}
            _{L^2(U)\rightarrow\CH^2(U)}
            \;\Big|\;\Omega_R
        \right]
        \leq
        \left(1+\frac{k_\CL C_B}{\gamma}\right)
        C\left(
            \frac{1+\alpha K}{\gamma^2N^{1-\beta}}
            +\frac{1}{\gamma N^{1/2}}
        \right),
    \end{equation}
    where $C_B=C_B(\sigma,\CD,\mu_0)>0$ is the constant of \eqref{eq:bound_TBstar} and $C=C(\sigma,\CD,\mu_0,\delta,C_c,C_R,C_z)>0$.
\end{lemma}

\begin{proof}
    On the event $\Omega_R$ of
    Lemma~\ref{prop:NTK_limit_bounds_z}, the bounds in
    \eqref{eq:lem:bounds_z} hold for all $i=1,\dots,N$
    and all $k=0,\dots,K$, together with
    \eqref{eq:proof:omega_transfer}.

    \emph{Step 0: Reduction to the PDE--NTK operators.}
    By the definitions of $T_{\CB^N_{\theta_k}}$ and
    $T_{\CB^*_0}$ in \eqref{eq:NNTK} and
    \eqref{eq:T_NNTK*}, respectively,
    \begin{equation}
        \label{eq:proof:lem:differenceNNTK_finite_infinite:1}
    \begin{split}
        T_{\CB^N_{\theta_k}}-T_{\CB^*_0}
        &= \frac{1}{\gamma}T_{B_{\theta_k}^N} \left(\Id_{L^2(U)} + \frac{1}{\gamma}\CL T_{B_{\theta_k}^N} \right)^{-1} - \frac{1}{\gamma}T_{B^*_0} \left(\Id_{L^2(U)} + \frac{1}{\gamma} \CL T_{B^*_0} \right)^{-1} \\
        &=
        \frac{1}{\gamma}\left(
            T_{B^N_{\theta_k}}
            -T_{B^*_0}
        \right)
        \left(
            \Id_{L^2(U)}+\frac{1}{\gamma}\CL T_{B^N_{\theta_k}}
        \right)^{-1}
        \\
        &\quad+
        \frac1\gamma T_{B^*_0}
        \left[
            \left(
                \Id_{L^2(U)}+\frac{1}{\gamma}\CL T_{B^N_{\theta_k}}
            \right)^{-1}
            -
            \left(
                \Id_{L^2(U)}+\frac1\gamma\CL T_{B^*_0}
            \right)^{-1}
        \right].
    \end{split}
    \end{equation}
    By \eqref{eq:lem:bounds_z}, $D^N_{\theta_k}$ is
    symmetric positive definite on $\Omega_R$.
    Hence, Lemma~\ref{lem:estimate_zeta} shows that
    both inverses exist and have
    $L^2(U)\rightarrow L^2(U)$ operator norm at most one.
    Using this together with the operator identity $A^{-1}-\widetilde{A}^{-1}=\widetilde{A}^{-1}(\widetilde{A}-A)A^{-1}$,
    which reads in our case
    \begin{equation}
    \begin{split}
        &
        \left(
            \Id_{L^2(U)}+\frac{1}{\gamma}\CL T_{B^N_{\theta_k}}
        \right)^{-1}
        -
        \left(
            \Id_{L^2(U)}
            +
            \frac{1}{\gamma}\CL T_{B^*_0}
        \right)^{-1}
        \\
        &\qquad=
        \left(
            \Id_{L^2(U)}
            +
            \frac{1}{\gamma}\CL T_{B^*_0}
        \right)^{-1}
        \frac{1}{\gamma}\CL
        \left(
            T_{B^*_0}
            -
            T_{B^N_{\theta_k}}
        \right)
        \left(
            \Id_{L^2(U)}+\frac{1}{\gamma}\CL T_{B^N_{\theta_k}}
        \right)^{-1},
    \end{split}
    \end{equation}
    we obtain with triangle inequality
    \begin{equation}
        \label{eq:proof:lem:differenceNNTK_finite_infinite:2}
    \begin{split}
        \Nbig{
            T_{\CB^N_{\theta_k}}-T_{\CB^*_0}
        }_{L^2(U)\rightarrow\CH^2(U)}
        &\leq
        \frac{1}{\gamma}\N{
            T_{B^N_{\theta_k}}
            -
            T_{B^*_0}
        }_{L^2(U)\rightarrow\CH^2(U)}
        +
        \frac{k_\CL C_B}{\gamma^2}
        \N{
            T_{B^N_{\theta_k}}
            -
            T_{B^*_0}
        }_{L^2(U)\rightarrow\CH^2(U)}
        \\
        &=
        \frac{1}{\gamma}\left(
            1+\frac{k_\CL C_B}{\gamma}
        \right)
        \N{
            T_{B^N_{\theta_k}}
            -
            T_{B^*_0}
        }_{L^2(U)\rightarrow\CH^2(U)},
    \end{split}
    \end{equation}
    where we used that $\N{T_{B^*_0}}_{L^2\rightarrow\CH^2}\leq C_B$ by Lemma~\ref{lem:standardNTKB:operatornorm}, specifically \eqref{eq:bound_TBstar}, and the bound $\N{\CL v}_{L^2(U)}\leq k_\CL\N{v}_{\CH^2(U)}$, which follows from Assumption~\ref{asm:CL} and linearity of $\CL$, as recorded in \eqref{eq:L_H2_L2} in Remark~\ref{rem:CL_linearity}.

    It therefore remains to estimate the difference between the PDE--NTK operators.
    Since integral operators depend linearly on their kernels,
    it holds that
    $T_{B^N_{\theta_k}}-T_{B^*_0}
    =T_{B^N_{\theta_k}-B^*_0}$.
    Moreover, the difference kernel
    $B^N_{\theta_k}-B^*_0$ and its weak derivatives
    with respect to $x$ up to order two belong to
    $L^2(U\times U)$ by the kernel estimates established
    in the proofs of
    Lemmas~\ref{lem:standardNTKB:operatornorm}
    and~\ref{lem:finitewidthNTKB:operatornorm}.
    Applying Lemma~\ref{lem:kernel_to_operator}, specifically \eqref{eq:kernel_to_operator}, to this difference kernel, and then using  triangle inequality after inserting mixed terms,
    gives
    \begin{equation}
        \label{eq:proof:thm:NTK_limit:discrepancy_10:aux:0}
    \begin{split}
        \N{
            T_{B^N_{\theta_k}}-T_{B^*_0}
        }_{L^2(U)\rightarrow\CH^2(U)}
        &\leq
        \sum_{\abs{\nu}\leq2}
        \N{
            D^x_\nu
            \left[B^N_{\theta_k}-B^*_0\right]
        }_{L^2(U\times U)}
        \\
        &\leq T_1+T_2+T_3,
    \end{split}
    \end{equation}
    where, using \eqref{eq:kernel_finite},
    \begin{align*}
        T_1
        &\coloneqq
        \sum_{\abs{\nu}\leq2}
        \N{
            \frac1N\sum_{i=1}^N
            \big[D_\nu(\chi s_{\theta_k^i})\big](\dummy)
            \left(
                \gamma(d^N_{\theta_k^i})^{-1}
                -\Id_{d+2}
            \right)
            \big(\big[\CL(\chi s_{\theta_k^i})\big](\dummy')\big)^\top
        }_{L^2(U\times U)},
        \\
        T_2
        &\coloneqq
        \sum_{\abs{\nu}\leq2}
        \N{
            \frac1N\sum_{i=1}^N
            \left(
                \big[D_\nu(\chi s_{\theta_k^i})\big](\dummy)
                \big(\big[\CL(\chi s_{\theta_k^i})\big](\dummy')\big)^\top
                -
                \big[D_\nu(\chi s_{\theta_0^i})\big](\dummy)
                \big(\big[\CL(\chi s_{\theta_0^i})\big](\dummy')\big)^\top
            \right)
        }_{L^2(U\times U)},
        \\
        T_3
        &\coloneqq
        \sum_{\abs{\nu}\leq2}
        \N{
            \frac1N\sum_{i=1}^N
            \big[D_\nu(\chi s_{\theta_0^i})\big](\dummy)
            \big(\big[\CL(\chi s_{\theta_0^i})\big](\dummy')\big)^\top
            -\big[D^x_\nu B^*_0\big](\dummy,\dummy')
        }_{L^2(U\times U)}.
    \end{align*}
    Here the first factor in each kernel is evaluated at $x$, and the factor containing $\CL$ is evaluated at $x'$.
    This is indicated by the notations $\dummy$ and $\dummy'$.

    \emph{Step 1: Regularization error $T_1$.}
    For the Gauss--Newton method, $d^N_{\theta_k^i}=\gamma\Id_{d+2}$ by \eqref{eq:def_d}, and hence $T_1=0$.
    For the Newton method, \eqref{eq:def_d} gives
    \begin{equation*}
    \begin{split}
        \gamma(d^N_{\theta_k^i})^{-1}-\Id_{d+2}
        &=
        \gamma(d^N_{\theta_k^i})^{-1}
        \left(
            \Id_{d+2}-\frac1\gamma d^N_{\theta_k^i}
        \right)
        \\
        &=
        \gamma(d^N_{\theta_k^i})^{-1}
        \left(
            \frac{1}{\gamma N^{1-\beta}}
        \int_U
        \left(h(x)-\CL g^N_{\theta_k}(x)\right)
        \big[\CL\big(\chi\,h_{\theta_k^i}\big)\big](x)\,dx
        \right).
    \end{split}
    \end{equation*}

    On $\Omega_R$, the bound $\abs{c_k^i}\leq C_c$ from
    \eqref{eq:lem:bounds_z} ensures that
    Lemma~\ref{lem:estimate_S_operatornorm_i} applies with
    $\theta=\theta_k$, $\theta^i=\theta_k^i$, and
    $\widetilde C=C_c$.
    Using the preceding identity, submultiplicativity of
    the operator norm,
    together with
    $\gamma\Nbig{(d^N_{\theta_k^i})^{-1}}_{\mathrm{op}}\leq2$
    from \eqref{eq:lem:bounds_z} and 
    Lemma~\ref{lem:estimate_S_operatornorm_i}, we obtain
    \begin{equation}
        \label{eq:proof:I_matrix}
    \begin{split}
        \N{
            \gamma(d^N_{\theta_k^i})^{-1}-\Id_{d+2}
        }_{\mathrm{op}}
        &\leq
        \frac{2C_S}{\gamma N^{1-\beta}}
        \left(1+\N{w_k^i}_2^2\right)
        \N{h-\CL g^N_{\theta_k}}_{L^2(U)}
        \\
        &\leq
        \frac{6C_SC_R}{\gamma}
        \frac{\omega_i}{N^{1-\beta}},
    \end{split}
    \end{equation}
    where the last inequality uses \eqref{eq:proof:omega_transfer} together with the residual bound in \eqref{eq:lem:bounds_z}.
    Here, $C_S$ denotes the constant from Lemma~\ref{lem:estimate_S_operatornorm_i} with $\widetilde C=C_c$.

    We can now estimate $T_1$ using the triangle inequality in the first step
    and Lemma~\ref{lem:rank_one_kernel}, specifically \eqref{eq:proof:rank_one}, in the second, applied with $\vartheta=\vartheta'=\theta_k^i$, $\widetilde C=C_c$, and $M=\gamma(d^N_{\theta_k^i})^{-1}-\Id_{d+2}$,
    to obtain
    \begin{equation}
        \label{eq:proof:thm:NTK_limit:discrepancy_10:aux:1}
    \begin{split}
        T_1
        &\leq
        \frac1N\sum_{i=1}^N
        \sum_{\abs{\nu}\leq2}
        \N{
            \big[D_\nu(\chi s_{\theta_k^i})\big](\dummy)
            \left(
                \gamma(d^N_{\theta_k^i})^{-1}-\Id_{d+2}
            \right)
            \big(\big[\CL(\chi s_{\theta_k^i})\big](\dummy')\big)^\top
        }_{L^2(U\times U)}
        \\
        &\leq
        \frac{C\vol{U}}{N}
        \sum_{i=1}^N
        \N{
            \gamma(d^N_{\theta_k^i})^{-1}-\Id_{d+2}
        }_{\mathrm{op}}
        \left(1+\N{w_k^i}_2^2\right)^2
        \\
        &\leq
        \frac{54C\vol{U}C_SC_R}{\gamma N^{1-\beta}}
        \frac1N\sum_{i=1}^N\omega_i^3
        \\
        &\leq
        \frac{C}{\gamma N^{1-\beta}}
        \frac1N\sum_{i=1}^N\omega_i^3,
    \end{split}
    \end{equation}
    where the third inequality uses \eqref{eq:proof:I_matrix} and $1+\N{w_k^i}_2^2\leq3\omega_i$ from \eqref{eq:proof:omega_transfer}.
    In the last inequality, we enlarged $C$ to absorb $54\vol{U}C_SC_R$.

    \emph{Step 2: Parameter-drift error $T_2$.}
    For every $i$, inserting a mixed term gives
    \begin{equation*}
    \begin{split}
        &
        \chi s_{\theta_k^i}
        \big(\CL(\chi s_{\theta_k^i})\big)^\top
        -
        \chi s_{\theta_0^i}
        \big(\CL(\chi s_{\theta_0^i})\big)^\top
        =
        \chi
        \big(s_{\theta_k^i}-s_{\theta_0^i}\big)
        \big(\CL(\chi s_{\theta_k^i})\big)^\top
        +
        \chi s_{\theta_0^i}
        \left(
            \CL\big(\chi(s_{\theta_k^i}-s_{\theta_0^i})\big)
        \right)^\top.
    \end{split}
    \end{equation*}
    On $\Omega_R$, the bounds
    $\abs{c_k^i}\leq C_c$ from \eqref{eq:lem:bounds_z}
    and $\abs{c_0^i}\leq C_{c,0}\leq C_c$ ensure that
    the difference estimate in Lemma~\ref{lem:rank_one_kernel}
    applies with $\widetilde C=C_c$.
    We can now estimate $T_2$ using the triangle inequality
    in the first step and Lemma~\ref{lem:rank_one_kernel},
    specifically \eqref{eq:proof:rank_one2}, in the second,
    applied with
    $\vartheta_1=\vartheta'_1=\theta_k^i$,
    $\vartheta_2=\vartheta'_2=\theta_0^i$, and
    $M=\Id_{d+2}$, to obtain
    \begin{equation}
        \label{eq:proof:thm:NTK_limit:discrepancy_10:aux:2}
    \begin{split}
        T_2
        &\leq
        \frac1N\sum_{i=1}^N
        \sum_{\abs{\nu}\leq2}
        \N{
            \big[D_\nu(\chi s_{\theta_k^i})\big](\dummy)
            \big(\big[\CL(\chi s_{\theta_k^i})\big](\dummy')\big)^\top
            -
            \big[D_\nu(\chi s_{\theta_0^i})\big](\dummy)
            \big(\big[\CL(\chi s_{\theta_0^i})\big](\dummy')\big)^\top
        }_{L^2(U\times U)}
        \\
        &\leq
        \frac{C\vol{U}}{N}\sum_{i=1}^N
        \left(1+\N{w_k^i}_2^2+\N{w_0^i}_2^2\right)
        \left(2+\N{w_k^i}_2^2+\N{w_0^i}_2^2\right)
        \N{\theta_k^i-\theta_0^i}_2
        \\
        &\leq
        \frac{16C\vol{U}}{N}\sum_{i=1}^N
        \omega_i^2\N{\theta_k^i-\theta_0^i}_2
        \leq
        \frac{16\alpha C\vol{U}}{N}\sum_{i=1}^N
        \omega_i^2\sum_{\ell=0}^{k-1}
        \N{\big(z_\ell^N\big)^i}_2
        \\
        &\leq
        \frac{16\alpha KC\vol{U}C_z}{\gamma N^{1-\beta}}
        \frac1N\sum_{i=1}^N\omega_i^3
        \\
        &\leq
        \frac{\alpha KC C_z}{\gamma N^{1-\beta}}
        \frac1N\sum_{i=1}^N\omega_i^3,
    \end{split}
    \end{equation}
    where the third inequality uses
    $1+\N{w_k^i}_2^2\leq3\omega_i$ from
    \eqref{eq:proof:omega_transfer} and
    $1+\N{w_0^i}_2^2=\omega_i$.
    The fourth inequality follows from the telescoping identity
    $\theta_k^i-\theta_0^i =\sum_{\ell=0}^{k-1}(\theta_{\ell+1}^i-\theta_\ell^i) =\alpha\sum_{\ell=0}^{k-1}(z_\ell^N)^i$
    and triangle inequality.
    The next-to-last inequality uses the bound $\N{(z_\ell^N)^i}_2 \leq C_z\omega_i/(\gamma N^{1-\beta})$ from \eqref{eq:lem:bounds_z}, together with $k\leq K$.
    In the last inequality, we enlarged $C$ to absorb $16\vol{U}$.

    \emph{Step 3: Finite-neuron error $T_3$.}
    For every multi-index $\abs{\nu}\leq2$, define the random function
    \begin{equation}
        \label{eq:def_Zalpha}
        Z_\nu^i
        \coloneqq
        \big[D_\nu(\chi s_{\theta_0^i})\big](\dummy)
        \left(\big[\CL(\chi s_{\theta_0^i})\big](\dummy')\right)^\top
        -
        \big[D^x_\nu B_0^*\big](\dummy,\dummy').
    \end{equation}
    By Assumption~\ref{asm:NN_mu0}, the parameters $\theta_0^1,\dots,\theta_0^N$ are i.i.d.
    Thus, for almost every fixed $(x,x')\in U\times U$, the scalar random variables $Z_\nu^1(x,x'),\dots,Z_\nu^N(x,x')$ are i.i.d.
    To verify that they are centered, note that
    $\abs{c_0^i}\leq C_{c,0}$, so
    Lemma~\ref{lem:bound_s} applies with
    $\widetilde C=C_{c,0}$ and an application of Cauchy--Schwarz inequality shows
    \begin{equation*}
    \begin{split}
        \abs{
            \big[D_\nu(\chi s_{\theta_0^i})\big](x)
            \left(\big[\CL(\chi s_{\theta_0^i})\big](x')\right)^\top
        }
        &\leq
        \N{\big[D_\nu(\chi s_{\theta_0^i})\big](x)}_2
        \N{\big[\CL(\chi s_{\theta_0^i})\big](x')}_2
        \leq C\omega_i^2,
    \end{split}
    \end{equation*}
    where the right-hand side is integrable since $\bbE[\omega_i^2] \leq 2\big(1+\bbE\N{w_0^i}_2^4\big) \leq 2\big(1+\CM_q^{4/q}\big)$.
    This bound holds for almost every $(x,x')$, with a constant independent of $(x,x')$, for every
    $\abs{\nu}\leq2$.

    By Fubini's theorem, for almost every fixed $x'$ the preceding domination holds almost surely, uniformly in $x$.
    For such $x'$, the integrand defining $B_0^*(x,x')$
    in \eqref{eq:NTK*} is twice continuously differentiable
    in $x$.
    The preceding integrable bounds on its derivatives
    allow differentiation under the expectation:
    the mean value theorem bounds the local difference
    quotients, and dominated convergence applies
    successively to the first and second derivatives.
    Thus,
    \begin{equation*}
    \begin{split}
        \big[D_\nu^x B_0^*\big](x,x')
        &=
        D_\nu^x
        \bbE\left[
            \chi(x)s_{\theta_0^i}(x)
            \left(\big[\CL(\chi s_{\theta_0^i})\big](x')\right)^\top
        \right]
        =
        \bbE\left[
            \big[D_\nu(\chi s_{\theta_0^i})\big](x)
            \left(\big[\CL(\chi s_{\theta_0^i})\big](x')\right)^\top
        \right],
    \end{split}
    \end{equation*}
    and consequently $\bbE[Z_\nu^i(x,x')]=0$.
    For the second moment of $Z_\nu^i(x,x')$,
    we use the variance identity
    $\bbE\abs{X-\bbE X}^2
    =\bbE\abs{X}^2-\abs{\bbE X}^2
    \leq\bbE\abs{X}^2$
    for square-integrable real random variables $X$
    to obtain
    \begin{equation}
        \label{eq:Z_moments}
    \begin{split}
        \bbE\left[\abs{Z_\nu^i(x,x')}^2\right]
        &\leq
        \bbE\left[
            \abs{
                \big[D_\nu(\chi s_{\theta_0^i})\big](x)
                \left(\big[\CL(\chi s_{\theta_0^i})\big](x')\right)^\top
            }^2
        \right]
        \\
        &\leq
        C\bbE[\omega_i^4]
        \leq8C\left(1+\CM_q^{8/q}\right)
        \leq C,
    \end{split}
    \end{equation}
    where the second inequality follows Cauchy--Schwarz inequality and uses Lemma~\ref{lem:bound_s}, applied with $\widetilde C=C_{c,0}$.
    The third inequality uses $q\geq8$, and we enlarged the generic constant $C=C(\sigma,\CD,\mu_0)>0$ as needed.
    
    Independence and centering imply $\bbE[Z_\nu^i(x,x')Z_\nu^j(x,x')]=0$ for $i\neq j$.
    Expanding the square, an application of Fubini--Tonelli's theorem therefore yields
    \begin{equation}
        \label{eq:proof:thm:NTK_limit:discrepancy_10:aux:variance}
    \begin{split}
        \bbE\N{
            \frac1N\sum_{i=1}^N Z_\nu^i
        }_{L^2(U\times U)}^2
        &=
        \int_{U\times U}
        \bbE\left[
            \abs{\frac1N\sum_{i=1}^N Z_\nu^i(x,x')}^2
        \right]\,dx\,dx'
        \\
        &=
        \frac1{N^2}\sum_{i=1}^N
        \int_{U\times U}
        \bbE\left[\abs{Z_\nu^i(x,x')}^2\right]\,dx\,dx'
        \\
        &=
        \frac1N
        \int_{U\times U}
        \bbE\left[\abs{Z_\nu^1(x,x')}^2\right]\,dx\,dx'
        \leq\frac{C}{N},
    \end{split}
    \end{equation}
    where the inequality in the last step follows from
    \eqref{eq:Z_moments}, after absorbing
    $\vol{U}^2$ into $C$.
    
    Finally, the definition of $T_3$ and Jensen's
    inequality give
    \begin{equation}
        \label{eq:proof:thm:NTK_limit:discrepancy_10:aux:3}
    \begin{split}
        \bbE[T_3]
        &=
        \sum_{\abs{\nu}\leq2}
        \bbE\N{
            \frac1N\sum_{i=1}^N Z_\nu^i
        }_{L^2(U\times U)}
        \leq
        \sum_{\abs{\nu}\leq2}
        \left(
            \bbE\N{
                \frac1N\sum_{i=1}^N Z_\nu^i
            }_{L^2(U\times U)}^2
        \right)^{1/2}
        \leq\frac{C}{N^{1/2}},
    \end{split}
    \end{equation}
    where we enlarged $C$ to absorb the finite number
    of multi-indices $\abs{\nu}\leq2$.
    
    \emph{Step 4: Conclusion.}
    Since $\bbP(\Omega_R)\geq1-\delta\geq1/2$,
    every non-negative random variable $X$ satisfies
    \begin{equation}
        \label{eq:conditioning}
        \bbE\left[X\mid\Omega_R\right]
        =
        \frac{\bbE\left[X\mathbf{1}_{\Omega_R}\right]}
             {\bbP(\Omega_R)}
        \leq
        \frac{\bbE X}{\bbP(\Omega_R)}
        \leq2\bbE X.
    \end{equation}
    
    We first note that, using \eqref{eq:conditioning}, the cubic weight average can be bounded as
    \begin{equation}
        \label{eq:proof:eq:proof:thm:NTK_limit:cubicweightaverage}
    \begin{split}
        \bbE\left[
            \frac1N\sum_{i=1}^N\omega_i^3
            \;\Big|\;\Omega_R
        \right]
        &\leq
        \frac2N\sum_{i=1}^N\bbE[\omega_i^3]
        =
        2\bbE[\omega_1^3]
        \leq
        8\left(1+\bbE\N{w_0^1}_2^6\right)
        \leq
        8\left(1+\CM_q^{6/q}\right),
    \end{split}
    \end{equation}
    where we used that $(1+a)^3\leq4(1+a^3)$ for $a\geq0$ and $q\geq8\geq6$ in the last step.

    The estimates for $T_1$ and $T_2$ from \eqref{eq:proof:thm:NTK_limit:discrepancy_10:aux:1} and \eqref{eq:proof:thm:NTK_limit:discrepancy_10:aux:2} hold simultaneously
    for all $k=0,\dots,K$ on $\Omega_R$, with right-hand
    sides independent of $k$.
    Moreover, $T_3$ depends only on the initialization and is therefore the same for every $k$.
    Thus \eqref{eq:proof:thm:NTK_limit:discrepancy_10:aux:0} together with \eqref{eq:proof:thm:NTK_limit:discrepancy_10:aux:1} and \eqref{eq:proof:thm:NTK_limit:discrepancy_10:aux:2}
    implies,
    on $\Omega_R$,
    \begin{equation}
        \label{eq:proof:eq:proof:thm:NTK_limit:discrepancy_10:aux:pre_final}
        \max_{0\leq k\leq K}
        \N{
            T_{B^N_{\theta_k}}-T_{B^*_0}
        }_{L^2(U)\rightarrow\CH^2(U)}
        \leq
        \frac{C(1+\alpha K)}{\gamma N^{1-\beta}}
        \frac1N\sum_{i=1}^N\omega_i^3
        +T_3.
    \end{equation}
    Taking conditional expectations in \eqref{eq:proof:eq:proof:thm:NTK_limit:discrepancy_10:aux:pre_final} and using \eqref{eq:proof:eq:proof:thm:NTK_limit:cubicweightaverage} for the first as well as  \eqref{eq:proof:thm:NTK_limit:discrepancy_10:aux:3} together with \eqref{eq:conditioning} for the second term,
    gives
    \begin{equation}
        \label{eq:proof:thm:NTK_limit:discrepancy_10:aux:final}
    \begin{split}
        &\bbE\left[
            \max_{0\leq k\leq K}
            \N{
                T_{B^N_{\theta_k}}-T_{B^*_0}
            }_{L^2(U)\rightarrow\CH^2(U)}
            \;\bigg|\;\Omega_R
        \right]
        \leq
        C\left(
            \frac{1+\alpha K}{\gamma N^{1-\beta}}
            +\frac1{N^{1/2}}
        \right),
    \end{split}
    \end{equation}
    for a constant $C=C(\sigma,\CD,\mu_0,\delta,C_c,C_R,C_z)>0$.

    Finally, taking the maximum over $k$ in \eqref{eq:proof:lem:differenceNNTK_finite_infinite:2} and then conditional expectations yields
    \begin{equation*}
    \begin{split}
        &\bbE\left[
            \max_{0\leq k\leq K}
            \N{
                T_{\CB^N_{\theta_k}}-T_{\CB^*_0}
            }_{L^2(U)\rightarrow\CH^2(U)}
            \;\Big|\;\Omega_R
        \right]
        \leq
        \left(1+\frac{k_\CL C_B}{\gamma}\right)
        C\left(
            \frac{1+\alpha K}{\gamma^2N^{1-\beta}}
            +\frac1{\gamma N^{1/2}}
        \right),
    \end{split}
    \end{equation*}
    which proves the statement.
\end{proof}

\subsection{Proof of Theorem~\ref{thm:NTK_limit}}
\label{sec:app:proof_thm_NTK_limit}

We now turn to the proofs of Theorem~\ref{thm:NTK_limit},
which establishes convergence to the infinite-width limit.
As an auxiliary ingredient,
let us first establish an integral representation of the Taylor remainder and bounds on its spatial derivatives up to order two.

\begin{lemma}[Bounds on remainder $R^N_k$]
    \label{lem:taylor_remainder}
    Let $\theta_k,z^N_k\in\bbR^{N(d+2)}$
    and let $\theta_{k+1}$ as in \eqref{eq:regularizedNewton}.
    Assume that $\abs{c^i_k},\abs{c^i_{k+1}}\leq\widetilde C$ for every $i=1,\dots,N$.
    Then, for every $x\in U$, it holds
    \begin{equation}
    \begin{split}
        g^N_{\theta_{k+1}}(x)
        &=
        g^N_{\theta_k}(x)
        +\alpha\,\chi(x)J^N_{\theta_k}(x)z^N_k
        +R^N_k(x)
    \end{split}
    \end{equation}
    with remainder
    \begin{equation}
    \label{eq:taylor_integral}  
    \begin{split}
        R^N_k(x)
        &=
        \alpha^2\int_0^1(1-t)\,
        (z^N_k)^\top
        \chi(x)H^N_{\theta_k+t\alpha z^N_k}(x)
        z^N_k\,dt.
    \end{split}
    \end{equation}
    The remainder satisfies
    $R^N_k\in\CC^2(\overline U)\cap\CH^2_{(0)}(U)$.
    For every spatial multi-index $\abs{\nu}\leq2$,
    \begin{equation}
        \label{eq:taylor_diff_under_integral}
        D_\nu R^N_k(x)
        =
        \alpha^2\int_0^1(1-t)\,
        (z^N_k)^\top
        \big[D_\nu
            \big(\chi H^N_{\theta_k+t\alpha z^N_k}\big)
        \big](x)
        z^N_k\,dt.
    \end{equation}
    
    In particular, under the hypotheses of
    Lemma~\ref{prop:NTK_limit_bounds_z}, it holds
    \begin{equation}
        \label{eq:taylor_bound_conditional}
        \bbE\left[
            \max_{0\leq k<K}
            \N{R^N_k}_{\CH^2(U)}
            \;\Big|\;\Omega_R
        \right]
        \leq
        \frac{C\alpha^2}{\gamma^2N^{1-\beta}}
    \end{equation}
    for a constant
    $C=C(\sigma,\CD,\mu_0,C_c,C_z)>0$.
\end{lemma}

\begin{proof}
    For fixed $x\in U$ define
    $\varphi(t)\coloneqq f^N_{\theta_k+t\alpha z^N_k}(x)$
    for $t\in[0,1]$.
    By \eqref{eq:NN} and
    Assumption~\ref{asm:NN_sigma'},
    $\varphi\in\CC^2([0,1];\bbR)$, with
    \begin{equation*}
    \begin{split}
        \varphi'(t)
        =
        \alpha J^N_{\theta_k+t\alpha z^N_k}(x)z^N_k
        \quad\text{ and }\quad
        \varphi''(t)
        =
        \alpha^2(z^N_k)^\top
        H^N_{\theta_k+t\alpha z^N_k}(x)z^N_k.
    \end{split}
    \end{equation*}
    A Taylor expansion with integral remainder gives
    \[
        \varphi(1)
        =
        \varphi(0)+\varphi'(0)
        +\int_0^1(1-t)\varphi''(t)\,dt.
    \]
    Multiplying by $\chi(x)$ proves
    \eqref{eq:taylor_integral}.

    For the spatial regularity, write $\widetilde\theta_t^i \coloneqq\theta_k^i+t\alpha(z^N_k)^i =(1-t)\theta_k^i+t\theta_{k+1}^i$,
    for which it holds $\absbig{(\widetilde\theta_t^i)_c} \leq (1-t)\abs{c_k^i}+t\abs{c_{k+1}^i} \leq\widetilde C$ for all $t\in[0,1]$.
    Hence Lemma~\ref{lem:bound_h} applies along each
    neuron-wise parameter segment.
    Its joint continuity statement, together with
    $H^N_{\widetilde\theta_t^i}=N^{-\beta}h_{\widetilde\theta_t^i}$
    for each single-neuron block, shows that
    $(t,x)\mapsto \big[D_\nu \big(\chi H^N_{\theta_k+t\alpha z^N_k}\big) \big](x)$
    is continuous on $[0,1]\times\overline U$
    for every $\abs{\nu}\leq2$.
    Since this set is compact, these functions are bounded.
    The parameter path is independent of $x$, so
    differentiation under the integral is justified by
    the dominated convergence theorem.
    This proves \eqref{eq:taylor_diff_under_integral}.
    The same continuity and domination show that the
    derivatives extend continuously to $\overline U$.
    Thus $R^N_k\in\CC^2(\overline U)\subset\CH^2(U)$.

    Moreover, since the integral representation extends to $\overline U$ and contains the factor $\chi(x)$,
    $R^N_k=0$ on $\partial U$.
    Since $\partial U$ is Lipschitz by Assumption~\ref{asm:D},
    the trace theorem gives $R^N_k\in\CH^1_0(U)$, and hence $R^N_k\in\CH^2_{(0)}(U)$.

    Since the Hessian and its spatial derivatives are
    block-diagonal, it holds
    \begin{equation*}
    \begin{split}
        D_\nu R^N_k(x)
        &=
        \alpha^2\int_0^1(1-t)
        \sum_{i=1}^N
        \big((z^N_k)^i\big)^\top
        \big[D_\nu\big(
            \chi H^N_{\widetilde\theta_t^i}
        \big)\big](x)
        (z^N_k)^i\,dt.
    \end{split}
    \end{equation*}
    Bounding each quadratic form by the operator norm
    yields
    \begin{equation}
        \label{eq:taylor_bound}
    \begin{split}
        \sup_{x\in U}\abs{D_\nu R^N_k(x)}
        &\leq
        \alpha^2\int_0^1(1-t)
        \sum_{i=1}^N
        \sup_{x\in U}
        \N{
            \big[D_\nu\big(\chi H^N_{\widetilde\theta_t^i}\big)\big](x)
        }_{\mathrm{op}}
        \N{(z^N_k)^i}_2^2\,dt
        \\
        &\leq
        \frac{\alpha^2}{2}
        \sum_{i=1}^N
        \sup_{t\in[0,1]}
        \sup_{x\in U}
        \N{
            \big[D_\nu\big(\chi H^N_{\widetilde\theta_t^i}\big)\big](x)
        }_{\mathrm{op}}
        \N{(z^N_k)^i}_2^2,
    \end{split}
    \end{equation}
    where we used $\int_0^1(1-t)\,dt=\frac12$.

    We now assume the hypotheses of
    Lemma~\ref{prop:NTK_limit_bounds_z} and work on
    the event $\Omega_R$.
    For every $k=0,\dots,K-1$, the bounds in
    \eqref{eq:lem:bounds_z} hold at both endpoints
    $\theta_k^i$ and $\theta_{k+1}^i$.
    By convexity, $\absbig{(\widetilde\theta_t^i)_c}\leq C_c$ and $\Nbig{\widetilde\theta_t^i-\theta_0^i}_2\leq1$ for all $t\in[0,1]$.
    Thus \eqref{eq:proof:omega_transfer} ensures
    $1+\Nbig{(\widetilde\theta_t^i)_w}_2^2\leq3\omega_i$.
    Lemma~\ref{lem:bound_h}, applied with
    $\widetilde C=C_c$, therefore implies
    \begin{equation*}
        \sup_{t\in[0,1]}
        \sup_{x\in U}
        \N{
            \big[D_\nu\big(\chi H^N_{\widetilde\theta_t^i}\big)\big](x)
        }_{\mathrm{op}}
        \leq
        \frac{3C\omega_i}{N^\beta}
    \end{equation*}
    for a constant $C=C(\sigma,\CD,C_c)>0$.
    Using this estimate in
    \eqref{eq:taylor_bound} together with the update bound \eqref{eq:lem:bounds_z} shows
    \begin{equation}
        \label{eq:proof:taylor_bound_weighted}
    \begin{split}
        \sup_{x\in U}\abs{D_\nu R^N_k(x)}
        &\leq
        \frac{\alpha^2}{2}
        \sum_{i=1}^N
        \frac{3C\omega_i}{N^\beta}
        \frac{C_z^2\omega_i^2}
             {\gamma^2N^{2-2\beta}}
        =
        \frac{3CC_z^2\alpha^2}
             {2\gamma^2N^{1-\beta}}
        \frac1N\sum_{i=1}^N\omega_i^3.
    \end{split}
    \end{equation}
    By the norm convention in
    Remark~\ref{rem:CL_linearity},
    we hence have
    \[
        \N{R^N_k}_{\CH^2(U)}
        =
        \sum_{\abs{\nu}\leq2}
        \N{D_\nu R^N_k}_{L^2(U)}
        \leq
        \sqrt{\vol{U}}
        \sum_{\abs{\nu}\leq2}
        \sup_{x\in U}\abs{D_\nu R^N_k(x)}.
    \]
    Summing \eqref{eq:proof:taylor_bound_weighted}
    over the finitely many multi-indices and taking
    the maximum over $k=0,\dots,K-1$ proves
    \begin{equation}
        \label{eq:taylor_bound_trajectory}
        \max_{0\leq k<K}
        \N{R^N_k}_{\CH^2(U)}
        \leq
        \frac{C\alpha^2}{\gamma^2N^{1-\beta}}
        \frac1N\sum_{i=1}^N\omega_i^3
    \end{equation}
    after absorbing the additional constants into $C$.

    Finally, $\bbP(\Omega_R)\geq1-\delta\geq1/2$,
    so \eqref{eq:conditioning} gives
    \begin{equation*}
    \begin{split}
        \bbE\left[
            \frac1N\sum_{i=1}^N\omega_i^3
            \;\Big|\;\Omega_R
        \right]
        &\leq
        \frac2N\sum_{i=1}^N\bbE[\omega_i^3]
        =
        2\bbE[\omega_1^3]
        \leq
        8\left(1+\bbE\N{w_0^1}_2^6\right)
        \leq
        8\left(1+\CM_q^{6/q}\right),
    \end{split}
    \end{equation*}
    where we used $(1+a)^3\leq4(1+a^3)$ for $a\geq0$ and $q\geq8$.
    Taking conditional expectations in \eqref{eq:taylor_bound_trajectory} and increasing $C$ to absorb this moment bound proves \eqref{eq:taylor_bound_conditional}.
\end{proof}

\begin{proof}[Proof of Theorem~\ref{thm:NTK_limit}]
    Let $\Omega_R$ denote the event of Lemma~\ref{prop:NTK_limit_bounds_z},
    on which \eqref{eq:lem:bounds_z} holds and $\bbP(\Omega_R)\geq1-\delta$.
    All path-wise estimates below are established on $\Omega_R$.

    \emph{Step 1: Recursion for approximation error in $\CH^2(U)$.}
    By Lemma~\ref{lem:eta_H10}, $g^N_{\theta_k}\in\CH^2_{(0)}(U)$ for all $k=0,\dots,K$.
    Moreover, Lemma~\ref{lem:standardNTKB:operatornorm} and the limit recursion \eqref{eq:NNfunction_evolution_limit}, starting from $g_0^*\equiv0$, imply inductively that $g_k^*\in\CH^2_{(0)}(U)$  for all $k=0,\dots,K$.
    Hence their difference belongs to $\CH^2_{(0)}(U)$,
    and Remark~\ref{rem:CL_linearity}, specifically \eqref{eq:L_H2_L2}, ensures that its image under
    $\CL$ belongs to $L^2(U)$.

    On $\Omega_R$, the hypotheses of
    Lemma~\ref{lem:finitewidthNTKB:operatornorm}
    hold at every $\theta_k$, by
    \eqref{eq:lem:bounds_z} and
    \eqref{eq:proof:omega_transfer}.
    Thus both $T_{\CB^N_{\theta_k}}$ and
    $T_{\CB^*_0}$ map $L^2(U)$ into
    $\CH^2_{(0)}(U)$.
    Finally, Lemma~\ref{lem:taylor_remainder}
    applies with $\widetilde C=C_c$ and gives
    $R^N_k\in\CH^2_{(0)}(U)$ for $k=0,\dots,K-1$.

    Subtracting \eqref{eq:NNfunction_evolution_limit} from \eqref{eq:NNfunction_evolution_prelimit}
    and inserting a mixed term gives the identity
    \begin{equation}
        \label{eq:proof:thm:NTK_limit:recursion}
    \begin{split}
        g^N_{\theta_{k+1}}-g^*_{k+1}
        &=
        \big(g^N_{\theta_k}-g^*_k\big)
        +\alpha
        \big(T_{\CB^N_{\theta_k}}-T_{\CB^*_0}\big)
        \big(h-\CL g^N_{\theta_k}\big)
        -\alpha T_{\CB^*_0}\CL
        \big(g^N_{\theta_k}-g^*_k\big)
        +R^N_k,
    \end{split}
    \end{equation}
    where all terms are in $\CH^2(U)$.
    For notational simplicity,
    denote
    $X_k\coloneqq \Nbig{g^N_{\theta_k}-g^*_k}_{\CH^2(U)}$.
    Taking norms in
    \eqref{eq:proof:thm:NTK_limit:recursion},
    using \eqref{eq:L_H2_L2},
    \eqref{eq:bound_TBstar}, and the residual bound
    in \eqref{eq:lem:bounds_z}, yields
    \begin{equation}
        \label{eq:proof:thm:NTK_limit:discrepancy_10}
        X_{k+1}
        \leq
        \left(1+\frac{\alpha k_\CL C_B}{\gamma}\right)X_k
        +\alpha C_R
        \Nbig{
            T_{\CB^N_{\theta_k}}-T_{\CB^*_0}
        }_{L^2(U)\rightarrow\CH^2(U)}
        +\N{R^N_k}_{\CH^2(U)}
    \end{equation}
    for all $k=0,\dots,K-1$.
    Since the coefficient of $X_k$ in
    \eqref{eq:proof:thm:NTK_limit:discrepancy_10}
    is at least one and the remaining terms
    are non-negative, that inequality implies
    \begin{equation}
        \label{eq:proof:thm:NTK_limit:discrepancy_10_2}
    \begin{split}
        \max_{0\leq j\leq k+1}X_j
        &\leq
        \left(1+\frac{\alpha k_\CL C_B}{\gamma}\right)
        \max_{0\leq j\leq k}X_j
        +\alpha C_R
        \Nbig{
            T_{\CB^N_{\theta_k}}-T_{\CB^*_0}
        }_{L^2(U)\rightarrow\CH^2(U)}
        +\N{R^N_k}_{\CH^2(U)}
    \end{split}
    \end{equation}
    for all $k=0,\dots,K-1$.
    Further denote $Y_k\coloneqq \bbE\big[ \max_{0\leq j\leq k}X_j \;\big|\;\Omega_R \big]$
    and define $\kappa_N\coloneqq \frac{1}{N^{\beta-1/2}} + \frac{1}{N^{1-\beta}}$.
    After taking the conditional expectation in \eqref{eq:proof:thm:NTK_limit:discrepancy_10_2} and employing Lemmas~\ref{lem:differenceNNTK_finite_infinite} and \ref{lem:taylor_remainder}, specifically \eqref{eq:lem:differenceNNTK} and \eqref{eq:taylor_bound_conditional},
    we obtain
    \begin{equation}
        \label{eq:proof:conditional_recursion}
        Y_{k+1}
        \leq
        \left(1+\frac{\alpha k_\CL C_B}{\gamma}\right)Y_k
        +C\kappa_N
    \end{equation}
    for all $k=0,\dots,K-1$ and for a constant $C=C(\alpha,\gamma,K,\delta,\beta,q,
    \sigma,\CD,\mu_0,C_{R,0})>0$.
    Here we used
    $N^{-1/2}\leq N^{-(1-\beta)}$
    for $N\geq1$, since $\beta>1/2$,
    and absorbed all $N$-independent factors into $C$.

    \emph{Step 2: Discrete Gr\"onwall's inequality.}
    Since $g_0^*\equiv0$, Lemma~\ref{lem:estimate_NNfunction_init},
    together with \eqref{eq:conditioning} and
    Cauchy--Schwarz inequality, gives
    \begin{equation}
        \label{eq:proof:e0}
    \begin{split}
        Y_0
        &=
        \bbE\left[
            \N{g^N_{\theta_0}}_{\CH^2(U)}
            \;\Big|\;\Omega_R
        \right]
        \leq
        2\left(
            \bbE\N{g^N_{\theta_0}}_{\CH^2(U)}^2
        \right)^{1/2}
        \leq
        \frac{2\sqrt{C}}
             {N^{\beta-1/2}}
    \end{split}
    \end{equation}
    for a constant $C=C(\sigma,\CD,\mu_0)>0$.
    An application of the discrete Gr\"onwall's inequality shows
    \begin{equation}
        \label{eq:proof:gronwall}
    \begin{split}
        Y_k
        &\leq
        \left(1+\frac{\alpha k_\CL C_B}{\gamma}\right)^kY_0
        +C\kappa_N
        \sum_{j=0}^{k-1}
        \left(1+\frac{\alpha k_\CL C_B}{\gamma}\right)^j
        \\
        &\leq
        \left(1+\frac{\alpha k_\CL C_B}{\gamma}\right)^K
        \big(Y_0+KC\kappa_N\big)
        \leq C\kappa_N
    \end{split}
    \end{equation}
    for every $k=0,\dots,K$.
    In the last inequality, we used \eqref{eq:proof:e0}
    and increased $C$ to absorb the remaining
    $N$-independent factors.

    \emph{Step 3: High-probability estimate.}
    Let $\Omega_N$ denote the event on which
    \[
        \max_{k=0,\dots,K}
        \Nbig{g^N_{\theta_k}-g_k^*}_{\CH^2(U)}
        >
        \frac{C}{\delta}\kappa_N,
    \]
    where
    $C=C(\alpha,\gamma,K,\delta,\beta,q,
    \sigma,\CD,\mu_0,C_{R,0})>0$
    is the constant from \eqref{eq:proof:gronwall}.
    Since all required updates are uniquely defined
    on $\Omega_R$,
    conditional Markov's inequality gives
    \begin{equation}
    \begin{split}
        \bbP(\Omega_N)
        &=
        \bbP(\Omega_N\cap\Omega_R)
        +\bbP(\Omega_N\cap\Omega_R^c)
        \\
        &\leq
        \bbP(\Omega_N\mid\Omega_R)\bbP(\Omega_R)
        +\bbP(\Omega_R^c)
        \\
        &\leq
        \frac{
            \bbE\left[
                \max_{k=0,\dots,K}X_k
                \;\big|\;\Omega_R
            \right]
        }{(C/\delta)\kappa_N}
        +\delta
        =
        \frac{Y_K}{(C/\delta)\kappa_N}+\delta
        \\
        &\leq2\delta.
    \end{split}
    \end{equation}
    Here we used
    $\bbP(\Omega_R^c)\leq\delta$
    in the third line and
    \eqref{eq:proof:gronwall} with $k=K$
    in the last line.
    Consequently,
    \eqref{eq:thm:NTK_limit} holds.
    Since $\bbP(\Omega_N^c)\geq1-2\delta$,
    this proves the asserted high-probability estimate.

    Since, for every fixed $\delta\in(0,1/2]$, the width
    conditions hold for all sufficiently large $N$ and
    $\kappa_N\rightarrow0$, the arbitrariness of $\delta$
    implies convergence in probability, uniformly over
    $k=0,\dots,K$.
\end{proof}

\section{Proof Details for Section~\ref{sec:globalconvergence}}
\label{sec:proof:globalconvergence}

In this section, we prove the spectral properties and convergence results stated in Section~\ref{sec:globalconvergence}.
We first establish the properties of $\CL T_{B_0^*}$ and then prove global convergence and the geometric rate for finite spectral projections.

\subsection{Proof of Lemma~\ref{lem:spectralpropertiesTB}}
\label{subsec:proof:spectralpropertiesTB}

We first verify the spectral properties of $\CL T_{B_0^*}$ and the orthogonality of its kernel to $\CL\big(\CH^2_{(0)}(U)\big)$.

\begin{proof}[Proof of Lemma~\ref{lem:spectralpropertiesTB}]
    Recalling the definition of the limit PDE--NTK integral operator (for gradient descent) $T_{B_0^*}$ from \eqref{eq:T_NTK*}, it is straightforward to verify that $\CL T_{B_0^*}$ is Hilbert--Schmidt, self-adjoint, and positive semi-definite, see \cite[Lemma~4.4]{jiang2023global} for a detailed proof.
    Hence, the first part of the statement follows.
    Since $\CL T_{B_0^*}$ is self-adjoint and compact, the spectral theorem applies, guaranteeing the existence of an orthonormal basis of $L^2(U)$ consisting of eigenvectors of $\CL T_{B_0^*}$.
    With $\CL T_{B_0^*}$ being positive semi-definite, its eigenvalues must be non-negative, verifying the second part of the statement.
    The positive eigenvalues can be ordered non-increasingly, counted with multiplicity.
    We refer to \cite[Lemma~4.5]{jiang2023global} for a detailed proof.
    Note, however, that, as made precise here, either $I$ or $I_0$ may be finite or empty, while at least one of them must be infinite, since $L^2(U)$ is infinite-dimensional.
    The last statement follows from \cite[Lemma~4.6]{jiang2023global}.
\end{proof}

\subsection{Proof of Theorem~\ref{thm:convergence}  and Corollary~\ref{cor:solution_convergence}}
\label{subsec:proof:globalconvergence}

We use the spectral decomposition of the residual evolution operator $\CL T_{\CB_0^*}$ to describe the evolution of the PDE residual.
Its orthogonality to the kernel and the dominated convergence theorem then yield convergence of the residual to zero.

\begin{proof}[Proof of Theorem~\ref{thm:convergence}]
    Let $r_k^* = h - \CL g_{k}^*\in L^2(U)$ denote the PDE residual at iteration $k$.
    Since $\CL g_{k+1}^* = \CL g_{k}^* + \alpha \CL T_{\CB_{0}^*} \left(h - \CL g_{k}^*\right)$ by \eqref{eq:NNfunction_evolution_limit},
    we obtain the recursion
    \begin{equation}
        \label{eq:proof:thm:convergence:10}
        \begin{split}
            r_{k+1}^*
            &= h - \CL g_{k+1}^*
            = \left(\Id_{L^2(U)} - \alpha \CL T_{\CB_{0}^*}\right) \left(h - \CL g_{k}^*\right) \\
            &= \left(\Id_{L^2(U)} - \alpha \CL T_{\CB_{0}^*}\right)r_{k}^*
        \end{split}
    \end{equation}
    and hence
    \begin{equation}
        \label{eq:proof:thm:convergence:12}
        r_{k}^* = \left(\Id_{L^2(U)} - \alpha \CL T_{\CB_{0}^*}\right)^{k}r_{0}^*.
    \end{equation}

    By Lemma~\ref{lem:spectralpropertiesTB}, the eigenfunctions of the operator $\CL T_{B_{0}^*}$ form an orthonormal basis of $L^2(U)$.
    We denote this basis by $\{e_\ell\}_{\ell\in I} \cup \{e^0_\ell\}_{\ell\in I_0}\subset L^2(U)$ with $I\subset\bbN$ and $I_0\subset\bbN$ denoting at most countable index sets such that
    \begin{equation}
        \label{eq:proof:thm:convergence:15}
        \CL T_{B_{0}^*} e_\ell = \lambda_\ell e_\ell
        \quad \text{ for all } \ell\in I,
        \quad\text{ and }\quad
        \CL T_{B_{0}^*} e_\ell^0 = 0
        \quad \text{ for all } \ell\in I_0,
    \end{equation}
    where $\{\lambda_\ell\}_{\ell\in I}$ is the family of strictly positive eigenvalues of the operator $\CL T_{B_{0}^*}$,
    counted with multiplicity and ordered in non-increasing order.

    Let $u\in\CH^2_{(0)}(U)$ be a solution of
    the PDE \eqref{eq:PDE}, which exists by
    Assumption~\ref{asm:solution}.
    Since $g_0^*\equiv0$ and $T_{\CB_0^*}$ maps
    $L^2(U)$ into $\CH^2_{(0)}(U)$ by
    Lemma~\ref{lem:standardNTKB:operatornorm},
    the recursion \eqref{eq:NNfunction_evolution_limit}
    implies inductively that $g_k^*\in\CH^2_{(0)}(U)$
    for every $k\in\bbN_0$.
    Consequently, $u-g_k^*\in\CH^2_{(0)}(U)$.
    Clearly, $h=\CL u$ and therefore
    \begin{equation}
        \label{eq:proof:thm:convergence:20}
        r_k^* = \CL(u-g_k^*) \in \operatorname{range}(\CL) \subset L^2(U)
    \end{equation}
    for all $k\in\bbN_0$.
    Since $\{e_\ell\}_{\ell\in I} \cup \{e^0_\ell\}_{\ell\in I_0}$ is an orthonormal basis of $L^2(U)$ according to Lemma~\ref{lem:spectralpropertiesTB},
    \begin{equation}
        \label{eq:proof:thm:convergence:25}
        r_k^* = \sum_{\ell\in I} c_\ell(k)e_\ell + \sum_{\ell\in I_0} c_\ell^0(k)e_\ell^0
    \end{equation}
    with coefficients $c_\ell(k) = \langle e_\ell,r_k^*\rangle_{L^2(U)}$ and $c_\ell^0(k) = \langle e_\ell^0,r_k^*\rangle_{L^2(U)}$.

    According to the last part of Lemma~\ref{lem:spectralpropertiesTB} (applied with $e=e_\ell^0\in L^2(U)$ satisfying $\CL T_{B_0^*}e_\ell^0=0$ and $f=u-g_k^*\in \CH_{(0)}^2$) it holds
    \begin{equation}
        \label{eq:proof:thm:convergence:30}
        \left\langle e_\ell^0,r_k^*\right\rangle_{L^2(U)} = 0
        \quad
        \text{ for all }k\in\bbN_0
        \text{ and all }\ell\in I_0
    \end{equation}
    and hence $c_\ell^0(k)=0$ for all $k\in\bbN_0$ and all $\ell\in I_0$.
    Consequently, \eqref{eq:proof:thm:convergence:25} reads
    \begin{equation}
        \label{eq:proof:thm:convergence:35}
        r_k^* = \sum_{\ell\in I}c_\ell(k)e_\ell.
    \end{equation}
    Note that if $I=\varnothing$, then \eqref{eq:proof:thm:convergence:35} immediately implies that $r_k^*=0$ and hence $L_k^*=0$ for all $k\in\bbN_0$.
    So in the following we may assume that $I\neq\varnothing$.

    By Remark~\ref{rem:spectralrelationNNTK}, the eigenfunctions in
    \eqref{eq:proof:thm:convergence:15} are also eigenfunctions of
    $\CL T_{\CB_0^*}$, with
    \begin{equation}
        \label{eq:proof:thm:convergence:40}
        \CL T_{\CB_{0}^*} e_\ell
        = \widetilde\lambda_\ell e_\ell
        = \frac{\lambda_\ell}
        {\gamma+\lambda_\ell}e_\ell
        \quad \text{ for all } \ell\in I,
        \quad\text{ and }\quad
        \CL T_{\CB_{0}^*} e_\ell^0 = 0
        \quad \text{ for all } \ell\in I_0.
    \end{equation}
    Since $\CL T_{\CB_0^*}$ is bounded, the powers of
    $\Id_{L^2(U)}-\alpha\CL T_{\CB_0^*}$ can be applied
    term-wise to the series \eqref{eq:proof:thm:convergence:35},
    which converges in $L^2(U)$.
    Using the expansion \eqref{eq:proof:thm:convergence:35} for $k=0$ in the representation \eqref{eq:proof:thm:convergence:12}
    together with
    \eqref{eq:proof:thm:convergence:40}, we obtain
    \begin{equation}
        \begin{split}
            r_k^*
            &= \left(\Id_{L^2(U)}-\alpha\CL T_{\CB_0^*}\right)^k
            \sum_{\ell\in I}c_\ell(0)e_\ell 
            = \sum_{\ell\in I}
            c_\ell(0)
            \left(1-\alpha\widetilde\lambda_\ell\right)^k
            e_\ell.
        \end{split}
    \end{equation}
    Comparing this expression with the expansion
    \eqref{eq:proof:thm:convergence:35} for some $k$ shows
    \begin{equation}
        c_\ell(k)
        =
        c_\ell(0)
        \left(1-\alpha\widetilde\lambda_\ell\right)^k
        \quad\text{ for all }\ell\in I.
    \end{equation}

    Since $\CL T_{B_0^*}$ is a compact, self-adjoint, and positive semi-definite operator as of Lemma~\ref{lem:spectralpropertiesTB}, its eigenvalues are bounded from above.
    Recalling the definition $\widetilde\lambda_{\max} =\|\CL T_{\CB_0^*}\|_{\mathrm{op}}$, we have $\widetilde\lambda_{\max} =\sup_{\ell\in I}\widetilde\lambda_\ell$, since $I\neq\varnothing$.
    As discussed in Remark~\ref{rem:spectralrelationNNTK}, $\widetilde\lambda_{\max} =\lambda_{\max}/(\gamma+\lambda_{\max})$.    Thanks to the assumption on the learning rate $\alpha$, it holds
    \begin{equation}
        0 < \alpha\widetilde\lambda_\ell
        \leq \alpha\widetilde\lambda_{\max}
        < 2
        \quad \text{ for all }\ell\in I,
    \end{equation}
    and hence
    \begin{equation}
        \left| 1-\alpha\widetilde\lambda_\ell \right|
        <1
        \quad \text{ for all }\ell\in I.
    \end{equation}
    Therefore, for every $\ell\in I$, we have
    \begin{equation}
        \begin{split}
            \lim_{k\to\infty} \left|c_\ell(k)\right|^2
            &= \lim_{k\to\infty} \left|c_\ell(0)\right|^2
            \left|1-\alpha\widetilde\lambda_\ell\right|^{2k}
            = 0.
        \end{split}
    \end{equation}
    Moreover, $\left|c_\ell(k)\right|^2 \leq \left|c_\ell(0)\right|^2$ for all $\ell\in I$ and all $k\in\bbN_0$, so the coefficients at $k=0$ provide a dominating sequence. Since $\sum_{\ell\in I}\left|c_\ell(0)\right|^2 = \left\|r_0^*\right\|_{L^2(U)}^2 < \infty$ by Parseval's identity, this dominating sequence is summable. The dominated convergence theorem therefore yields
    \begin{equation}
        \begin{split}
            \lim_{k\to\infty}\left\|r_k^*\right\|_{L^2(U)}^2
            &=
            \lim_{k\to\infty}
            \sum_{\ell\in I}
            \left|c_\ell(0)\right|^2
            \left|
                1-\alpha\widetilde\lambda_\ell
            \right|^{2k} \\
            &=
            \sum_{\ell\in I}
            \lim_{k\to\infty}
            \left|c_\ell(0)\right|^2
            \left|
                1-\alpha\widetilde\lambda_\ell
            \right|^{2k} 
            =
            0,
        \end{split}
    \end{equation}
    or equivalently
    \begin{equation}
        \lim_{k\to\infty}
        \N{h-\CL g_k^*}_{L^2(U)}
        =
        0,
    \end{equation}
    and hence, by the definition of the loss in \eqref{eq:loss_limit}, $\lim_{k\to\infty}L_k^* = 0$,
    which concludes the proof.
\end{proof}

We now use the additional stability estimate~\eqref{eq:PDE_stability} to establish convergence to the unique PDE solution in $\CH^1(U)$ (and in $\CH^2(U)$ under the corresponding stronger stability estimate), proving Corollary~\ref{cor:solution_convergence}.

\begin{proof}[Proof of Corollary~\ref{cor:solution_convergence}]
    Existence of a strong solution follows from
    Assumption~\ref{asm:solution}.
    If $u_1,u_2\in\CH^2_{(0)}(U)$ are two strong solutions,
    then $\CL(u_1-u_2)=0$.
    Applying \eqref{eq:PDE_stability} to $u_1-u_2$ gives
    \[
        \N{u_1-u_2}_{\CH^1(U)}
        \leq C
        \N{\CL(u_1-u_2)}_{L^2(U)}
        =0,
    \]
    proving uniqueness.

    Since $u-g_k^*\in\CH^2_{(0)}(U)$ and
    $\CL(u-g_k^*)=h-\CL g_k^*=r_k^*$,
    the stability estimate \eqref{eq:PDE_stability} yields
    \[
        \N{g_k^*-u}_{\CH^1(U)}
        \leq C\N{r_k^*}_{L^2(U)}.
    \]
    The right-hand side tends to zero by
    Theorem~\ref{thm:convergence}, completing the proof.

    If \eqref{eq:PDE_stability} holds with the $\CH^2(U)$ norm in
    place of the $\CH^1(U)$ norm, the same argument gives convergence
    in $\CH^2(U)$.
\end{proof}

\subsection{Proof of Theorem~\ref{lem:finite_mode_rate}}
\label{subsec:proof:finite_mode_rate}

We use the residual recursion to bound the contraction factors uniformly over the first $L$ spectral modes.
The spectral relation then gives the explicit rate for the full step $\alpha=1$.

\begin{proof}[Proof of Theorem~\ref{lem:finite_mode_rate}]
    For $\ell=1,\dots,L$, let $c_\ell(k)\coloneqq\langle r_k^*,e_\ell\rangle_{L^2(U)}$.
    By \eqref{eq:proof:thm:convergence:10}, the PDE residual satisfies
    \begin{equation}
        \label{eq:residual_recursion_limit}
        r^*_{k+1}=\big(\Id_{L^2(U)}-\alpha\,\CL T_{\CB^*_0}\big)r^*_k .
    \end{equation}
    Taking the inner product with $e_\ell$ and using that $\CL T_{\CB^*_0}$ is self-adjoint with $\CL T_{\CB^*_0}e_\ell=\widetilde\lambda_\ell e_\ell$, we obtain
    \begin{equation}
        \label{eq:coefficient_recursion}
        c_\ell(k+1)
        =\big(1-\alpha\widetilde\lambda_\ell\big)c_\ell(k),
        \qquad\text{ hence }\qquad
        c_\ell(k)=\big(1-\alpha\widetilde\lambda_\ell\big)^kc_\ell(0).
    \end{equation}

    By \eqref{eq:eigenvalue_ordering} and $\alpha>0$, the factors $1-\alpha\widetilde\lambda_\ell$ lie between $1-\alpha\widetilde\lambda_1$ and $1-\alpha\widetilde\lambda_L$.
    Their maximum absolute value is therefore attained at an endpoint:
    \begin{equation}
        \label{eq:max_at_endpoints}
        \max_{1\leq\ell\leq L}\absbig{1-\alpha\widetilde\lambda_\ell}
        =\max\Big\{\absbig{1-\alpha\widetilde\lambda_1},\ \absbig{1-\alpha\widetilde\lambda_L}\Big\}
        =q_L.
    \end{equation}
    Thanks to the assumption on $\alpha$, we have $0<\alpha\widetilde\lambda_\ell\leq\alpha\widetilde\lambda_{\max}<2$ for every $\ell=1,\dots,L$.
    Hence each factor has absolute value strictly less than one, and consequently $q_L<1$.

    Since $e_1,\dots,e_L$ are orthonormal, Parseval's identity together with \eqref{eq:coefficient_recursion} and \eqref{eq:max_at_endpoints} gives
    \begin{equation}
        \label{eq:parseval_rate}
        \Nbig{r^*_{k,L}}^2_{L^2(U)}
        =\sum_{\ell=1}^L\absbig{c_\ell(0)}^2\absbig{1-\alpha\widetilde\lambda_\ell}^{2k}
        \leq q_L^{2k}\sum_{\ell=1}^L\absbig{c_\ell(0)}^2
        =q_L^{2k}\Nbig{r^*_{0,L}}^2_{L^2(U)}.
    \end{equation}
    Taking square roots yields \eqref{eq:finite_mode_rate}.

    Finally, for $\alpha=1$, Remark~\ref{rem:spectralrelationNNTK} gives $\widetilde\lambda_\ell=\lambda_\ell/(\gamma+\lambda_\ell)\in(0,1)$.
    Thus all factors $1-\widetilde\lambda_\ell$ are positive, and
    \begin{equation}
        q_L
        =1-\widetilde\lambda_L
        =\frac{\gamma}{\gamma+\lambda_L}.
    \end{equation}
    Substituting this expression into \eqref{eq:finite_mode_rate} yields \eqref{eq:finite_mode_rate_fullstep}, which concludes the proof.
\end{proof}

% Acknowledgements and Disclosure of Funding should go at the end, before appendices and references
\section*{Acknowledgements}
\paragraph{Funding.}
This research project was supported by ``DMS-EPSRC: Asymptotic Analysis of Online Training Algorithms in Machine Learning: Recurrent, Graphical, and Deep Neural Networks'' (NSF DMS-2311500).
The authors would like to acknowledge the use of the University of Oxford Advanced Research Computing (ARC) facility in carrying out this work. 

\paragraph{Open access.}
For the purpose of Open Access, the authors have applied a CC BY public copyright license to any Author Accepted Manuscript (AAM) version arising from this submission.

\paragraph{AI Declaration.}
The authors used GPT-6 Astra, GPT-6.1 Sol and GPT-5.6 Sol, as well as Claude-5.5 Opus and Sonnet, to assist with literature search, routine calculations, mathematical proofs, proofreading the manuscript, and Python scripts for numerical simulations.
The authors independently derived and rigorously verified all mathematical statements, proofs, computations, and code.
We, the authors, take full responsibility for the accuracy, integrity, and originality of the final manuscript.

\addcontentsline{toc}{section}{References}
\bibliographystyle{abbrv}
\bibliography{biblio.bib}

\end{document}